\documentclass[reqno]{amsart}
\usepackage[margin=2cm,letterpaper]{geometry}

\usepackage{amssymb,enumerate,verbatim}
\usepackage{xr-hyper}
\newcommand{\cop}[1]{\begin{NoHyper}\cite[#1]{PavlovScholbach:Operads}\end{NoHyper}}
\newcommand{\csy}[1]{\begin{NoHyper}\cite[#1]{PavlovScholbach:Symmetry}\end{NoHyper}}
\usepackage[hypertex,backref]{hyperref}
\usepackage[all,cmtip]{xy}
\SilentMatrices % shut up, xy

\usepackage[utf8x]{inputenc} % for BibTeX from arXiv
\usepackage[T1]{fontenc} % ditto

\def\typeout#1{} % shut up, LaTeX

\makeatletter
\@addtoreset{equation}{section}
\@addtoreset{equation}{subsection}
\makeatother

\theoremstyle{definition}
\newtheorem{Defi}[equation]{Definition} \newcommand{\defi}{\begin{Defi}} \newcommand{\xdefi}{\end{Defi}} \newcommand{\refde}[1]{Definition~\ref{defi--#1}}
\newtheorem{Bsp}[equation]{Example} \newcommand{\exam}{\begin{Bsp}} \newcommand{\xexam}{\end{Bsp}} \newcommand{\refex}[1]{Example~\ref{exam--#1}}

\theoremstyle{remark}
\newtheorem{Bem}[equation]{Remark} \newcommand{\rema}{\begin{Bem}} \newcommand{\xrema}{\end{Bem}} \newcommand{\refre}[1]{Remark~\ref{rema--#1}}
\newtheorem{Nota}[equation]{Notation} \newcommand{\nota}{\begin{Nota}} \newcommand{\xnota}{\end{Nota}} 

\theoremstyle{plain}
\newtheorem{Theo}[equation]{Theorem} \newcommand{\theo}{\begin{Theo}} \newcommand{\xtheo}{\end{Theo}} \newcommand{\refth}[1]{Theorem~\ref{theo--#1}}
\newtheorem{Satz}[equation]{Proposition} \newcommand{\prop}{\begin{Satz}} \newcommand{\xprop}{\end{Satz}} \newcommand{\refpr}[1]{Proposition~\ref{prop--#1}}
\newtheorem{Lemm}[equation]{Lemma} \newcommand{\lemm}{\begin{Lemm}} \newcommand{\xlemm}{\end{Lemm}} \newcommand{\refle}[1]{Lemma~\ref{lemm--#1}}
\newtheorem{Coro}[equation]{Corollary} \newcommand{\coro}{\begin{Coro}} \newcommand{\xcoro}{\end{Coro}} \newcommand{\refco}[1]{Corollary~\ref{coro--#1}}

\newcommand{\refchap}[1]{\S\ref{chap--#1}}
\newcommand{\refsect}[1]{\S\ref{sect--#1}}
\newcommand{\refit}[1]{(\ref{item--#1})}
\newcommand{\refeq}[1]{(\ref{eqn--#1})}

\newcommand{\pf}{\begin{proof}} \newcommand{\xpf}{\end{proof}}

\let\ppar\endgraf % a paragraph inside a proof

\newcommand{\category}[1]{\mathbf{#1}}
\newcommand{\sSeq}{\category{sSeq}} % symmetric sequences
\newcommand{\sColl}{\category{sColl}} % symmetric collections
\newcommand{\Set}{\category{Set}} % sets
\newcommand{\simpl}{\category s} % simplicial ...
\newcommand{\sSet}{\simpl\category{Set}} % simplicial sets
\newcommand{\Top}{\category{Top}} % topological spaces
\newcommand{\PSh}{\category{PSh}} % presheaves
\newcommand{\sPSh}{\simpl\category{PSh}} % simplicial presheaves
\newcommand{\sMod}{\simpl\category{Mod}} % simplicial modules
\newcommand{\Sm}{\category{Sm}} % smooth varieties
\newcommand{\SmAn}{\category{SmAn}} % smooth complex manifolds
\newcommand{\Mod}{\category{Mod}} % modules
\newcommand{\sOper}{\category{sOper}} % symmetric operads
\newcommand{\Modsi}{\category{Mod}^\simpl} % modules with some simplicial aspect
\newcommand{\ModsAb}{\category{Mod}^{\simpl\Ab}} % modules with some simplicial abelian group aspect
\newcommand{\ModcCh}{\Mod^{\cCh}} % modules with connective chain complexes
\newcommand{\ModCh}{\Mod^\Ch} % modules with chain complexes
\newcommand{\Alg}{\mathbf{Alg}} % algebras
\newcommand{\Ch}{\category{Ch}} % chain complexes
\newcommand{\cCh}{\Ch_+} % connective chain complexes
\newcommand{\Fun}{\category{Fun}} % functors
\newcommand{\Ab}{\category{Ab}} % abelian groups

\newcommand{\proj}{\mathrm{pro}} % projective
\newcommand{\inje}{\mathrm{in}} % injective
\newcommand{\posi}{+} % positive
\newcommand{\ppos}{{[\posi]}} % positive or nonpositive
\newcommand{\stab}{\mathrm{s}} % stable
\newcommand{\CO}{\mathrm{CO}} % cofibrant objects
\newcommand{\cof}{\mathrm{C}} % cofibrations
\newcommand{\we}{\mathrm{W}} % weak equivalences
\newcommand{\AC}{\mathrm{AC}} % acyclic cofibrations
\newcommand{\cell}{\mathrm{cell}} % cellular closure
\def\mcr{{\rm Q}} % cofibrant replacement
\def\mfr{{\rm R}} % fibrant replacement
\def\ldf{{\rm L}} % left derived functor
\def\rdf{{\rm R}} % right derived functor
\newcommand{\id}{\mathrm{id}} % identity

\newcommand{\ws}[1]{\mathop{\rm cof}({#1})} % weak saturation
\newcommand{\inj}[1]{\mathop{\rm inj}({#1})} % maps with right lifting property
\newcommand{\colim}{\operatornamewithlimits{colim}} % colimit
\newcommand{\coeq}{\operatorname{coeq}} % coequalizer
\newcommand{\coker}{\operatorname{coker}} % cokernel
\newcommand{\cone}{\operatorname{cone}} % cone
\def\id{{\rm id}} % identity
\def\ev{{\rm ev}} % evaluation
\def\Ev{{\rm Ev}} % evaluation
\def\op{{\rm op}} % opposite category
\def\Mor{\mathop{\rm Mor}\nolimits} % the set of morphisms between two objects
\def\Hom{\mathop{\rm Hom}\nolimits} % internal hom
\def\Env{\mathop{\rm Env}\nolimits} % enveloping operad
\def\IHom{\Hom} % another macro for internal hom
\def\hHom{\mathop{\rdf \rm Hom}\nolimits}
\def\hMap{\mathop{\rdf \rm Map}\nolimits} % the derived mapping space
\def\Conv{\mathop{\rm Conv}\nolimits} % the convolution monoid
\def\hConv{\mathop{\rdf \rm Conv}\nolimits} % the derived convolution monoid
\def\Ho{\category{Ho}} % the set of morphisms in the homotopy category

\def\pp{\mathbin{\Box}} % pushout product binary operation
\mathchardef\ppdom"1400 % boxdot: domain of the pushout product
\newcommand{\ppop}{\mathop{\mathchoice{\doppop\Huge}{\doppop\Large}{\doppop\normalsize}{\doppop\small}}\displaylimits}
\newcommand{\doppop}[1]{\vcenter{#1\kern.2ex\hbox{$\Box$}\kern.2ex}}

\newcommand{\Sym}{{\rm Sym}} % symmetric powers

\def\Ax{\Sigma} % some automorphism group
\def\Z{{\bf Z}} % integers and the free simplicial abelian group
\def\NN{{\bf N}} % natural numbers
\def\Q{{\bf Q}} % rationals
\def\RR{{\bf R}} % real numbers
\def\CC{{\bf C}} % complex numbers
\def\A{{\bf A}} % affine space
\def\P{{\bf P}} % projective space
\def\Gm{\mathbf {G}_\mathrm m} % multiplicative group
\newcommand{\twi}[1]{\{#1\}} % twists

\def\H{{\rm H}} % (co)homology
\def\HH{{\bf H}} % hyper(co)homology
\def\HD{\H_{\rm D}} % Deligne cohomology
\def\N{{\rm N}} % normalized Moore chain complex
\def\BL{{\rm L}} % left Bousfield localization
\def\LL{{\rm L}} % left derived functor

\def\Spec{\mathop{\rm Spec}} % spectrum
\newcommand{\Comm}{\mathrm{Comm}} % the commutative operad
\newcommand{\sOpOp}{\mathrm{sOper}} % operad of symmetric operads
\def\Ei{{\rm E_\infty}} % the Barratt--Eccles operad
\def\Ai{{\rm A_\infty}} % A-infty operad
\def\E{{\rm E}} % the total space of the universal principal bundle
\newcommand{\an}{\mathrm{an}} % analytic space
\newcommand{\Y}{\mathcal{Y}} % some class of morphisms
\newcommand{\C}{\mathcal{C}} % the main category in which everything is constructed
\newcommand{\V}{\mathcal{V}} % enrichment
\newcommand{\D}{\mathcal{D}} % another category
\newcommand{\cA}{\mathcal{A}} % some abelian category
\newcommand{\cO}{O}

\let\x\times
\let\ol\overline
\renewcommand{\t}{\otimes}
\renewcommand{\r}{\rightarrow}
\newcommand{\lr}{\longrightarrow}

\def\matrix#1{\null\,\vcenter{\normalbaselines
    \ialign{\hfil$##$\hfil&&\quad\hfil$##$\hfil\crcr
      \mathstrut\crcr\noalign{\kern-\baselineskip}
      #1\crcr\mathstrut\crcr\noalign{\kern-\baselineskip}}}\,}
\def\vcd#1{\def\normalbaselines{\baselineskip20pt\lineskip1pt\lineskiplimit0pt } \harrowsize#1 \matrix}

\newdimen\harrowsize
\def\mapright#1{\smash{\mathop{\hbox to\harrowsize{\rightarrowfill}}\limits^{#1}}}

\def\mapdown#1{\Big\downarrow\rlap{$\vcenter{\hbox{$\scriptstyle#1$}}$}}

\catcode`@=11
\let\over\@@over
\let\atop\@@atop
\let\above\@@above
\let\overwithdelims\@@overwithdelims
\let\atopwithdelims\@@atopwithdelims
\let\abovewithdelims\@@abovewithdelims
\def\eqalign#1{\null\,\vcenter{\openup\jot\m@th
  \ialign{\strut\hfil$\displaystyle{##}$&$\displaystyle{{}##}$\hfil
      \crcr#1\crcr}}\,}
\newskip\xcentering
\def\eqalignno#1{\displ@y \tabskip\xcentering
  \halign to\displaywidth{\hfil$\@lign\displaystyle{##}$\tabskip\z@skip
    &$\@lign\displaystyle{{}##}$\hfil\tabskip\xcentering
    &\llap{$\@lign##$}\tabskip\z@skip\crcr
    #1\crcr}}
\def\cases#1{\left\{\,\vcenter{\normalbaselines\m@th
    \ialign{$##\hfil$&\quad##\hfil\crcr#1\crcr}}\right.}
\def\@writetocindents{}
\def\eqlabel#1{\refstepcounter{equation}\label{eqn--#1}\ifmmode\ifinner\else\eqno\fi\fi\hbox{\@eqnnum}} % Make equation labels work inside $$

\catcode`@=12

\begin{document}

\title{Symmetric operads in abstract symmetric spectra}

\centerline{\font\tfont=cmss17 \tfont\shorttitle}

\bigskip
{\tabskip0pt plus 1fil
\halign to\hsize{&\hfil#\hfil\cr
{\bf Dmitri Pavlov\/}\cr
Faculty of Mathematics, University of Regensburg\cr
Department of Mathematics and Statistics, Texas Tech University\cr
\href{https://dmitripavlov.org/}{https:/\negthinspace/dmitripavlov.org/}\vadjust{\medbreak}\cr
{\bf Jakob Scholbach\/}\cr
Mathematical Institute, University of M\"unster\cr
\href{https://wwwmath.uni-muenster.de/u/jakob.scholbach/}{https:/\negthinspace/wwwmath.uni-muenster.de/u/jscho\_04/}\cr
}}

\begin{abstract}
This paper sets up the foundations for derived algebraic geometry, Goerss--Hopkins obstruction theory, and the construction of commutative ring spectra
in the abstract setting of operadic algebras in symmetric spectra in an (essentially) arbitrary model category.

We show that one can do derived algebraic geometry a la To\"en--Vezzosi in an abstract category of spectra.
We also answer in the affirmative a question of Goerss and Hopkins by showing that the obstruction theory for operadic algebras in spectra
can be done in the generality of spectra in an (essentially) arbitrary model category.
We construct strictly commutative simplicial ring spectra representing a given cohomology theory
and illustrate this with a strictly commutative motivic ring spectrum representing higher order products on Deligne cohomology.

These results are obtained by first establishing Smith's stable positive model structure for abstract spectra
and then showing that this category of spectra possesses excellent model-theoretic properties:
we show that all colored symmetric operads in symmetric spectra valued in a symmetric monoidal model category are admissible, i.e., algebras over such operads carry a model structure.
This generalizes the known model structures on commutative ring spectra and $\Ei$-ring spectra in simplicial sets or motivic spaces.
We also show that any weak equivalence of operads in spectra gives rise to a Quillen equivalence of their categories of algebras.
For example, this extends the familiar strictification of $\Ei$-rings to commutative rings in a broad class of spectra, including motivic spectra.
We finally show that operadic algebras in Quillen equivalent categories of spectra are again Quillen equivalent.
%This paper is also available at arXiv:1410.5699v2.
\end{abstract}

\makeatletter\@setabstract\makeatother

\tableofcontents

\numberwithin{equation}{section}
\section{Introduction}

Ever since Brown's representability theorem, spectra occupy a central place in a variety of areas.
They are the objects representing cohomology theories,
i.e., for some cohomology theory $\H^*(-)$, one can find a spectrum~$E$ such that the cohomology of all spaces~$X$ is given by homotopy classes of morphisms of spectra
from the infinite suspension of~$X$ to a suspension of the spectrum:
$$\H^n(X) = [\Sigma^\infty X, \Sigma^n E].$$
Most cohomology theories in algebraic topology, algebraic geometry, and beyond carry a commutative and associative product
$$\H^m(X) \t \H^n(X) \r \H^{m+n}(X).$$
This makes it desirable to refine the multiplicative structure on the cohomology to one on the representing spectrum.
Ideally, one would like a strictly commutative and associative product
$$E \wedge E \r E$$
that gives back the above product.
In this case $E$ is called a commutative ring spectrum.
The following theorem is the basis of the homotopy theory of commutative ring spectra and spectra with a much more general multiplicative structure,
namely algebras over symmetric colored operads:

\theo (See \refth{operad.admissible.spectra}.)
\label{theo--introduction.admissibility}
Suppose $\C$ is a symmetric monoidal model category satisfying some mild additional assumptions (see \refde{nice.model.category} for the precise list),
$R$ is a commutative monoid in symmetric sequences in~$\C$,
$O$ is an operad in symmetric $R$-spectra (i.e., $R$-modules in symmetric sequences in~$\C$).
Then the stable positive model structure on $R$-spectra exists and gives rise to a model structure on $O$-algebras in $R$-spectra.
\xtheo

For example, this applies to $\C=\sSet_\bullet$, $R_n=(S^1)^{\wedge n}$ or $\C = \Top$.
If $O$ is the commutative operad (i.e., $O_n=S^0$),
$R$-modules are just symmetric spectra and $O$-algebras are \emph{commutative symmetric ring spectra}.
If $O$ is the Barratt--Eccles operad (i.e., $O_n=\E\Sigma_n$), then $O$-algebras are \emph{symmetric $\Ei$-ring spectra}.
(For the purposes of this introduction, we highlight algebras over these two operads, but all our results are true for arbitrary symmetric colored operads and their algebras.)

Another example is the category~$\C$ of pointed simplicial presheaves $\sPSh_\bullet(\Sm/S)$ or presheaves of complexes of abelian groups $\PSh(\Sm / S, \Ch)$ on the site of smooth varieties over a scheme~$S$,
equipped with the projective, flasque, or injective model structure, or any localization thereof
(such as the Nisnevich $\A^1$-localization), and $R_n=(\P^1)^{\wedge n}$.
In this case $R$-modules are known as motivic symmetric $\P^1$-spectra and commutative monoids are (strictly) \emph{commutative motivic symmetric ring spectra}.

In addition to the above existence result, we also give a supplementary condition that guarantees, for example, that the underlying spectrum of a cofibrant commutative ring spectrum is nonpositively cofibrant (see \refth{operad.strongly.admissible.spectra} for the precise statement).

In practice, it is often hard to construct \emph{strictly} commutative ring spectra.
Often it is the case that we instead can construct an algebra over an operad weakly equivalent to the commutative operad $\Comm$,
for example, the Barratt--Eccles operad~$\Ei$.
Essentially, this means that instead of defining a single product, there is a whole space of binary products and more generally $n$-ary products.
The following theorem says in particular that a multiplication whose space of $n$-ary operations is contractible, can be strictified to a strictly commutative and associative product.
Again, this is just one example, the result is true for algebras over arbitrary operads.

\theo \label{theo--Quillen.invariance.spectra} (See \refth{rectification.spectra}.)
With $\C$ and $R$ as above, any morphism $f\colon O\to P$ of operads in $R$-spectra
induces a Quillen adjunction between $O$-algebras and $P$-algebras,
which is a Quillen equivalence if $f$ is a weak equivalence.
\xtheo

\theo (See \refth{quasirect.spectra}.)
With $\C$ and $R$ as above, any simplicial operad~$O$
induces an equivalence of quasicategories between the underlying quasicategory of algebras over~$O$ in $R$-spectra
and the quasicategory of quasicategorical algebras over the operadic nerve of~$O$ in the underlying quasicategory of $R$-spectra.
\xtheo

We also study operadic algebras in spectra with values in Quillen equivalent categories (\refth{transport.algebras.spectra}).
As a special case we obtain the following Quillen invariance:

\theo (See \refco{transport.algebras.spectra}.)
For a weak equivalence $\varphi\colon R \stackrel \sim \r S$ of commutative monoids in $\Sigma \C$, and any levelwise fibrant operad $P$ in $S$-spectra and any levelwise cofibrant operad $O$ in $R$-spectra, there are Quillen equivalences
$$\eqalign{\varphi_* : \Alg_O^{\stab, \posi} (\Mod_R) & \leftrightarrows \Alg_{S \t_R O}^{\stab, \posi} (\Mod_S) : \varphi^*,\cr
\varphi_* : \Alg_{\varphi^* P}^{\stab, \posi} (\Mod_R) & \leftrightarrows \Alg_P^{\stab, \posi} (\Mod_S) : \varphi^*.\cr}$$
\xtheo

We give various applications of the above theorems.
The first can be paraphrased by saying that one can do derived algebraic geometry \`a la To\"en--Vezzosi \cite{ToenVezzosi:HomotopicalII} over ring spectra in the generality of $R$-spectra mentioned above:

\theo (See \refth{Toen.Vezzosi}.)
For any $\C$ and $R$ as in \refth{introduction.admissibility} such that $R_1$ is the suspension of some object $B \in \C$,
the category $\Mod_R^{\stab, \posi}$ of $R$-spectra with the positive stable model structure is a homotopical algebraic context.
\xtheo

Our next application is to show that one can do Goerss--Hopkins obstruction theory \cite{GoerssHopkins:Commutative, GoerssHopkins:Structured} in $R$-spectra.
This answers in the affirmative a question raised by Goerss and Hopkins.

\theo (See \refde{Goerss.Hopkins} and \refth{Goerss.Hopkins}.)
For any $\C$ and $R$ as in \refth{introduction.admissibility} such that $R_1$ is the suspension of some object $B \in \C$, the category $\Mod_R^{\stab, \posi}$ is a Goerss--Hopkins context.
\xtheo

Our third application is the construction of simplicial commutative ring spectra representing cohomologies given by commutative dgas.
For simplicity of the introduction, we just outline the statement here.

\theo (See \refth{construction.spectra}.)
Given a commutative ring spectrum $H$ with coefficients in $\Ch(\A)$ (for some abelian category~$\A$) there is a commutative ring spectrum with coefficients in~$\simpl\A$
representing the same cohomology as~$H$, with all higher products (such as Massey products) preserved.
\xtheo

The precise statement of \refth{construction.spectra} is flexible enough so that it is easily applicable in concrete cases.
For example, it allows to replace the monoid~$R$ with respect to which the stabilization of the model structure is done by a weakly equivalent one.
This is often necessary in practice and requires a Quillen equivalence of ring spectra valued in Quillen equivalent model categories.
Among the theorems above, this is by far the most involved fact in the homotopy theory of operadic algebras.
Moreover, the Dold--Kan equivalence uses the rectification theorem \ref{theo--rectification.spectra}.

As an example, we choose Deligne cohomology because of its interest in arithmetic and beyond.
Recall that Beilinson's conjecture \cite{Beilinson:Higher} expresses special $L$-values in terms of Deligne cohomology and algebraic $K$-theory.
Massey products on Deligne cohomology are related to special values of certain $L$-functions \cite{Deninger:Higher}.

\theo (See \refth{ultimate}.)
There is a strictly commutative motivic $\P^1$-spectrum representing Deligne cohomology with integral coefficients,
including the product structure and all higher product operations such as Massey products.
\xtheo

After a few recollections on general model categories in \refsect{preliminaries}, we study the stable model structure on $R$-spectra in \refsect{stable}.
Very briefly, our exposition uses a fairly flexible list of axioms on model structures in  symmetric sequences ($\Sigma \C$),
which is then transferred to the category of $R$-spectra, i.e., $R$-modules in $\Sigma\C$.
The next step is the usual stabilization using the technique of Bousfield localization.
We can summarize this by the following chain of Quillen adjunctions.
The middle adjunction is a Bousfield localization, while the other two adjunctions serve to transfer the model structure on the left to the right.
The superscripts indicate the precise choice of model structure: ``$\posi$'' and ``$\stab,\posi$'' refer to the positive and stable positive structures.
Underneath we indicate the place where the model structure in question is defined.
$$\begin{array}{ccccccc}
\Sigma^{\posi}\C&\leftrightarrows&\Mod_R^{\posi}&\leftrightarrows&\Mod_R^{\stab,\posi}&\leftrightarrows&\Alg_O(\Mod_R^{\stab,\posi}).\cr
\ref{defi--admissible}&&\ref{defi--Mod.R.unstable}&&\ref{theo--stable.R.general}&&\ref{theo--operad.admissible.spectra}\cr
\end{array}
\eqlabel{long.Quillen.chain}$$
For various special cases, the existence and properties of the stable positive model structure are known, but do not appear in the literature in the generality needed for our applications.
In \refth{operad.admissible.spectra}, we show that the stable positive model structure on $\Mod_R$ can be transferred to a model structure on algebras over arbitrary operads in $R$-spectra.
The key argument, due to Smith, is very simple: for a stable positive acyclic cofibration $f$ (whose level $0$ is by definition an \emph{isomorphism}), the $n$-fold pushout product $f^{\pp n}$ has very good properties.
For example, for any spectrum $X$ with a $\Sigma_n$-action, $X \t_{\Sigma_n} f^{\pp n}$ is a couniversal weak equivalence, since $\Sigma_n$ acts freely (roughly speaking) on $f^{\pp n}$ so that modding out the $\Sigma_n$-action on $X \t f^{\pp n}$ preserves the homotopical properties of $f$ and $X$.
Being a couniversal weak equivalence is weaker than being an acyclic cofibration, but joint with a mild condition on compact generators (called pretty small below) enough to establish the model structure on operadic algebras.
Along the way we prove the \emph{monoid axiom} for the stable model structures on $R$-spectra.

We go on to proving the operadic rectification result cited above (see \refth{rectification.spectra})
using the notion of symmetric flatness
which again holds for the stable \emph{positive} model structure on $R$-spectra.
The key point here is that for a weak equivalence $y\colon Y \r Y'$ of $\Sigma_n$-equivariant spectra and a positively cofibrant $R$-spectrum $X$, $y \t_{\Sigma_n} X^{\t n}$ is again a stable weak equivalence.

Sections \ref{sect--DAG}--\ref{sect--Deligne} are devoted to the above-mentioned applications.
An appendix explains how the (weak) assumptions on our basic model category $\C$ can be weakened further.

It is a pleasure to acknowledge the wealth of ideas that have helped to shape this paper.
For us, a starting point was an observation by Lurie that guarantees both the existence of a model structure on commutative monoids in a model category~$\C$
and a rectification result \cite[\S4.5.4]{Lurie:HA}.
It requires that $f^{\pp n}$ is a $\Sigma_n$-projective acyclic cofibration for all acyclic cofibrations $f\in\C$.
Roughly, this means that $\Sigma_n$ acts freely on the complement of the image of this iterated pushout product.
This is a harder condition than just asking that $f^{\pp n}/\Sigma_n$ is an acyclic cofibration.
In fact, Lurie's condition is rarely satisfied in practice.
It holds for chain complexes over a field of characteristic zero, but fails for the categories of simplicial sets.

The positive model structure on spectra is due to Smith.
It was studied in the context of topological spaces by Mandell, May, Schwede, and Shipley,
who showed the existence of model structures on commutative ring spectra and noted the rectification of $\Ei$-ring spectra in topological spaces
\cite[Theorem~15.1, Remark~0.14]{MandellMaySchwedeShipley:Model}.
The positive model structure on symmetric spectra with values in an arbitrary model category has been studied by Gorchinskiy and Guletski\u\i~\cite{GorchinskiyGuletskii:Positive}.
They showed the homotopy orbits property (under a strong assumption related to Lurie's condition mentioned above).
This property is a key step in the operadic rectification.
Harper also proved a rectification result as in \refth{Quillen.invariance.spectra} \cite[Theorem~1.4]{Harper:Symmetric} for $\C=\sSet_\bullet$,
which was generalized to $\C$ being the category of simplicial presheaves with the injective model structure by Hornbostel \cite[Theorem~3.6]{Hornbostel:Preorientations}.
These two model categories possess special features that substantially simplify the proof,
one of them being the fact that all objects are cofibrant.
Pereira has independently studied cofibrancy properties for operads in the (positive) stable model categories on spectra in $\sSet_\bullet$ with a particular focus on cofibrancy questions.
For example, his \cite[Theorem~1.3]{Pereira:Cofibrancy} is closely related to the case $\C=\sSet_\bullet$ of our \refth{operad.strongly.admissible.spectra}, see below.

In another direction, Harper showed the existence of a model structure on algebras over operads~\cite[Theorem~1.4]{Harper:Monoidal}
under the assumption that every symmetric sequence is projectively cofibrant.
Again, this is a strong assumption, which applies to such special categories as chain complexes over a field of characteristic zero.
In this case, rectification goes back to Hinich~\cite{Hinich:Homological}.
A recent application was the construction of motives (with rational coefficients) over general bases by Cisinski and D\'eglise \cite[Theorem~4.1.8]{CisinskiDeglise:Triangulated}.
In fact, our paper grew out from the desire to construct a convenient (i.e., fibrant) ring spectrum representing (higher) algebraic cobordism groups.
We plan to present such applications in a separate paper.

We thank Denis-Charles Cisinski, John Harper, Birgit Richter, and Brooke Shipley for helpful conversations.
This work was partially supported by the SFB 878 grant.

\numberwithin{equation}{subsection}
\section{Preliminaries on spectra}
\label{sect--model.structures.spectra}

\subsection{Model-categorical preliminaries}
\label{sect--preliminaries}

This paper uses the language of model categories.
Very briefly, we recall the less standard notions developed in \csy{\S\ref{chap--model} and~\S\ref{chap--symmetric}}.
A \emph{pretty small} model category $\C$ has, by definition, another model structure on the same underlying category that has the same weak equivalences, but fewer cofibrations, which are required to be generated by a set of maps whose domain and codomain is compact.
In this paper, a \emph{symmetric monoidal model category} is a model category that is also symmetric monoidal, such that the pushout product axiom \cite[Definition~4.2.6.1]{Hovey:Model} holds.
The unit axiom will always be explicitly mentioned separately.

According to Batanin and Berger \cite[Definition~1.7]{BataninBerger:Homotopy}, a symmetric monoidal model category $\C$ is called \emph{h-monoidal}, if for any object $Y \in \C$ and any (acyclic) cofibration $s$, $Y \t s$ is an (acyclic) h-cofibration.
If $\C$ is left proper, as are all the model categories in this paper, being an h-cofibration is equivalent to requiring that pushouts along this map are homotopy pushouts (and that it is moreover a weak equivalence in the acyclic case) \cite[Proposition~1.5]{BataninBerger:Homotopy}.
Moreover, any cofibration is an h-cofibration in this case.
The point of h-monoidality is that $Y \t s$ in many model categories~$\C$ is not an (acyclic) cofibration (for arbitrary~$Y$), but just an h-cofibration.
However, h-cofibrations are sufficiently well behaved for establishing model structures:
for example, combining h-monoidality with a compact generation condition yields the monoid axiom, see \cite[Proposition~2.5]{BataninBerger:Homotopy} or \csy{\refle{i.monoidal.monoid.axiom}}.

We call $\C$ \emph{flat} if for any cofibration $s$ and any weak equivalence $y \in \C$, the pushout product $s \pp y$ is a weak equivalence.
Again, at least for $s\colon \emptyset \r X$, a cofibrant object, this notion has appeared in various places in the literature.

To treat symmetric operads and their algebras, we also need symmetric enhancements of these concepts.
A class $S$ of morphisms in a symmetric monoidal model category $\C$ is \emph{(acyclic) symmetric h-monoidal}, if for any finite multi-index $n = (n_1, \ldots, n_e)$, $n_i \ge1$, and any object $Y \in \Sigma_n \C := (\prod_i \Sigma_{n_i}) \C$, and any finite family of maps $s = (s_i)$ in $S$, the map
$$Y \t_{\Sigma_n} s^{\pp n} := (Y \t s_1^{\pp n_1} \pp \cdots \pp s_e^{\pp n_e})_{\Sigma_n}$$
is an (acyclic) h-cofibration.
Here (for a single map $f$) $f^{\pp t}$ is the $t$-fold pushout product of $f$.
It agrees with the notation $Q^t_{t-1}$ in \cite[Definition~4.13]{Harper:Monoidal}, say.
A related condition is called \emph{(acyclic) $\Y$-symmetroidality} of~$S$: it requires that for any map $y$ in a fixed class of morphisms $\Y_n \subset \Mor \Sigma_n \C$ (for example all injective cofibrations), and any finite family $s$ of maps in $S$, the map
$$y \pp_{\Sigma_n} s^{\pp n}$$
is an (acyclic) cofibration.

Finally, a class $S$ of morphisms in $\C$ is called \emph{symmetric flat}, if for any weak equivalence $y \in \Sigma_n \C$ (i.e., a $\Sigma_n$-equivariant map that is a weak equivalence in $\C$) and any finite family of maps $(s_i)$ in $S$, the map
$$y \pp_{\Sigma_n} s^{\pp n} \eqlabel{symmetric.flat}$$
is a weak equivalence.

We call $\C$ symmetric flat, symmetric h-monoidal, or $\Y$-symmetroidal, if the cofibrations are symmetric h-monoidal symmetric flat, respectively the (acyclic) cofibrations are (acyclic) symmetric h-monoidal, respectively (acyclic) $\Y$-symmetroidal.

\exam
We explain what these notions mean in the special case of a single cofibration of the form $s\colon \emptyset \r X$ (where $X$ is a cofibrant object).
The case of an arbitrary cofibration is morally similar (and necessary to actually carry out the proofs).

Symmetric flatness asks that $(Y \r Y') \t_{\Sigma_n} X^{\t n}$ is a weak equivalence for any $\Sigma_n$-equivariant weak equivalence $Y \r Y'$.
This condition ensures a rectification result for operadic algebras, when $Y$ and $Y'$ are the $n$th levels of weakly equivalent operads.

The nonacyclic part of symmetric h-monoidality requires that $\emptyset \r Y \t_{\Sigma_n} X^{\t n}$ is an h-cofibration (and similarly in the acyclic case).
This condition arises when establishing the model structure on operadic algebras.

Finally, an important special case of symmetroidality is when $\Y_n$ consists just of $\emptyset \r1_\C$, and $s$ is a (single) map $\emptyset \r X$ to a cofibrant object $X$.
Then $\Y$-symmetroidality requires that $(X^{\t n})_{\Sigma_n}$ is cofibrant.
This special case gives a hint that symmetroidality ensures that cofibrant commutative algebras (for example, free algebras $\coprod_{n \ge0} X^{\t n}_{\Sigma_n}$ on a cofibrant object~$X$)
are cofibrant as objects of~$\C$.
A different class $\Y$ would correspond to similar cofibrancy properties for operads whose $n$th level is  constructed from the maps~$\Y_n$.
\xexam

\rema
\label{rema--explain.multiindices}
The three preceding definitions make use of multi-indices.
There are two reasons for that: one is that it is necessary for the homotopy theory of colored operads.
The more important reason, though, is the observation that multi-indices are necessary to ensure that the property in question
is stable under weak saturation (most importantly, under composition of morphisms).
We explain this for symmetric h-monoidality: suppose we know that in some model category~$\C$
maps of the form $Y \t_{\Sigma_n} s^{\pp n}$ are h-cofibrations for any (single) \emph{generating} cofibration $s$ and any object $Y \in \Sigma_n \C$
(and similarly for acyclic cofibrations).
This condition can be verified by hand in basic model categories such as simplicial sets, simplicial abelian groups, or rational chain complexes \csy{\refchap{examples}}.
Since only \emph{generating} cofibrations are concerned, this property is easily preserved under transfer \csy{\refth{transfer.symmetric.monoidal}} and,
with a bit more work for the acyclic part, under left Bousfield localizations \csy{\refth{symmetric.monoidal.localization}}.
In this sense, the condition for generating cofibrations is easy to establish for many model categories.
To obtain a similar statement for \emph{arbitrary} cofibrations, one has to consider (among other things) expressions of the form $Y \t_{\Sigma_n} (s_1 \circ s_2)^{\pp n}$,
where $s_1$, $s_2$ are generating cofibrations.
Unwinding this expression
leads to considering terms of the form $Y \t_{\Sigma_{n_1} \x \Sigma_{n_2}} s_1^{\pp n_1} \pp s_2^{\pp n_2}$ for $n=n_1+n_2$ (see for example \cite[Lemma~13]{GorchinskiyGuletskii:Symmetric}, \cite[Lemma~4.7]{Pereira:Cofibrancy}, \cite[Lemma A.1]{White:Model}, or \csy{\refle{combinatorial}}).
Thus, even if we are eventually only interested in algebras over single-colored operads,
where we need to know that $Y \t_{\Sigma_n} s^{\pp n}$ is an acyclic h-cofibration for a (single) arbitrary acyclic cofibration $s$ (see \refre{explain.symmetric.h.monoidal}),
we have to include multi-indices in the definition in order to ensure that symmetric h-monoidality only has to be tested on generating acyclic cofibrations,
which is crucial in showing that the property is stable under transfers and left Bousfield localizations.
Similar remarks apply for symmetric flatness and symmetroidality, see also \refre{explain.symmetric.flat}.
\xrema

While symmetric h-monoidality is satisfied for many model categories $\C$, symmetroidality and in particular symmetric flatness are more rare.
For example, simplicial sets are symmetric h-monoidal and symmetroidal, but not symmetric flat.
Simplicial presheaves with the projective model structure are symmetric h-monoidal, but not symmetroidal.
For a commutative ring $R$ the category of chain complexes of $R$-modules
is symmetric h-monoidal, symmetroidal, and symmetric flat precisely if $R$ contains~$\Q$, but none of these properties hold otherwise.
These and further basic examples are discussed in \csy{\refchap{examples}}.
A more sophisticated example is the
positive stable model structure on $R$-modules in symmetric sequences, i.e., \emph{symmetric $R$-spectra} with values in an abstract model category $\C$ (subject to some mild conditions).
This category is symmetric h-monoidal, symmetroidal, and symmetric flat, see \refpr{stable.R.positive} for the precise statement.

A \emph{monoidal left Bousfield localization} $\LL_S^\t \C$ of a symmetric monoidal model category $\C$ with respect to a class~$S$ is the left Bousfield localization in the bicategory of symmetric monoidal model categories.
Its underlying model category can be computed as $\LL_{S^\t} \C$, where $S^\t$ denotes the monoidal saturation of~$S$ in~$\C$, which can be computed as $S^\t = \{\mcr S \t X, X \in \C \text{ cofibrant}\}$.
If $\C$ is tractable, $S^\t$ can be taken to be $\mcr S \t \text{(co)dom} (I)$, where $I$ is some set of generating cofibrations with cofibrant source.
This notion is due to Barwick \cite[Proposition~4.47]{Barwick:Left}, the term monoidal localization was coined by White \cite{White:Monoidal}.

For a finite group $G$ and a subgroup $H$ and some object $X$ with a left $H$-action, we write $G \cdot_H X := \colim_H (\coprod_{G} X)$.
It carries a natural left $G$-action.

\subsection{The category of $R$-spectra}
\label{sect--symmetric.spectra}
\label{sect--symseq}

The first step in the construction of symmetric spectra in a category $\C$ is the category
$$\Sigma\C := \Fun(\Sigma, \C)$$
of \emph{symmetric sequences}, i.e., functors from the category $\Sigma$ of finite sets and bijections or, equivalently, its skeleton.
Prominent examples are given in \refex{monoids}.

There is an obvious adjunction
$$G_n : \Sigma_n\C\leftrightarrows \Sigma\C : \ev_n,\eqlabel{adjunction.Sigma.C}$$
where $\ev_n$ is the evaluation on $n$ and $G_n(X)(m)$ is $X$ for $m=n$ and the initial object of~$\C$ else.
For some fixed $k\ge0$, these assemble to an adjunction
$$G_{\ge k} : \prod_{n\ge k} \Sigma_n \C \leftrightarrows \Sigma \C : \ev.\eqlabel{adjunction.symmetric.sequences}$$
For $k=0$ this is an equivalence of categories, but we will mostly be interested in $k=1$ in the sequel.

The category $\Sigma \C$ is equipped with the monoidal structure, denoted~$\t$, coming from the disjoint union of finite sets \cite[Definition~2.1.3]{HoveyShipleySmith:Symmetric}.
It satisfies
$$G_n (X) \t G_{n'}(X') = G_{n+n'}(\Sigma_{n+n'} \cdot_{\Sigma_n \x \Sigma_{n'}} X \t X').\eqlabel{product.Sigma.C}$$

Let $R$ be a commutative monoid in~$\Sigma\C$.
We denote the category of $R$-modules in~$\Sigma\C$ by~$\Mod_R$ and refer to it as the \emph{category of $R$-spectra}.
See \cite[\S2.2]{HoveyShipleySmith:Symmetric} or \cite{Shipley:HZ-algebras} for more details, where this category is denoted by $\category{Sp}^\Sigma(\C, R)$.
$R$-spectra form a symmetric monoidal category with the tensor product of $R$-modules $M$~and~$N$ being
$$M\t_R N=\coeq(M\t R\t N\rightrightarrows M\t N),$$
where the tensor products on the right are computed in~$\Sigma\C$.

\subsection{Model structures on symmetric sequences}
\label{sect--model.structures.symseq}

To do homotopy theory, we use the language of model categories.
Throughout, we work with a basic category $\C$ satisfying the following assumptions.
In \refsect{appendix}, we will explain that these assumptions can be weakened even further.

\defi
\label{defi--nice.model.category}
A model category~$\C$ is {\it nice\/} if it is
pretty small \csy{\refde{pretty.small}},
symmetric monoidal, h-monoidal, flat \csy{\refde{flat}}, and tractable \cite[Definition~1.21]{Barwick:Left}.
\xdefi

\rema
\label{rema--combinatorial}
This list may look intimidating, but in practice these assumptions are both mild and robust.
We emphasize that no conditions such as \emph{symmetric} h-monoidality or \emph{symmetric} flatness are imposed on~$\C$.
Examples of nice model categories include simplicial sets and chain complexes of abelian groups, as well as presheaves with values in these categories \csy{\refchap{examples}}.
Moreover, if $\C$ is nice, then so is any monoidal left Bousfield localization $\LL_S^\t\C$, as well as any model structure that is transferred from~$\C$ to~$\D$,
provided that the adjunction has good monoidal properties, see \csy{\S\ref{chap--transfer} and \S\ref{chap--localization}} for precise statements.
In particular, the localization of the injective or the projective model structure (on simplicial presheaves or presheaves of cochains) with respect to some topology is also a nice model category.
Topological spaces are not nice in the sense above, but the majority of our results in this paper does apply for $\C = \Top$.
See \refsect{appendix}.
\xrema

Depending on the model category $\C$, there are typically different interesting model structures on $\Sigma \C$, which later give rise to model structures on $R$-spectra (Definitions \ref{defi--Mod.R.unstable}, \ref{defi--stable}).
The classical choice, for example used by Hovey \cite[Theorem~8.2]{Hovey:Spectra}, is the projective model structure on $\Sigma \C$.
Somehow at the other extreme, there is the injective model structure on $\Sigma \C$.
For $\C=\sSet$, the injective model structure on $\Sigma \C$ yields to a model structure on $S^1$-spectra called the level $S$-model structure \cite[Proposition~2.2]{Shipley:Convenient}.
Both injective and projective model structures their individual advantages and disadvantages:
the projective model structure is admissible, but not usually strongly admissible (\refle{projective.admissible}).
The injective model structure is (up to a minor technical point, see \refle{injective.admissible}) strongly admissible, but its fibrations are practically impossible to describe, which makes computations of derived mapping spaces very hard.
We therefore present a short list of axioms for model structures on $\Sigma \C$ that are used in the sequel.
As a guide for the reader, we point out that admissible model structures on $\Sigma \C$ will ultimately guarantee that for all operads $O$ in the stable positive model structure on symmetric $R$-spectra, there is a model structure on $O$-algebras (\refth{operad.admissible.spectra}), i.e., $O$ is admissible.
The strong admissibility of the model structure on $\Sigma \C$ will be the essential ingredient guaranteeing that an operadic algebra in spectra, which is cofibrant as an algebra, is also cofibrant as a spectrum (see \refth{operad.strongly.admissible.spectra}).

Recall the notation $G \cdot_H -$ from \refsect{preliminaries} and the notation $\Sigma_n$, $h^{\pp n}$ for a multi-index $n = (n_1, \ldots, n_k)$ from \csy{\refde{multi}}.
We write $\cof$ ($\AC$) for the class of (acyclic) cofibrations of a model category.

\defi \label{defi--admissible}
Let $k$ be either $0$ (referred to as nonpositive case below) or $1$ (positive case).
A collection of model structures on $\Sigma_m \C$, $m \ge k$ is called \emph{admissible} if the following conditions are satisfied:
\begin{enumerate}[(1)]
\item \label{item--admissible.tractable} Each $\Sigma_m \C$ is a tractable model category.

\item
\label{item--admissible.we} The weak equivalences do not depend on the group action, i.e., $\we_{\Sigma_m \C} = \varphi^{-1}(\we_\C)$, where $\varphi$ denotes the functor that forgets any action of a finite group on some object in $\C$.

\item \label{item--admissible.injective}
The cofibrations satisfy:
$$\Sigma_m \cdot \cof_\C \subset \cof_{\Sigma_m \C} \subset \varphi^{-1} (\cof_\C).$$

\item
\label{item--admissible.product.2}
For any decomposition $m = m'+m''$, $m',m''\ge k$, $f \in \cof_{\Sigma_m \C}$ and $f' \in \cof_{\Sigma_{m'} \C}$, the map
$$\Sigma_m \cdot_{\Sigma_{m'} \x \Sigma_{m''}} (f \pp f')\eqlabel{admissible.product.2}$$
is a cofibration in $\Sigma_m \C$.

\medskip\vbox{
The following property will eventually be used to show that some operad $O$ is strongly admissible, i.e., cofibrant $O$-algebras in spectra have cofibrant underlying objects (\refth{operad.strongly.admissible.spectra}).
We will only apply it in the positive case.
Let $\Y_{n}$ be a class of morphisms in $\Sigma_n \Sigma \C$,
where $n$ is a multi-index.
We suppose that for any $y \in \Y_{n}$, $y \pp - $ preserves (acyclic) cofibrations in $\Sigma_n^\inje \Sigma \C$, i.e., those $\Sigma_n$-equivariant maps that are (acyclic) cofibrations in $\Sigma \C$.
The purpose of the maps $\Y_n$ is that the $n$th level of an operad $O$ is constructed (by means of pushouts, transfinite compositions, and retracts) from the maps $\Y_n$.
We have two main examples in mind: $\Y_{n}$ consists just of the single map $\emptyset \r1_{\Sigma \C} (=G_0(1_\C))$ with the trivial $\Sigma_n$-action.
In this case strong $\Y$-admissibility will be relevant for the strong admissibility of the commutative operad.
Another example is the class $\Y_{n}$ of injective cofibrations in $\Sigma_n \Sigma \C$, i.e., $\Sigma_n$-equivariant maps that are cofibrations in $\Sigma \C$.
In this case strong $\Y$-admissibility will yield the strong admissibility of levelwise injectively cofibrant operads.
See \refth{operad.strongly.admissible.spectra} for the precise statement.

With these motivations, we call the model structure \emph{strongly $\Y$-admissible} if, in addition to \refit{admissible.tractable}--\refit{admissible.product.2}, the following condition holds.
}

\item
\label{item--admissible.wedge.2}
For any multi-index $n \ge1$, any multi-index (of the same size) $t \ge1$, any $y \in \Y_{n}$, any finite family of generating (acyclic) cofibrations $h \in \Sigma_t \C$ (i.e., $h_i \in \Sigma_{t_i} \C$), the expression
$$\Sigma_{u+tn} \cdot_{\Sigma_u \x (\Sigma_n \rtimes \Sigma_t^n)} y_u \pp h^{\pp n}$$
is an (acyclic) cofibration in $\Sigma_{tn+u} \C$, where $y_u \in \Sigma_n \Sigma_u \C$ is the $u$th component of $y$.
\end{enumerate}
\xdefi

\nota
The model structure on $\Sigma \C$ that is the product of these model structures is just denoted by $\Sigma \C$ in the nonpositive case and $\Sigma^\posi \C$ in the positive case.
Concretely a cofibration (respectively fibration, weak equivalence) in $\Sigma \C$ is a map $f=(f_m)_{m \ge0}$ whose levels $f_m$ lie in the corresponding class of $\Sigma_m \C$ for all $m \ge0$.
In $\Sigma^\posi \C$ a map $f=(f_m)_{m \ge0}$ is a cofibration (respectively fibration, weak equivalence) if the $f_m$ are in the corresponding class of $\Sigma_m \C$ for $m \ge1$ and if $f_0$ is an isomorphism (and arbitrary for fibrations and weak equivalences).
In order to indicate some particular choice of the model structures on $\Sigma_m \C$, we write $\Sigma^{\posi,\proj}$ etc.
The notation $\Sigma^\ppos \C$ refers to either $\Sigma \C$ or $\Sigma^\posi \C$.
\xnota

\rema
Since $\C$ is combinatorial by \ref{defi--admissible}\refit{admissible.tractable}, both the projective model structure $\Sigma \C^{\proj}$ and the injective model structures $\Sigma \C^{\inje}$ exist (see below).
Thus, the above condition \ref{defi--admissible}\refit{admissible.injective} is equivalent to requiring that the identity functors $\Sigma \C^{\proj} \to \Sigma \C \to \Sigma \C^{\inje}$ be left Quillen functors.
In order to extend our results to certain noncombinatorial model categories (\refsect{appendix}), the above definition was phrased without explicitly referring to these two model structures.
\xrema

We now give examples of (strongly) admissible model structures.
\refle{injective.admissible} shows that the injective model structure on $\Sigma \C$ is strongly admissible, except, possibly, for the quasitractability.
(This is also the reason why \ref{defi--admissible}\refit{admissible.product.2} and \refit{admissible.wedge.2} do not contain any requirements on the acyclic cofibrations).
The quasitractability requirement \ref{defi--admissible}\refit{admissible.tractable} is primarily of technical importance.
It will be used to carry through monoidal properties to the stabilization of $R$-modules,
which is helpful to prove the monoid axiom for the stable structure on $R$-modules (\refpr{symseq.general}).
Ignoring this necessity, the injective model structure can be used in the sequel.
However, fibrancy is very difficult to check in this model structure.
A strongly admissible structure with controlled cofibrations (and therefore, acyclic fibrations) is provided by \refpr{minimal.admissible}.

\lemm \label{lemm--injective.admissible}
Let $\C$ be a combinatorial, symmetric monoidal model category.
The injective model structure on $\Sigma \C$ exists and satisfies the following properties.
\begin{enumerate}[(1)]
\item
\label{item--injective.admissible.1}
For $m = m'+m''$, the functor $\Sigma_m \cdot_{\Sigma_{m'} \x \Sigma_{m''}} - \t - \colon \Sigma_{m'} \C \x \Sigma_{m''} \C \r \Sigma_m \C$ is a left Quillen bifunctor.
\item
\label{item--injective.admissible.2}
For any finite family $h = (h_1, \ldots, h_e)$ of (acyclic) cofibrations in $\Sigma_{t_i} \C$ with $t_i>0$, and any (injective) cofibration $y \in \Sigma_n \Sigma_u \C$, the following map is a cofibration
$$\Sigma_{tn+u} \cdot_{\Sigma_n \rtimes (\Sigma_u \x \Sigma_t^n)} y \pp h^{\pp n} \in \Sigma_{tn+u} \C.$$
\item
\label{item--injective.admissible.conclusion}
The injective model structure on $\Sigma \C$ is strongly $\Y$-admissible for the class $\Y_n$ of $\Sigma_n$-injective cofibrations in $\Sigma_n \Sigma \C$, if we replace tractability by combinatoriality.
\end{enumerate}
\xlemm

\pf
The existence is a special case of \cite[Proposition A.2.8.2]{Lurie:HTT}.
\refit{injective.admissible.1} holds because of the noncanonical, but functorial isomorphism $\varphi(\Sigma_m \cdot_{\Sigma_{m'} \x \Sigma_{m''}} -) \cong \Sigma_m / (\Sigma_{m'} \x \Sigma_{m''}) \cdot \varphi( -)$ \csy{\refle{dirty}}.

\refit{injective.admissible.2}~Using the notation of \ref{defi--admissible}\refit{admissible.wedge.2}, $h^{\pp n}$ is an (acyclic) cofibration in $\C$ by the pushout product axiom and therefore $y \pp h^{\pp n}$ is again a cofibration in $\C$ by the assumption on $\Y$.
Finally apply \csy{\refle{dirty}} to the subgroup $\Sigma_n \rtimes (\Sigma_u \x \Sigma_t^n) \subset \Sigma_{tn+u}$.

\refit{injective.admissible.conclusion}~In addition to the preceding points, it remains to note that the injective structure is combinatorial \cite[Proposition~A.2.8.2]{Lurie:HTT} and that
the first bifunctor in \refeq{admissible.product.2} is left Quillen since the pushout product commutes with $\varphi$ and $\C$ is monoidal.
\xpf

\rema \label{rema--injective.admissible}
The injective model structure on  $\Sigma \C$ is tractable if every object of $\C$ is cofibrant.
This applies, for example, for simplicial sets or for simplicial presheaves with the \emph{injective} model structure.
\xrema

\rema
In the case of symmetric spectra in $\C=\sSet_\bullet$ and $R_n = S^n$, the $n$-sphere,
not every cofibration in $\Sigma^{\posi,\inje} \C$ (positive injective structure) is a power cofibration:
the object $(R \t G_1(*_+))^{\t_R2} = R \t G_2 (\Sigma_2 \cdot *_+)$ is cofibrant in $\Sigma_2^\proj \Mod_R$
because its evaluation in degree $2$ is $(* \sqcup *)_+$ on which both copies of $\Sigma_2$ act by permutation.
This object is not cofibrant in $\Sigma_2^\proj \Sigma_2^\inje \sSet_\bullet$.
\xrema

Between the injective and the projective model structure, there is the following minimal strongly $\Y$-admissible model structure.

\prop \label{prop--minimal.admissible}
Suppose that $\Y$ is a \emph{set} (as opposed to a class) of morphisms.
Also suppose that $\C$ is a nice model category (\refde{nice.model.category}).
Then there is a strongly $\Y$-admissible positive model structure.
We call it the \emph{canonical strongly $\Y$-admissible positive model structure}.
\xprop

\pf
We use \cite[Proposition~A.2.6.13]{Lurie:HTT} to construct a combinatorial model structure on each $\Sigma_m\C$ for $m \ge1$.
The weak equivalences will always be $W := \varphi^{-1} (\we_\C)$, as required by \ref{defi--admissible}\refit{admissible.we}.
This is a perfect class (in the sense of loc.~cit.)\ since $\C$ is pretty small \csy{\refle{sequential}}.
In addition we need to define a set $I_m$ of maps in~$\Sigma_m\C$ for $m \ge1$.
These will be the generating cofibrations of a model structure on $\Sigma_m \C$ provided that two conditions are met.
(1)~Any $f\in I_m$ is an h-cofibration in $\Sigma_m\C$.
This will be satisfied as soon as $I_m$ consists of injective cofibrations.
(2)~The class $\inj{I_m}$ is contained in $W$.
This will be satisfied provided that $I_m$ contains $\Sigma_m \cdot I_\C$, where $I_\C$ is a set of generating cofibrations of $\C$, since all maps in $\inj{\Sigma_m \cdot I_\C}$ are projective acyclic fibrations in $\Sigma_m \C$, in particular weak equivalences in $\C$.
\ppar
Starting with $I_1 := I_\C$, we inductively construct $I_m$ by
$$I_m := \Sigma_m \cdot I_\C \cup
\bigcup_{m = m'+m''} (\Sigma_m \cdot_{\Sigma_{m'} \x \Sigma_{m''}} I_{m'} \pp I_{m''})
\cup \bigcup_{m = tn+u, y} \Sigma_m \cdot_{\Sigma_u \x (\Sigma_n \rtimes \Sigma_t^n)} y_u \pp I_{t}^{\pp n}.\eqlabel{I0}$$
The first union runs over partitions of $m$ into \emph{positive} parts.
The second union runs over all multi-indices (of the same size) $t \ge1$, $n \ge1$ where at least one entry $n_i>1$, all $u\ge0$, and
all $y \in \Y_{n}$ (which is a set by assumption; $y_u$ is the $u$th level of $y$).
As usual, we have abbreviated $I_t^{\pp n} := I_{t_1}^{\pp n_1} \pp \cdots \pp I_{t_e}^{\pp n_e}$.
Note that $m'$, $m''$, and the $t_i$ are all strictly less than $m$.
Therefore, $I_{m'}$ etc.\ is defined.
By \refle{injective.admissible}\refit{injective.admissible.1} and \refit{injective.admissible.2}, $I_m$ consists of injective cofibrations.
Moreover, $\Sigma_m \cdot I_\C \subset I_m$, as requested above.
Hence, $W$ and $I$ define a combinatorial model structure on $\Sigma_m \C$.
By design, the functor in \ref{defi--admissible}\refit{admissible.product.2} is a left Quillen bifunctor (note that since we are establishing a positive model structure $m', m''\ge1$ there).
Again by design, the strong admissibility requirement \ref{defi--admissible}\refit{admissible.wedge.2} is met for those multi-indices $n$ where at least one $n_i$ is at least~2.
If all $n_i=1$, then the expression in \ref{defi--admissible}\refit{admissible.wedge.2} reduces to $\Sigma_{t + u} \cdot_{\Sigma_u \x \Sigma_t} y \pp h$ (where $t+u := \sum t_i +u$), which is the $(t+u)$th level of $G_u(y) \pp G_t(h)$.
The latter map is a cofibration in $\Sigma \C$ by the assumption on $\Y$ made in \refde{admissible}.
\ppar
This also shows that the tractability of $\Sigma_{m'} \C$ etc.\ carries over to the one of the newly minted model structure on $\Sigma_m\C$.
\xpf

\lemm \label{lemm--projective.admissible}
If $\C$ is a cofibrantly generated model category (for example, a nice model category in the sense of \refde{nice.model.category}), then the projective model structure on $\Sigma \C$ exists. It is admissible (\refde{admissible}).

If every cofibration $c$ in $\C$ is a \emph{power cofibration} (i.e., $c^{\pp n}$ is a $\Sigma_n$-projective cofibration, see \cite[\S4.5.4]{Lurie:HA}),
then the projective model structure is strongly $\Y$-admissible with respect to the class $\Y_n$ consisting of $\Sigma_n$-projective cofibrations in $\Sigma_n \Sigma \C$.
\xlemm
\pf
The existence and admissibility is standard, see for example~\csy{\refpr{projective.G}}.
As for strong admissibility, the generating projective cofibrations of $\Sigma_t \C$ are given by $\Sigma_t \cdot I_\C$.
The following chain of inclusion shows our claim for generating projective cofibrations in $\Sigma_n \Sigma_u \C$.
The general case follows from this using \csy{\refle{symmetroidal.Y.weakly.sat}}.
$$\eqalignno{\Sigma_u \cdot z \pp (\Sigma_t \cdot I_\C)^{\pp n}& =\Sigma_u \cdot \Sigma_t^n \cdot z \pp (I_\C)^{\pp n}\cr
& \subset \Sigma_u \cdot \Sigma_t^n \cdot z \pp \cof_{\Sigma_n^\proj \C}&\eqlabel{eqnarra1} \cr
& = \Sigma_u \cdot \Sigma_t^n \cdot z \pp \ws{\Sigma_n \cdot I_\C}\cr
&\subset \ws{\Sigma_u \cdot (\Sigma_n \rtimes \Sigma_t^n) \cdot z \pp I_\C}&\eqlabel{eqnarra3} \cr
&\subset \cof_{\Sigma_u \x \Sigma_n \rtimes (\Sigma_t^n)^\proj \C}.\cr}$$
The inclusion \refeq{eqnarra1} holds by assumption.
For \refeq{eqnarra3}, observe that $\Sigma_n$ acts on $\Sigma_t^n$ by permutation.
\xpf

\section{The stable model structure on $R$-spectra}
\label{sect--stable}

\subsection{Definitions}
With this supply of model structures on $\Sigma^\ppos \C$, we turn to model structures on $R$-spectra. To save space, we introduce the following definition:

\defi
\label{defi--context}
A {\it spectral context\/} is a pair $(\Sigma \C,R)$,
where $\C$ is a nice model category (\refde{nice.model.category}),
$\Sigma \C$ is endowed with an admissible model structure (\refde{admissible}),
and $R$ is a commutative monoid in~$\Sigma\C$.
\xdefi

\exam \label{exam--monoids}
In many applications,~$R$ is the free commutative monoid on~$G_1(A)$ for some object $A\in\C$, i.e., $R_n=A^{\t n}$ with $\Sigma_n$ acting by permutations.
In \refpr{I.spaces} we discuss the case $A=1_\C$, the monoidal unit.
More specifically, for the category $\C=\sSet_\bullet$ of pointed simplicial sets (alternatively, $\Top_\bullet$, see \refsect{appendix}) and the pointed circle $A=S^1$, $\Mod_R$ is the category of (simplicial or topological) symmetric $S^1$-spectra.
For this category, the positive stable model structure has been established in \cite[Theorem~14.2]{MandellMaySchwedeShipley:Model}.

The model category used in motivic homotopy theory is $\C=\sPSh_\bullet(\Sm/S)$ (pointed simplicial presheaves on the site of smooth schemes over some base scheme~$S$),
for which we take the pointed projective line $A=(\P^1_S,\infty)$ or, alternatively, $A=\A^1/(\A^1\setminus\{0\})$~\cite{Jardine:Motivic}.
The category~$\Mod_R$ is known as the category of motivic $\P^1$-spectra.
In the projective model structure on \emph{pointed} simplicial presheaves (or any localization thereof), $(\P^1_S,\infty)$ is \emph{not} cofibrant.
This is why we avoid imposing any cofibrancy hypotheses on~$R$, unlike Hovey \cite[\S8]{Hovey:Spectra}.
The flatness of $\C$ ensures that the category of $R$-spectra is replaced by a Quillen equivalent category if $R$ is replaced by a weakly equivalent commutative monoid,
see~\cite[Theorem~4.3]{SchwedeShipley:Algebras}.
This is used in \refsect{Deligne} to construct a strictly commutative $\P^1$-spectrum representing Deligne cohomology.

If $\C$ consists of the Nisnevich $\A^1$-localization of simplicial presheaves with the injective model structure,
the existence of the stable positive model structure has been shown by Hornbostel \cite[Theorem~3.4]{Hornbostel:Preorientations}
in the case where the chosen model structure is the injective model structure.
\xexam

\defi
\label{defi--Mod.R.unstable}
If $(\Sigma\C,R)$ is a spectral context,
then the \emph{unstable model structure} $\Mod_R$ is the one transferred from~$\Sigma\C$ along
the adjunction
$$R\t- : \Sigma\C\leftrightarrows \Mod_R : U.\eqlabel{adjRMod}$$
It is occasionally referred to as the unstable nonpositive model structure.
The \emph{unstable positive model structure} $\Mod_R^\posi$ is the one transferred from $\Sigma^\posi \C$.
As above, $\Mod_R^\ppos$ stands for either $\Mod_R$ or $\Mod_R^\posi$.
\xdefi

The next step is to localize the unstable model structure on $R$-modules to obtain the stable model structure.
Consider the adjunction
$$F_n : \C \leftrightarrows \Mod_R : \Ev_n\eqlabel{adjunction.F.ev}$$
obtained by composing the adjunctions $\Sigma_n \cdot - \colon \C \leftrightarrows \Sigma_n \C$, \refeq{adjunction.Sigma.C} and \refeq{adjRMod}.
The right adjoint evaluates at the $n$th level
(after forgetting the $R$-module structure and the $\Sigma_n$-action).
The left adjoint is given by $F_n(X)=G_n(\Sigma_n \cdot X) \t R$.

\defi \label{defi--stable}
Suppose $(\Sigma \C,R)$ is a spectral context.
The \emph{stable (positive) model structure} $\Mod_R^{\stab, \ppos}$
is the symmetric monoidal left Bousfield localization
(i.e., left Bousfield localization in the bicategory of symmetric monoidal model categories, see \refsect{preliminaries})
of the (positive) unstable model structure~$\Mod_R^{\ppos}$ on $R$-modules with respect to the set
$$\xi^R := \{ \xi_n^R \colon F_n(\mcr R_n) \to R\mid \ n\ge0 \}.\eqlabel{class.xi}$$
Here $\mcr$ is the cofibrant replacement functor in~$\C$.
\xdefi

\rema \label{rema--stable}
For $n\ge0$ (in the positive case for $n\ge1$),
the map $\xi_n$ above is the derived adjoint of the identity map~$R_n := \Ev_n(R)\to\Ev_n(R)\in \C$ with respect to the Quillen adjunction \refeq{adjunction.F.ev}.
The derived adjoint can be computed by precomposing with the weak equivalence $\mcr R_n\to R_n$ so that the source is cofibrant
and taking the ordinary adjoint.

If $\C$ is $\V$-enriched, then $\Mod_R^{\stab, \ppos}$ is the $\V$-enriched monoidal localization by \csy{\refre{monoidal.enriched}}.
The name ``stable model structure'' for this model structure is standard, even though this model structure is not stable for all $R$, for example for $R_n=1_\C$ (see the discussion following \refpr{I.spaces}).
See, however, \refpr{stable.R.stable} for a sufficient criterion for stability.

Suppose $R$ is the free commutative monoid on $G_1(R_1)$, i.e., $R_n = R_1^{\t n}$.
Suppose further that $R_1$ is either cofibrant in $\C$ or there is a cofibration $1 \r R_1$.
Then the above localization agrees with the one with respect to $\xi_1$ only, since
$F_1(\mcr R_1)^{\t n} = F_n((\mcr R_1)^{\t n})$ and $(\mcr R_1)^{\t n} \sim \mcr (R_1^{\t n})$ by \cop{\refle{cofibrant.replacement}}.

For the projective structure on $\Sigma \C$ and $R = \Sym(G_1(R_1))$ with a cofibrant object $R_1 \in \C$, the (nonpositive) stable model structure $\Mod_R^\stab$ has been defined by Hovey in \cite[Definition~8.7]{Hovey:Spectra} as the localization (in the bicategory of mere model categories, i.e., disregarding the monoidality and $\V$-enrichment of $\Mod_R$) with respect to the set of maps
$$\zeta_n(C)\colon F_{n+1}(C \t R_1) \r F_n (C)$$
adjoint to the map $C \t R_1 \r \Ev_{n+1} F_n(C) = \Sigma_{n+1} \cdot C \t R_1$ given by the identity element of $\Sigma_{n+1}$.
Here $n\ge0$ and $C$ runs through the (co)domains of generating cofibrations of $\C$.
Hovey's definition agrees with the one above.
Indeed, %by \csy{\refpr{Bousfield}},
the monoidal localization with respect to $\xi_1 = \zeta_0(1)$ is the (ordinary) localization with respect to the set $F_n(C) \t_R \mcr \zeta_0(1)$, which is equivalent to the one by $F_n(C) \t_R \zeta_0(1)$ by the flatness of $\Mod_R$.
One checks that this map is just $\zeta_n (C)$.
The objects $F_n(\C)$ are precisely the (co)domains of generating cofibrations of the (projective, nonpositive) model structure $\Mod_R$.

For the same type of commutative monoid Gorchinskiy and Guletski\u\i~ define the stable positive structure to be the localization with respect to Hovey's class, but for $n\ge1$.
Both their definition and \refde{stable} have the property that positive stable weak equivalences agree with nonpositive stable equivalences \cite[Theorem~9]{GorchinskiyGuletskii:Positive},
\refpr{stable.we}, so that the model structure in loc.~cit.\ is Quillen equivalent to the one defined above.
\xrema

\subsection{Existence}

\theo
\label{theo--stable.R.general}
For any spectral context $(\Sigma\C,R)$ the model category $\Mod_R^{\stab,\ppos}$ exists.
Its fibrant objects are those objects~$W$ that are fibrant in $\Mod_R^{\ppos}$ and such that the derived internal hom in $\Mod_R^{\ppos}$,
$$\mfr \IHom(\xi_n, W)$$ is a weak equivalence for all $n\ge0$.
\xtheo

\pf
We write $\cof$ and $\AC$ for the cofibrations and acyclic cofibrations of a model category.
The first step for the existence of $\Mod_R^{\stab, \ppos}$ is the monoid axiom for $\Sigma^\ppos \C$.
For the pushout product axiom, it is enough to check $I \pp I \subset \cof(\Sigma^\ppos\C)$ and $I \pp J \cup J \pp I \subset \AC(\Sigma^\ppos\C)$.
Here $I$ ($J$) are the generating (acyclic) cofibrations of~$\Sigma^\ppos \C$.
They are of the form $G_n (f)$, where $n\ge0$ (respectively $n\ge1$ in the positive case) and $f \in \Sigma_n \C$ is a generating (acyclic) cofibration.
Using \refeq{product.Sigma.C}, we obtain the monoidality of the model structure by \refde{admissible}\refit{admissible.product.2}.

The category $\Sigma^{\ppos} \C$ is h-monoidal, i.e., for any $X \in \Sigma^\ppos \C$ and any (acyclic) cofibration $u \in \Sigma^\ppos \C$, $X \t u$ is an (acyclic) h-cofibration in $\Sigma^\ppos \C$.
This has to be checked only for generating (acyclic) cofibrations $u=G_t(h)$, $h \in \Sigma_t \C$ and also only for $X=G_t(Z)$.
Here we use that (acyclic) h-cofibrations in an h-monoidal model category are stable under finite coproducts \cite[Lemma~1.3]{BataninBerger:Homotopy} and therefore, using the pretty smallness and
\csy{\refle{sequential}}, under countable coproducts.
As $\Sigma_u \x \Sigma_t$ is a subgroup of $\Sigma_{u+t}$, we obtain a noncanonical isomorphism \csy{\refle{dirty}}
$$\varphi (\Sigma_{t+u} \cdot_{\Sigma_u \x \Sigma_{t}} Z \t h^{\pp n})
={\Sigma_{t+u}\over \Sigma_u \x \Sigma_{t}} \cdot \varphi \left (Z \t h \right),$$
so the h-monoidality of $\C$ implies the one of $\Sigma^\ppos \C$, since $\varphi$ detects h-cofibrations.
Then, $\Sigma^\ppos \C$ satisfies the monoid axiom by pretty smallness and h-monoidality \csy{\refle{i.monoidal.monoid.axiom}}.

The unstable model structure $\Mod_R^\ppos$ exists by the monoid axiom for $\Sigma^\ppos \C$~\cite[Theorem~4.1(2)]{SchwedeShipley:Algebras}.
It is tractable: by definition, $\C$ and hence $\Sigma \C$ and thus also $\Mod_R$ are locally presentable.
Moreover, $\Sigma_m \C$ is tractable and cofibrantly generated by assumption (see \refde{admissible}\refit{admissible.tractable}) and therefore so is $\Sigma^\ppos \C$ and hence $\Mod_R^\ppos$.

The stable model structure $\Mod_R^{\stab, \ppos}$ exists since $\Mod_R^\ppos$ is combinatorial. %, see e.g. \csy{\refpr{Bousfield}}.

The description of fibrant objects of $\Mod_R^{\stab, \ppos}$ is an application of
\cite[4.46.4]{Barwick:Left}.
%\csy{\refle{fibrant.monoidal.localization}}.
\xpf

As an example, we examine the special case $R = E$, where $E$ is the free commutative monoid in~$\Sigma\C$ on the monoidal unit.
Its levels are given by $E_n=1_\C$, the monoidal unit (with the trivial $\Sigma_n$-action).
In this case, $E$-modules coincide with $I$-spaces, as defined by Sagave and Schlichtkrull \cite{SagaveSchlichtkrull:Diagram}.
By definition, these are functors from the category $I$ of finite sets and injections to $\C$.
Indeed, an $E$-module~$X$ is the same as a sequence of objects $X_n \in \Sigma_n \C$ with a $\Sigma_n$-equivariant bonding map $X_n \cong X_n\t1 \r X_{n+1}$.
This datum is equivalent to specifying an $I$-space whose value on objects and isomorphisms $\sigma \in \Sigma_n$ is given by the $X_n$ and whose value on injections is given by compositions of bonding maps.
What is more, the stable model structure on $I$-spaces defined in loc.~cit.\ agrees with the stable model structure on $\Mod_E$:

\prop \label{prop--I.spaces}
Let $\C$ be a nice model category (\refde{nice.model.category}).
We equip $\Sigma \C$ with the projective model structure and consider the resulting (un)stable (positive) projective model structures on $E$-modules.
The (un)stable (positive) projective structures on $\Mod_E$ and the category $I\C$ of $I$-spaces coincide, i.e., all five classes of maps are preserved under the above equivalence.
\xprop

\pf
\ppar
The unstable positive projective model structures on $E$-modules and $I$-spaces coincide since they are both transferred from $\prod_{n\ge1} \C$.

For the stable structures it is enough to prove that stable weak equivalences of $I$-spaces correspond to stable weak equivalences of $E$-modules.
Both model structures are left Bousfield localizations, so it is sufficient to establish
that the stably fibrant $E$-modules are exactly the stably fibrant $I$-spaces.
By \refth{stable.R.general}, stably positively fibrant $E$-modules are precisely those $E$-modules~$X$ that are unstably positively fibrant
and $\hHom(F_n(\mcr1) \r E, X)$ is a weak equivalence in $\Mod_E^{\posi}$ for all $n\ge0$, or, equivalently, the $r$th level ($r\ge1$) of this is a weak equivalence in~$\C$.
Here and below $\hHom$ is the derived internal hom with respect to the positive stable model structure on $E$-modules.
This condition is automatic for $n=0$ since $F_0(\mcr1) \r E=F_0(1)$ is always a weak equivalence:
in spectral level $r$, the map is $E_r \t \mcr1 =1 \t \mcr1 \stackrel \cong \r \mcr1 \stackrel \sim \r1$.
A cofibrant replacement of $E$ in the positive model structure is given by $E^\posi$ that is $\emptyset$ (respectively $1_\C$) in spectral level~0 (respectively $\ge1$).
We can compute the $r$th spectral level of $\hHom(E, X)$ as $\IHom(E^\posi, X)_r = X_r$.
Moreover, for $n\ge1$, $F_n(\mcr1)$ is positively cofibrant, so that there are weak equivalences $\hHom(F_n(\mcr1),X)_r = \IHom(F_n(\mcr1),X)_r = X_{r+n}$.
(The second weak equivalence uses the unit axiom of~$\C$ (implied by the flatness of~$\C$) in the equivalent formulation of \cite[Lemma 4.2.7]{Hovey:Model}.)
In other words stably positively fibrant $E$-modules are those unstably fibrant $E$-modules such that $X_r \r X_{r+n}$ is a weak equivalence for all $n \ge0$ and all $r \ge1$.
These are exactly the stably positively fibrant $I$-spaces \cite[\S3.1]{SagaveSchlichtkrull:Diagram}.
The same works in the nonpositive case.
\xpf

\subsection{Invariance properties}

In this section, we establish a few Quillen equivalences between different model structures on $R$-spectra.
With a slightly different definition (see \refre{stable}), \refpr{stable.we} is due to Gorchinskiy and Guletski\u\i~\cite[Theorem~9]{GorchinskiyGuletskii:Positive}.

\prop
\label{prop--stable.we}
For any spectral context $(\Sigma \C,R)$
the class of stable positive weak equivalences $\we_{\stab,\posi} := \we_{\Mod_R^{\stab,\posi}}$ agrees with the stable (nonpositive) weak equivalences $\we_\stab := \we_{\Mod_R^\stab}$.
In particular, the categories $\Mod_R^{\stab}$ and $\Mod_R^{\stab, \posi}$ are Quillen equivalent.
\xprop

\exam
We continue with the notation of \refpr{I.spaces}.
By \cite[Theorem~9.1]{Hovey:Spectra}, $\Mod_E^{\stab}$ and therefore $\Mod_E^{\stab, \posi}$ is Quillen equivalent to $\C$.
For example, if $\C$ is not symmetric flat (such as $\C = \sSet$), it is nonetheless Quillen equivalent to $\Mod_E^{\stab, \posi}$ (or $I$-spaces in~$\C$),
which is symmetric flat.
This point of view goes back to Jeff Smith.
\xexam

\pf
We essentially reproduce the proof of \cite[Theorem~9]{GorchinskiyGuletskii:Positive}.
We will not explicitly mention that a model structure on $\Mod_R$ is nonpositive or unstable, but will always indicate positivity and/or stability where necessary.
Moreover, a superscript indicates a certain model-categorical operation related to the model category structure in question.
For example, $\mcr$ is the cofibrant replacement functor in $\Mod_R$, $\mcr^{\posi}$ the one of $\Mod_R^{\posi}$.
Similarly, $\hMap^{\stab, \posi}$ is the derived mapping space of $\Mod_R^{\stab, \posi}$.
By definition, there is a Quillen adjunction, where $\IHom$ denotes the internal hom:
$$F_1 \mcr (R_1) \t_R - : \Mod_R \leftrightarrows \Mod_R^{\posi} : \Theta_1 := \IHom(F_1 \mcr R_1, -).\eqlabel{adjunction.positive.unstable}$$
It localizes to a Quillen adjunction
$$F_1 \mcr (R_1) \t_R - : \Mod_R^{\stab} \leftrightarrows \Mod_R^{\stab, \posi} : \Theta_1.\eqlabel{adjunction.positive.stable}$$
In fact, $F_1 \mcr (R_1) \t_R^\LL \xi^R$ is weakly equivalent to
$F_1 \mcr (R_1) \t_R \xi^R$ by the flatness of $\Mod_R$.
As $F_1 \mcr (R_1)$ is cofibrant in $\Mod_R^{\posi}$,
the latter set is contained in the monoidal saturation of $\xi^R$ with respect to the model structure $\Mod_R^{\posi}$ (recall that the monoidal saturation can be computed as the class $\{\xi^R \t E, E \in \Mod_R^\posi \text{ cofibrant }\}$.)
Therefore the derived functor of the left adjoint sends $\xi^R$ to weak equivalences in $\Mod_R^{\stab, \posi}$, which shows that \refeq{adjunction.positive.stable} is a Quillen adjunction.

We first prove two preliminary claims.
The first claim is that any $f \in \we_{\posi}$ is a stable (nonpositive) weak equivalence.
Both $\we_\stab$ and $\we_{\posi}$ are preserved by (unstable nonpositive) fibrant replacement,
so that we may assume that $f$ is a map between nonpositively, a fortiori positively fibrant objects.
By Brown's lemma (applied to \refeq{adjunction.positive.unstable}), $\Theta_1 (f) \in \we \subset \we_\stab$.
Let $f \stackrel \sim \r f'$ be the fibrant replacement of $f$ in the stable structure.
In the following commutative diagram, $\sim$ indicates a stable equivalence.
$$\xymatrix{
f = \IHom (F_0 1, f) \ar[r]  \ar[d]^\sim
& \IHom (F_0 \mcr1, f) \ar[r] \ar[d]
& \Theta_1(f)
\ar[d]
\cr
f' = \IHom (F_0 1, f') \ar[r]^{*}_\sim
& \IHom (F_0 \mcr1, f') \ar[r]^{**}_\sim
& \Theta_1(f').
}$$
The map $*$ is a stable weak equivalence since $F_0 (\mcr1) \t_R Y\r F_0(1) \t_R Y$ is a weak equivalence in $\Mod_R$ (and therefore $\Mod_R^{\stab}$) for any cofibrant object $Y \in \Mod_R$ by the flatness of $\Mod_R$ (\refpr{symseq.general}).
The map $**$ is a stable weak equivalence by the very definition of this model structure.
Consequently, in the homotopy category $\Ho(\Mod_R^\stab)$, $f$ is a retract of the isomorphism $\Theta_1 (f)$, so that $f$ is also a stable weak equivalence.
This finishes the first claim.

The second claim is that for any fibrant object $Z \in \Mod_R^{\stab, \posi}$, the map $\IHom(\xi_1, Z)\colon Z \r \Theta_1 (Z)$ is an (unstable) weak equivalence in $\Mod_R^{\posi}$.
Indeed, for any $n \ge1$ and any cofibrant object $T \in \Sigma_n \C$,
$$\eqalign{\hMap^{\posi}(F_n(T), Z) & \stackrel{\xi_1} \lr \hMap^{\posi}(F_n(T) \t_R F_1(\mcr R_1), Z)\cr
& \sim \hMap^{\posi}(F_n(T), \Theta_1(Z))\cr}$$
are weak equivalences, the first by the definition of the \emph{monoidal} Bousfield localization, the second by (homotopy) adjunction.
Since the objects $F_n(T)$ are homotopy generators of $\Mod_R^{\posi}$, we are done with the second claim.

The first claim implies that there is a Quillen adjunction
$$\id : \Mod_R^{\stab, \posi} \leftrightarrows \Mod_R^{\stab} : \id.\eqlabel{adjunction.stable.identity}$$
Indeed, $\mcr^{\posi}(\xi_n^R)$ is positively and therefore (by the first claim) nonpositively stably weakly equivalent to $\xi_n^R$.
Therefore, any fibrant object $T \in \Mod_R^\stab$, is also fibrant in $\Mod_R^{\stab, \posi}$.
For any $X \in \Mod_R$,
the natural map of derived mapping spaces (in $\Mod_R^{\stab}$ and $\Mod_R^{\stab, \posi}$, respectively) induced by the transformation of cofibrant replacement functors $\mcr^{\posi} X \r \mcr X$,
$$\hMap^\stab(X, T) \r \hMap^{\stab, \posi}(X, T)$$
is a weak equivalence.
Indeed, $\mcr^{\posi} X \r \mcr X$ is a positive weak equivalence and therefore a stable (nonpositive) equivalence by the first claim.

We finally prove the proper statement.
For a morphism $f$ and an object $Z \in \Mod_R$, we consider the commutative diagram
whose horizontal maps stem from the Quillen adjunction \refeq{adjunction.stable.identity}:
$$\xymatrix{
\hMap^\stab(f,Z) \ar[r] \ar[d] &
\hMap^{\stab, \posi}(f, Z) \ar[d] \cr
\hMap^\stab (f,\Theta_1 Z) \ar[r] &
\hMap^{\stab, \posi}(f, \Theta_1 Z).
}$$
Suppose $f$ is in $\we_{\stab, \posi}$, so that $\hMap^{\stab, \posi}(f, -)$ is a weak equivalence.
For any fibrant object $Z \in \Mod_R^{\stab}$, the top horizontal map is a weak equivalence (of arrows, i.e., a weak equivalence of source and target) by the above consequence of the first claim.
Thus $\hMap^\stab(f,Z)$ is a weak equivalence, i.e., $f$ is in $\we_{\stab}$.

Conversely, suppose $f \in \we_{\stab}$ so that $\hMap^\stab(f, -)$ is a weak equivalence.
For any fibrant object $Z \in \Mod_R^{\stab, \posi}$, $\Theta_1(Z)$ is fibrant in $\Mod_R^{\stab}$ by \refeq{adjunction.positive.stable}.
Hence, by the consequence of the first claim, the bottom horizontal map is a weak equivalence.
By the second claim $Z \r \Theta_1(Z)$ is in $\we_{\posi} \subset \we_{\stab, \posi}$, hence the right vertical map is a weak equivalence.
We conclude that $\hMap^{\stab, \posi}(f, Z) $ is a weak equivalence so that $f \in \we_{\stab, \posi}$.
\xpf

\prop
\label{prop--stable.R.admissible.choice}
For any spectral context $(\Sigma \C,R)$,
different choices of admissible model structures on $\Sigma \C$ yield Quillen equivalent model structures on $\Mod_R^{\stab, \ppos}$.
\xprop

\pf
By \refde{admissible}\refit{admissible.we}, weak equivalences are the same in all admissible model structures,
so the same is true for the transferred model structures on $\Mod_R^{\ppos}$, which are therefore Quillen equivalent.
This localizes to a Quillen equivalence $\Mod_R^{\stab, \ppos}$ between different choices of admissible model structures
since they are monoidal localizations with respect to the weakly equivalent sets $\xi^R$ (but different because
of different choices of cofibrant replacements in the definition of~$\xi^R$) of morphisms.
\xpf

The next result compares spectra in different spectral contexts.
Its proof uses the following lemma.

\lemm\label{lemm--adjunction.monoidal.localization}
If $F : \C \leftrightarrows \C' : G$ is a Quillen equivalence of monoidal model categories such that $F$~is strong monoidal, then (assuming the left Bousfield localizations exist) there is a Quillen equivalence
$$F : \D := \BL_{S^\t} \C \leftrightarrows \D' := \BL_{\ldf F(S)^\t} \C' : G.$$
\xlemm

\pf
The class $F(\CO_\C)$ is a class of homotopy generators of~$\C'$.
Hence $\D'$ can be computed as the (nonmonoidal) localization with respect to the class $F(\CO_{\C}) \t^\LL \ldf F(S) = F(\CO_\C \t^\LL S)$.
Thus, by \cite[Proposition~3.3.18, Theorem~3.3.20]{Hirschhorn:Model}, the left Quillen functor $\C \r \C' \r \D'$ factors over a left Quillen functor $\D \r \D'$ since $\ldf F(\CO_\C \t^\LL S)$ consists of weak equivalences in~$\D'$.
Moreover, $\D \leftrightarrows \D'$ is a Quillen equivalence if $\C \leftrightarrows \C'$ is one.
\xpf

\prop
\label{prop--stable.R.transport.2}
Suppose
$$\varphi : (\Sigma \C, R) \rightleftarrows (\Sigma \D, S) : \psi$$
is an adjunction of spectral contexts, by which we mean an adjunction $\varphi : \C \rightleftarrows \D : \psi$ with $\varphi$ strong symmetric monoidal, and a map of commutative monoids $R \r \psi (S)$.

Suppose, moreover, that the adjunction $\Sigma \C \rightleftarrows \Sigma \D$ is a Quillen equivalence.
For example, we can take $\C=\D$ with the projective or the injective model structures.
Suppose in addition that $\varphi (\mcr1) \sim \varphi(1)$.

Finally, suppose that one of the following conditions is satisfied:
\begin{enumerate}[(a)]
\item
\label{item--case.b}
The right adjoint $\psi$ preserves all weak equivalences and $R\to\psi(S)$ is a weak equivalence  in~$\Sigma\C$.
\item
\label{item--case.a}
There is a weak equivalence $R' \r R$, where $R'$ is a commutative monoid in $\Sigma \C$ that is cofibrant as an object of $\Sigma \C$.
Moreover, the composite $\varphi(R') \r \varphi(R) \r S$ is a weak equivalence in $\Sigma \D$.
\end{enumerate}
Then there is a Quillen equivalence
$$\varphi_*=S\otimes_{\varphi(R)}\varphi(-) : \Mod_R^{\stab, \ppos} \rightleftarrows \Mod_{S}^{\stab, \ppos} : \psi.$$
\xprop

\exam
For $\C = \D$ and any weak equivalence $R \r S$, we obtain a Quillen equivalence $\Mod_R^{\stab, \ppos} \rightleftarrows \Mod_S^{\stab, \ppos}$. In other words, the assignment $R \mapsto \Mod_R^{\stab, \ppos}$ is ``derived''.

Another important case is $\varphi(R)=S$, in which case we obtain a criterion for replacing $\C$ with a Quillen equivalent model category.
\xexam

\pf
First note that $\varphi(R)$ and $\psi(S)$ is a commutative monoid since $\varphi$ is strong monoidal and $\psi$ is lax monoidal.
We claim that the Quillen equivalence $\Sigma \C \r \Sigma \D$ extends to a Quillen equivalence $\Mod_R \rightleftarrows \Mod_S$.
In Case \refit{case.b}, this follows from \cite[Theorem~3.12(2)]{SchwedeShipley:Equivalences} and the flatness of $\Mod_R$ (see the proof of \refpr{symseq.general}), which yields a Quillen equivalence $\Mod_R \rightleftarrows \Mod_{\psi (S)}$ \cite[Theorem~4.3]{SchwedeShipley:Algebras}.
In Case \refit{case.a}, again using flatness, we can replace $R$ by $R'$.
The forgetful functor $\Mod_R \r \Sigma\C$ preserves cofibrant objects,
which follows from \csy{\refpr{transfer.basic}\refit{G.preserves.cofibrations}}, using that $R$ is cofibrant as a symmetric sequence.
Furthermore, for any cofibrant $R$-module~$X$ the map $\varphi(X)\to\varphi_*(X)$ is a weak equivalence,
using cellular induction, flatness, and h-monoidality.
Then, the claim follows from Quillen's criterion for Quillen equivalences.

\ppar
We now upgrade this to a Quillen equivalence $\Mod_R^\stab \rightleftarrows \Mod_S^\stab$.
This also settles the claim for positive model structures by \refpr{stable.we}.
By \refle{adjunction.monoidal.localization}, $\Mod_R^\stab$ is Quillen equivalent to the monoidal localization $\BL^\t_{\ldf \varphi_*(\xi^R)} \Mod_S$.
It is enough to show that the natural map
$$\ldf \varphi_*(\xi_n^R) \r \xi_n^S$$
is a weak equivalence of maps in $\Mod_S$. Recall from \refeq{class.xi} that
$\xi_n^R \colon R \t G_n (\Sigma_n \cdot \mcr R_n) \r R$.

We first treat Case \refit{case.b}.
For the codomains of these maps, we use the assumption $R \stackrel \sim \r \psi(S)$, which computes $\rdf \psi(S)$, since $\psi$ preserves all weak equivalences.
This is equivalent to $\ldf \varphi(R) \stackrel \sim \r S$.
The domains of $\xi_n^R$ are cofibrant, so it is enough to show that $\varphi_*(R \t G_n (\Sigma_n \cdot \mcr R_n)) = S \t G_n (\Sigma_n \cdot \varphi(\mcr R_n)) \stackrel {S \t \alpha} \r S \t G_n (\Sigma_n \cdot \mcr S_n)$ is a weak equivalence.
For some cofibrant replacement $\mcr S \stackrel \sigma \r S$ of $S$ in $\Sigma \C$, $\mcr S \t \alpha$ is a weak equivalence by Brown's lemma. Moreover, tensoring $\sigma$ with the cofibrant (co)domains of $\alpha$ gives a weak equivalence by flatness of $\Sigma \C$.
Therefore $S \t \alpha$ is a weak equivalence.

We now do Case \refit{case.a}.
Applying Case \refit{case.b} to $\C = \D$ and $R'$ vs. $R$, we may assume $R$ is cofibrant as an object of $\Sigma \C$.
For the domains of the maps, $\ldf \varphi_*(R \t G_n(\Sigma_n \cdot \mcr R_n)) \sim S \t G_n (\Sigma_n \cdot \varphi (\mcr R_n)) \r S \t G_n (\Sigma_n \cdot \mcr S_n)$ is a weak equivalence:
as above, use flatness to replace $S \t - $ by $\mcr S \t -$, then use Brown's lemma and the assumption $\varphi(R_n) \stackrel \sim \r S_n$.
For the codomains of the maps, $\ldf \varphi_* (R)$ can be computed as $\varphi_* (R \t \mcr1)$ since $R$ is cofibrant as an object of $\Sigma \C$.
Now, $\varphi_*(R \t \mcr1)=S \t \varphi(\mcr1)$ is weakly equivalent to $S\t1 = S$:
this follows again using a cofibrant replacement $\mcr S \stackrel \sigma \r S$ and flatness as above, and also using that $\sigma\t1$ is an isomorphism, in particular a weak equivalence.
\xpf

\subsection{Monoidal and other basic properties}

\prop
\label{prop--stable.R.stable}
Suppose $(\Sigma \C,R)$ is a spectral context such that $\C$ is pointed.
We write $S^1 \in \C$ for some cofibrant representative of the suspension of the monoidal unit $1_\C$, i.e., the homotopy pushout $* \sqcup^h_{1_\C} *$.
Suppose that $R$ is such that $R_1$ is weakly equivalent to $S^1 \t B$ for some cofibrant object $B \in \C$.
Then the model structure $\Mod_R^{\stab, \ppos}$ is stable in the sense that it is pointed and the suspension and loop functors are inverse Quillen equivalences on $\Mod_R^{\stab, \ppos}$ \cite[Definition~2.1.1]{SchwedeShipley:Stable}.
\xprop

\pf
By \refpr{stable.we}, we only need to consider the nonpositive case.
For a cofibrant object $X \in \Mod_R^{\stab}$, the suspension $\Sigma X$ is weakly equivalent to $X \t S^1 = X \t F_0(S^1)$, where $F_n$ is defined in \refeq{adjunction.F.ev}.
As $F_1$ is a left Quillen functor, $F_0(S^1) \t F_1(B) = F_1 (S^1 \t B)$ is weakly equivalent to $F_1(\mcr(R_1)) = R \t G_1(\mcr(R_1))$, where $\mcr$ is the cofibrant replacement functor.
By definition of the stable model structure, this is stably weakly equivalent to $F_0(1_\C) = R$, which is the monoidal unit in $\Mod_R$.
\xpf

We have the following basic properties of the stable model structure.
For $\C = \sSet_\bullet$, the monoid axiom of $S^1$-spectra was shown independently by Hess and Shipley \cite[Notation~5.5]{HessShipley:Waldhausen}.
We do not address the question of \emph{right} properness of $\Mod_R^{\stab, \ppos}$.
For simplicial and topological spectra and motivic $\P^1$-spectra,
the right properness was shown by Schwede \cite[Theorem III.4.11]{Schwede:Spectra} and Jardine \cite[Theorem~4.15]{Jardine:Motivic}.
Both proofs use that these two examples are actually stable model categories, so any generalization of this statement might use \refpr{stable.R.stable}.

\prop
\label{prop--symseq.general}
Let $(\Sigma \C, R)$ be a spectral context (\refde{context}).
Then the stable (positive) model structure $\Mod_R^{\stab, \ppos}$ is nice (\refde{nice.model.category}).
In particular, it is left proper and satisfies the monoid axiom.
\xprop

\pf
We first show that $\Sigma^\ppos \C$ is nice, except that $\Sigma^\posi \C$ is not flat.
(For \emph{any} map $y \in \C$, $G_0(y)$ is a weak equivalence in $\Sigma^{\posi} \C$, but $y \pp G_1(c)$ is not.
However, $\Sigma \C$ is flat, as we show below.)

The tractability of $\Sigma^\ppos \C$ is a basic property of transfer \csy{\refpr{transfer.basic}}.
Similarly, $\C$ is pretty small, hence so is the projective model structure on $\Sigma_m \C$, and therefore the given admissible model structure on $\Sigma_m \C$.
The symmetric monoidality and h-monoidality of $\Sigma^\ppos \C$ were already shown in the proof of \refth{stable.R.general}.

Using the h-monoidality of~$\Sigma \C$ it is enough to check flatness of $\Sigma \C$ for generating cofibrations \csy{\refth{symmetric.weakly.saturated}}.
Thus we need to show $y \pp G_n(c)$ is a weak equivalence for any weak equivalence $y \in \Sigma \C$ and any cofibration $c$ in $\Sigma_n \C$, $n \ge0$.
We have $y \pp G_n(c) = \coprod_{r \ge0} G_{n+r}(\Sigma_{n+r} \cdot_{\Sigma_r \x \Sigma_n} y_r \pp c)$.
It is enough to see that $\Sigma_{n+r} \cdot_{\Sigma_r \x \Sigma_n} y_r \pp c$ is a weak equivalence.
Again by \csy{\refle{dirty}}, it is isomorphic, in $\C$, to a finite coproduct of copies of $y_r \pp c$, which is a weak equivalence in $\C$ by the flatness of $\C$.
Moreover, by the h-monoidality, $\text{(co)dom}(y_r) \t c$ is an h-cofibration, so that $y_r \pp c$ is a couniversal weak equivalence by \csy{\refle{preparation.flat.weakly.sat}}.
These are stable under finite coproducts in any model category.

The nicety properties of $\Sigma^\ppos \C$ transfer from~$\Sigma^\ppos \C$ to $\Mod_R^\ppos$ by \csy{\refth{commutative.monoid}} and localize to $\Mod_R^{\stab, \ppos}$ by \csy{%\refpr{Bousfield},
\refpr{monoidal.localization}}.

The monoid axiom for $\Mod_R^{\stab, \posi}$ follows from pretty smallness and h-monoidality \csy{\refle{i.monoidal.monoid.axiom}}.
\xpf

\prop
\label{prop--F.0.preserves.compact}
Let $(\Sigma \C, R)$ be a spectral context such that $R_n$ is monoidally compact, i.e., the derived internal hom $\rdf \Hom(R_n, -)$ preserves filtered homotopy colimits.
Then the derived left adjoint of the following adjunction (cf.~\refeq{adjunction.F.ev}) preserves homotopy compact objects:
$$F_0 : \C \leftrightarrows \Mod_R^{\stab, \ppos} : \Ev_0.$$
\xprop

\pf
Preservation of compact objects by $\ldf F_0$ is equivalent to preservation of filtered homotopy colimits by $\rdf \Ev_0$.
Since the functor $\C \gets \Sigma \C \gets \Mod_R : \Ev_0$ preserves all colimits and is automatically derived, it is enough to note that filtered homotopy colimits in $\Mod_R$ preserve local objects, i.e., $\Omega$-spectra~$W$ such that the derived internal hom (with respect to the unstable model structure) $\rdf \Hom(\xi_n^R, W)$ is a weak equivalence (\refth{stable.R.general}).
\xpf

\subsection{Symmetricity properties}
\label{sect--symmetricity.2}

We now establish the symmetricity properties of~$\Mod_R^{\stab,\posi}$.
These are the technical key points of this entire paper.
For example, the symmetric h-monoidality of $\Mod_R^{\stab, \posi}$ will give rise to the existence of model structures on operadic algebras,
while the symmetric flatness is responsible for the rectification of algebras over operads.
See also Remarks \ref{rema--explain.symmetric.h.monoidal} and \ref{rema--explain.symmetric.flat} for further explanation of the notion of symmetric h-monoidality and flatness.

In the generality stated below, these properties are new.
However, various aspects of this description are well known and, in fact, the basic principle of these claims is a very simple idea, due to Smith:
many model categories (say, $\C = \sSet$) are not symmetric flat: for a $\Sigma_n$-equivariant weak equivalence $y$ and a cofibrant object $X \in \C$, $y \t_{\Sigma_n} X^{\t n}$ is not usually a weak equivalence.
The point of introducing $\Sigma^\posi \C$ is to allow only those cofibrant objects $X \in \Sigma \C$ (and more generally, cofibrations $x$) that are trivial in level~0.
For example, for $X=G_1(X')$, $X^{\t n}=G_n(\Sigma_n \cdot X'^{\t n})$, so that modding out the $\Sigma_n$-action on $y \t X^{\t n}$ preserves the homotopical properties of $y$.
This property of $\Sigma^\posi \C$ is then propagated to various symmetricity properties of $\Mod_R^{\stab, \ppos}$.

A special case of symmetric flatness (namely the case where the weak
equivalence $y \in \Sigma_n \Mod_R$ is given by the projective cofibrant replacement of $1_{\Mod_R} = R$, $\E\Sigma_n \r R$)
is due to Gorchinskiy and Guletski\u\i~\cite[Theorem~11]{GorchinskiyGuletskii:Positive}.
They prove this statement under the assumption that every positive projective cofibration in $\Mod_R^{\posi}$ (i.e., the transfer of the positive projective structure on $\Sigma \C$ to $R$-modules) is a power cofibration.
As was explained in \refle{projective.admissible}, this condition ensures that the projective structure is strongly admissible (which only holds in very special cases).
The more general symmetric flatness will be used to show the operadic rectification (\refth{rectification.spectra}).

\prop
\label{prop--stable.R.positive}
For any spectral context $(\Sigma \C, R)$, the model category $\Mod_R^{\stab, \posi}$ is symmetric flat and symmetric h-monoidal.
If, moreover, the model structure on $\Sigma \C$ is strongly $\Y$-admissible for a collection $\Y = (\Y_n)$, $\Y_n \subset \Mor \Sigma_n \Sigma \C$,
then the (acyclic) cofibrations of $\Mod_R^{\stab, \posi}$ form an (acyclic) $\ws{R \t \Y}$-symmetroidal class in $\Mod_R^{\stab}$.
In particular, $\Mod_R^{\stab, \posi}$ is $\ws{R \t \Y}$-symmetroidal in this case.
\xprop

\pf
We first show the following statements:
\begin{enumerate}[(a)]
\item
\label{item--symseq.symmetric.h.monoidal}
The (acyclic) cofibrations of $\Sigma^{\posi} \C$ form an (acyclic) symmetric h-monoidal class in $\Sigma \C$:
for any finite family $v$ of positive (acyclic) cofibrations, and any object $Y \in \Sigma_n \Sigma \C$, $Y \t_{\Sigma_n} v^{\pp n}$ is an (acyclic) h-cofibration.

\item
\label{item--symseq.symmetric.flat}
The cofibrations of $\Sigma^{\posi} \C$ form a symmetric flat class in $\Sigma \C$:
for any finite family $v$ of positive (acyclic) cofibrations, and any weak equivalence $y \in \Sigma_n \Sigma \C$, $y \pp_{\Sigma_n} v^{\pp n}$ is a weak equivalence.

\item
\label{item--symseq.symmetroidal}
If the model structure on $\Sigma \C$ is \emph{strongly} $\Y$-admissible for some $\Y = (\Y_{n})$ as in \refde{admissible}, then the (acyclic) cofibrations in $\Sigma^{\posi} \C$ form a class that is (acyclic) $\Y$-symmetroidal in the model structure in $\Sigma \C$:
for any finite family $v$ of positive (acyclic) cofibrations, and any map $y \in \Y_n$, $y \pp_{\Sigma_n} v^{\pp n}$ is a cofibration.
In particular, $\Sigma^{\posi} \C$ is $\Y$-symmetroidal in this case.
\end{enumerate}

By \csy{\refth{symmetric.weakly.saturated}}, these statements only have to be checked on generating (acyclic) cofibrations.
The acyclic parts of these statements are proven by replacing the words ``cofibration'' and ``h-cofibration'' by their acyclic analogues, so that proof is omitted.
Let $v = (v_1, \ldots, v_e)$ be a finite family of generating cofibrations of $\Sigma^{\posi} \C$.
They are given by $v_i= G_{t_i}(h_i)$ for some generating cofibrations $h_i \in \Sigma_{t_i} \C$ and $t_i \ge1$.
Let $n = (n_i)$ be a multi-index with $n_i \ge1$.

For an object $Y = G_u (Z)$ in $\Sigma_n \Sigma \C$, where $Z \in \Sigma_n \Sigma_u \C$, we have
$$Y \t_{\Sigma_n} v^{\pp n}=\left (G_{t n+u} (\Sigma_{tn+u} \cdot_{\Sigma_u \x \Sigma_{t}^n} Z \t h^{\pp n}) \right)_{\Sigma_n} = G_{tn+u} (\Sigma_{tn+u} \cdot_{\Sigma_n \rtimes (\Sigma_u \x \Sigma_t^n)} Z \t h^{\pp n}).\eqlabel{proof.Y}$$
The group $\Sigma_n$ acts trivially on $\Sigma_u$ and $\Sigma_u$ acts trivially on $h^{\pp n}$.
In $\C$ (as opposed to $\Sigma_{t n+u} \C$), there is an isomorphism
$$\varphi (\Sigma_{t n+u} \cdot_{\Sigma_n \rtimes (\Sigma_u \x \Sigma_{t}^n)} Z \t h^{\pp n})
={\Sigma_{tn+u}\over\Sigma_n \rtimes (\Sigma_u \x \Sigma_{t}^n)} \cdot \varphi \left (Z \t h^{\pp n} \right),\eqlabel{proof.77}$$
by \csy{\refle{dirty}}.
This uses the positivity of at least one of the $t_i$, which implies that $\Sigma_n \rtimes (\Sigma_u \x \Sigma_t^n) = \Sigma_u \x (\Sigma_n \rtimes \Sigma_t^n)$ is a subgroup of~$\Sigma_{tn+u}$.

\refit{symseq.symmetric.h.monoidal}~We have to show that $Y \t_{\Sigma_n} v^{\pp n}$ is an h-cofibration in $\Sigma \C$ for all $Y \in \Sigma_n \Sigma \C$.
As in the proof of h-monoidality of $\Sigma \C$ (\refth{stable.R.general}), we may assume $Y = G_u(Z)$, where $u\ge0$ and $Z \in \Sigma_m \Sigma_u \C$ is arbitrary.
We show the stronger statement that the above map is an h-cofibration in $\C$ in all degrees.
The $h_i$ are cofibrations, so that $h^{\pp n}$ is also a cofibration (in $\C$, by the pushout product axiom).
Hence, $Z \t h^{\pp n}$ and therefore the right side of \refeq{proof.77} are h-cofibrations in $\C$, using the h-monoidality of $\C$.

\refit{symseq.symmetric.flat}~The symmetric flatness of the positive cofibrations in $\Sigma \C$ is proven similarly: replace $Y \t_{\Sigma_n} - $ by $y \pp_{\Sigma_n} -$ for any weak equivalence $y \in \Sigma_n \Sigma \C$.
Again, the reduction from a general weak equivalence~$y$ to $y = G_u(z)$, $u \ge0$, $z$ a weak equivalence in $\Sigma_n \Sigma_u \C$, is possible by pretty smallness.
Now, note that $z \pp h^{\pp n}$ is a weak equivalence in~$\C$ since $\C$ is flat.

\refit{symseq.symmetroidal}~To show that $y \pp_{\Sigma_n} v^{\pp n}$ is a cofibration in $\Sigma \C$ we may replace $y$ by $G_u(y_u)$, in which case the expression reads
$$\Sigma_{tn+u} \cdot_{\Sigma_n \rtimes (\Sigma_u \x \Sigma_t^n)} z \pp h^{\pp n}.$$
This is a cofibration in $\Sigma_{tn+u} \C$ exactly by the strong admissibility condition \ref{defi--admissible}\refit{admissible.wedge.2}.
The $\Y$-symmetroidality of $\Sigma^{\posi} \C$ also follows from this, noting that $tn + u \ge u \ge1$ in this case, so the previous expression is in addition an isomorphism in degree $0$, hence a positive cofibration.

We have finished the proof of \refit{symseq.symmetric.h.monoidal}--\refit{symseq.symmetroidal}.
The adjunction $\Sigma \C \rightleftarrows \Mod_R$ is a Hopf adjunction with strong monoidal left adjoint, so we can apply
\csy{\refth{transfer.symmetric.monoidal}}.
Therefore, the cofibrations of $\Mod_R^{\posi}$ are symmetric flat and symmetric h-monoidal in $\Mod_R$, and similarly for symmetroidality with respect to the class $R \t \Y$.

By \csy{\refth{symmetric.monoidal.localization}}, the positive cofibrations are also symmetric flat and symmetric h-monoidal in $\Mod_R^{\stab}$.
Since (acyclic) h-cofibrations depend only on the weak equivalences, the symmetric h-monoidality and symmetric flatness of a class of morphisms also depend only on the weak equivalences.
By \refpr{stable.we}, we therefore conclude that the stable positive model structure is symmetric flat and symmetric h-monoidal.

The nonacyclic part of $\ws{R \t \Y}$-symmetroidality of $\Mod_R^{\stab, \posi}$ follows immediately from the one of $\Mod_R^{\posi}$.
The acyclic part follows from a variant of \csy{\refth{symmetric.monoidal.localization}\refit{localization.symmetroidal}}, as follows: by \csy{\refth{symmetric.weakly.saturated}\refit{symmetroidal.weakly.sat}}, it is enough to show that the generating acyclic cofibrations of $\Mod_R^{\stab, \posi}$ are acyclic $\ws{R \t \Y}$-symmetroidal in $\Mod_R^{\stab}$.
By tractability, we may assume they have cofibrant source.
Thus they are acyclic $\Y$-symmetroidal in $\Mod_R^{\stab}$ by \csy{\refpr{acyclic}}.
\xpf

\rema
The symmetroidality of $\Sigma \C$ holds provided that $\C$ itself is symmetroidal.
This excludes the projective model structure on chain complexes of abelian groups, for example.
(See \csy{\refchap{examples}} for a discussion of concrete model categories (not) satisfying symmetric h-monoidality, symmetroidality and symmetric flatness.)
The positive structure $\Sigma^\posi \C$ does not require such an assumption.
Likewise, the (nonsymmetric) flatness of $\C$ promotes to the \emph{symmetric} flatness property \refit{symseq.symmetric.flat} of $\Sigma^{\posi} \C$ above.

The strong admissibility (as opposed to mere admissibility) is necessary to ensure the symmetroidality $\Sigma^\posi \C$.
For example, the argument above fails for the projective structure on~$\Sigma \C$: take $\C = \sSet$ and pick $t=1$, $v=G_1(h)$, where $h$ is some cofibration (i.e., monomorphism) in $\sSet$.
However, $\Sigma_n$ does not usually act freely on the complement of the image $h^{\pp n}$, so this map is not a projective cofibration in $\Sigma_n \sSet$.
\xrema

\numberwithin{equation}{section}
\section{Operadic algebras in symmetric spectra}
\label{sect--algebras.spectra}

We now exploit the excellent model-theoretic properties of the stable positive model structure $\Mod_R^{\stab, \posi}$ on symmetric $R$-spectra to study algebras over operads in this category.
A single-colored operad $O$ in $\Mod_R$ consists of an $R$-module $O_n$ with a $\Sigma_n$-action for each $n \ge0$.
It can be thought of as the space of $n$-ary operations.
For different $n$, they are connected by $\Sigma_{r_1} \x \cdots \x \Sigma_{r_n}$-equivariant multiplication maps (tensor products are over $R$)
$$O_n \t O_{r_1} \t \cdots \t O_{r_n} \r O_{r_1 + \cdots + r_n}.$$
Moreover, there is a unit map $\eta_O\colon R \r O_1$ (which we can rewrite as $\eta_O \colon R[1] \r O$, where $R[1]$ is $R$ in level~1 and $\emptyset$ otherwise).
Multiplication and unit map obey the usual rules.
Specifying an $O$-algebra structure on some $M \in \Mod_R$ amounts to specifying maps
$$O_n \t_{\Sigma_n} M^{\t n} \r M$$
which are again compatible in a suitable sense.
For example, the commutative operad $\Comm$ is such that $\Comm_n=1_{\Mod_R} = R$, so a $\Comm$-algebra is exactly a strictly commutative ring spectrum.

Since there are no essential additional difficulties, we actually work with \emph{$W$-colored operads} or just operads for short.
The set $W$ (called set of \emph{colors}) is fixed.
Instead of the indexing by $n \in \NN$ in single-colored operads, $W$-colored operads are indexed by tuples $(s, w)$ consisting of a map of sets $s\colon I \r W$ (the \emph{multisource}), where $I$ is a finite set and $w \in W$ (the \emph{target}).
Such tuples form a category $\sSeq_W$.
This category is a groupoid and the automorphism group of $(s, w)$ is given by $\Ax_s := \prod_{r \in W} \Sigma_{s^{-1}(r)}$.
The category $\sColl_W(\Mod_R) := \Fun(\sSeq_W, \Mod_R)$ of \emph{symmetric collections} is equipped with the \emph{substitution product}, denoted~$\circ$, which turns this into a monoidal category.
Its monoidal unit $R[1]$ is such that $R[1]_{s,w} = \emptyset$ except for $s\colon I = \{*\} \r W$, $s(*) = w$, in which case it is $R$, the monoidal unit of $\Mod_R$.

A symmetric $W$-colored operad is, by definition, a monoid in $(\sColl_W(\Mod_R), \circ)$.
They form a category denoted $\sOper_W(\Mod_R)$.
The multiplication $O \circ O \r O$ amounts to giving maps
$$O_{s,w} \t \bigotimes_{i \in I} O_{t_i, s(i)} \r O_{\bigcup_{i \in I} t_i, w}.$$
An $O$-algebra consists of $M_w \in \Mod_R$, for every $w \in W$, together with maps
$$O_{s,w} \t \bigotimes_{i \in I} M_{s(i)} \r M_w.$$
Of course, these are subject to appropriate associativity and unitality constraints.
For a slightly less short summary of operads and their algebras, the reader may consult \cop{\refsect{colored.collections}}.

We now turn to the model-theoretic properties of algebras over operads in $R$-spectra.
We show the admissibility of all operads (\ref{theo--operad.admissible.spectra}), give a criterion for (almost) strong  admissibility of levelwise cofibrant operads (\ref{theo--operad.strongly.admissible.spectra}), rectification of algebras over weakly equivalent operads (\ref{theo--rectification.spectra}), and Quillen equivalences of algebras over operads in different categories of spectra (\ref{theo--transport.algebras.spectra}) and finally the special case of $R$-spectra and $S$-spectra, where $R \sim S$ are weakly equivalent (\ref{coro--transport.algebras.spectra}).

The admissibility of operads in symmetric spectra is due to Elmendorf and Mandell for $\C=\Top$ \cite[Theorem~1.3]{ElmendorfMandell:Rings},
and Harper for $\C=\sSet_\bullet$ \cite[Theorem~1.1]{Harper:Symmetric}.
It was generalized by Hornbostel to the category~$\C$ of simplicial presheaves
with the injective model structure and the injective model structure on~$\Sigma\C$ \cite[Theorem~1.3]{Hornbostel:Preorientations}.
In the latter two cases, all objects are cofibrant.
This considerably simplifies the situation because all h-monoidality questions are trivial.
The assumption that every object in~$\C$ is cofibrant excludes the projective model structures on presheaves, which is a main motivating example for~us.
In fact, this paper grew out from an attempt to construct an algebraic cobordism spectrum,
as a \emph{fibrant} commutative ring spectrum.
The fibrancy is necessary to actually compute the homotopy groups of this spectrum
(i.e., the higher algebraic cobordism groups).
For the injective model structure on presheaves the fibrancy condition is practically impossible to check.

Most theorems in this section are immediate consequences of the corresponding abstract results on operadic algebras.
To help the reader, we highlight the main ideas of the proofs of these abstract results in the case of single-colored operads
in Remarks \ref{rema--explain.symmetric.h.monoidal}, \ref{rema--explain.symmetroidal}, \ref{rema--explain.symmetric.flat} and \ref{rema--explain.transport}.

Let us also point out that their proofs require some hard algebra to describe pushouts of operadic algebras
(Elmendorf--Mandell~\cite[\S12]{ElmendorfMandell:Rings}, Fresse \cite[Proposition~18.2.11]{Fresse:Modules}, Harper \cite[Proposition~7.12]{Harper:Symmetric})
and operads (Spitzweck \cite[Proposition~3.5]{Spitzweck:Operads}, Berger--Moerdijk \cite[Lemma~3.1]{BergerMoerdijk:Derived}, \cite[\S5.11]{BergerMoerdijk:Axiomatic}).
In addition, the techniques developed in \cite{PavlovScholbach:Symmetry, PavlovScholbach:Operads} can be seen as a way to conveniently handle model-categorical properties
that allow one to control the homotopical properties of the maps arising in such pushouts and other operadic constructions.

\theo \label{theo--operad.admissible.spectra}
Suppose $(\Sigma \C,R)$ is a spectral context.
Any (symmetric $W$-colored) operad $O$ in~$\Mod_R$ is admissible, i.e., the category of $O$-algebras
carries a model structure that is transferred along the adjunction
$$O \circ - : \Mod_R^{\stab,\posi}\leftrightarrows\Alg_O(\Mod_R) : U.$$
We refer to it as the \emph{stable positive model structure} and denote it by
$\Alg_O^{\stab,\posi}(\Mod_R)$.

For example, for $O = \Comm$, this gives a model structure on strictly commutative ring spectra.
For the operad $\sOpOp_W$ of $W$-colored operads, this gives a model structure on $W$-colored operads in spectra.
\xtheo

\pf
This follows from \cop{\refth{O.Alg}}, whose assumptions are satisfied by \refth{stable.R.general} and \refpr{stable.R.positive}.
\xpf

\rema
\label{rema--explain.symmetric.h.monoidal}
The key point in the proof of \cop{\refth{O.Alg}} is to ensure that a pushout $a\colon A \r A'$ of a map of free $O$-algebras $O \circ x$, where $x$ is an acyclic cofibration of $\Mod_R^{\stab, \posi}$, is again a weak equivalence of $O$-algebras.
The functor $U$ does not preserve pushouts, but instead, as first observed by Elmendorf and Mandell, $U(a)$ is a transfinite composition of maps of the form $\Env(O, A)_n \t_{\Sigma_n} x^{\pp n}$,
where $\Env(O,A)$ is the so-called enveloping operad (whose precise definition does not matter at this point).
The symmetric h-monoidality of $\Mod_R^{\stab, \posi}$ shown in \refpr{stable.R.positive} ensures that this map is a couniversal weak equivalence.

In fact, we only need one part of the condition of symmetric h-monoidality, namely the part about acyclic cofibrations,
and (in the uncolored case) only a single map~$x$, as opposed to finite families.
Finite families in the definition of symmetric h-monoidality ensure that symmetric h-monoidality is stable under composition of maps,
see \csy{\refth{symmetric.weakly.saturated}} and its proof, which (somewhat surprisingly) is the least formal part in showing that symmetric h-monoidality is stable under weak saturation.
This feature is crucial for carrying symmetric h-monoidality through transfer and left Bousfield localizations,
which was heavily used above (see the proof of \refpr{stable.R.positive}).
Similarly, the presence of the nonacyclic part in the definition of symmetric h-monoidality is to ensure the preservation of this property under left Bousfield localizations,
see the proofs of \csy{\refpr{monoidal.localization}\refit{localization.i.monoidal}, \refth{symmetric.monoidal.localization}\refit{localization.symmetric.i.monoidal}}.
See also \refre{explain.multiindices}.
\xrema

\exam
For $\C=\sSet_\bullet$ and $R$ given by $R_n=(S^1)^{\wedge n}$, $\Alg_\Comm(\Mod_R)$ is known as the category of commutative ring spectra (in simplicial sets).
Another example is the case $\C=\sPSh_\bullet(\Sm/S)$ of pointed simplicial presheaves on the site of smooth schemes over some base scheme $S$
and the monoid given by $R_n=(\P^1_S,\infty)^{\wedge n}$.
Any of the standard model structures, for example the projective model structure or any monoidal localization, such as the Nisnevich localization or the Nisnevich-$\A^1$ localization
yields a nice model category (\refde{nice.model.category}).
In this case $\Alg_\Comm(\Mod_R)$ is the category of (strictly) commutative motivic ring $\P^1$-spectra.
\xexam

The next result addresses the strong admissibility of operads \cop{\refde{operad.admissible}}, i.e., the question when the forgetful functor $U$ preserves cofibrant objects.
The main abstract result \cop{\refth{O.strongly.admissible}} states that the forgetful functor $\Alg_O^{\stab, \posi}(\Mod_R) \r \Mod_R^{\stab, \posi}$ preserves cofibrations with cofibrant domains for all operads $O$ whose levels $O_{s,w}$ are of the form $R \t $ some \emph{positively} cofibrant object (with respect to the chosen strongly admissible model structure on $\Sigma \C$).
This result is inapplicable to the commutative operad, for example, since $\Comm_0 = R$ is not \emph{positively} cofibrant.
We therefore provide a statement that requires a nonpositive cofibrancy condition on the levels of $O$, and consequently only obtain that positive cofibrations of $O$-algebras are nonpositive cofibrations of spectra.

For $\C = \sSet_\bullet$ and the injective model structure on $\Sigma \sSet_\bullet$ and the commutative operad, the statement is due to Shipley \cite[Proposition~4.1]{Shipley:Convenient}.
By \refle{injective.admissible}, the injective structure on $\Sigma \sSet_\bullet$ is strongly admissible, so our result generalizes Shipley's.
Recall the notion of a strongly $\Y$-admissible model structure from \refde{admissible} and also the construction of such model structures from \refpr{minimal.admissible}.
For $\C = \sSet$ the strong admissibility of (single colored) operads $O$ whose levels $O_n$ are nonpositively cofibrant was independently  obtained by \cite[Theorem~1.3]{Pereira:Cofibrancy}.

\theo
\label{theo--operad.strongly.admissible.spectra}
Let $(\Sigma \C, R)$ be a spectral context.
Suppose the model structure on $\Sigma \C$ is strongly $\Y$-admissible (\refde{admissible}\refit{admissible.wedge.2}) with respect to some class $\Y = (\Y_n \subset \Mor_{\Sigma_n \Sigma \C})$.
For single-colored operads $O$, suppose that for all $n \ge0$, the $n$th degree of the unit map of $O$,
$$(\eta_O)_n \in \ws{R \t \Y_n}.$$
For example, if $\Y_n = \{ \emptyset \r1\}$ for all $n$, this condition is satisfied for the commutative operad $\Comm$.
If $\Y_n$ is a set of generating $\Sigma_n$-injective cofibrations in $\Sigma_n \Sigma \C$, this condition is satisfied for any operad whose levels $O_n$ are obtained as a transfinite composition of maps of the form $R \t $ some injective cofibration in $\Sigma_n \Sigma \C$.

For arbitrary $W$-colored operads, suppose that for all $(s \colon I \r W, w \in W) \in \sSeq_W$,
$$(\eta_O)_{s,w} \in \ws{R \t \Y_n},$$
where $n$ is the finite multi-index given by $n_r = \sharp s^{-1}(r)$ for $r \in W$.
(Note that only finitely many $r$ appear since $I$ is finite.)

Then the forgetful functor
$$\Mod_R^{\stab} \leftarrow \Alg_O(\Mod_R^{\stab, \posi}) : U$$
preserves cofibrant objects and cofibrations between them.
(Note the $\posi$ at the right, which is not present at the left.)
\xtheo

\pf
By \cop{\refle{preparation.strongly.admissible}}, it is enough to notice that for any finite multi-index $n = (n_r)$, $n_r \ge1$, any multisource $s$ as in the statement, any $w \in W$, and any finite family $x = (x_r)$ of generating cofibrations of $\Mod_R^{\stab, \posi}$
$$(\eta_O)_{s,w} \pp_{\Sigma_n} x^{\pp n} := (\eta_O)_{s,w} \pp_{\prod_r \Sigma_{n_r}} \ppop_r x_r^{\pp n_r}$$
is a cofibration in $\Mod_R^{\stab}$ by \refpr{stable.R.positive}.
\xpf

\rema
\label{rema--explain.symmetroidal}
The main point in the proof of \cop{\refle{preparation.strongly.admissible}} is a description of the levels of cofibrant $O$-algebras:
the underlying $R$-module of a cofibrant $O$-algebra $A$ lies in the weak saturation (saturation under transfinite composition, pushouts, and retracts)
of maps of the form $(\eta_O)_n \pp_{\Sigma_n} x^{\pp n}$.
The strong admissibility of the model structure ensures that these maps are (nonpositive) cofibrations of $R$-modules.
\xrema

The following is a rectification result for algebras over weakly equivalent operads in spectra.
For $\C$ being the category of compactly generated topological spaces, it is due to Goerss and Hopkins \cite[Theorem~1.2.4]{GoerssHopkins:Structured}.
For $R$-spectra in spaces, where $R$ is the free commutative monoid on the monoidal unit~1 in degree~1, this is due to Sagave and Schlichtkrull~\cite[Proposition~9.12]{SagaveSchlichtkrull:Diagram}, see also \refpr{I.spaces}.

\theo
\label{theo--rectification.spectra}
Let $(\Sigma \C, R)$ be a spectral context.
Let $\psi\colon P \r Q$ be a map of operads in $\Mod_R$.
Then there is a Quillen adjunction
$$Q \circ_P - : \Alg_P^{\stab,\posi}(\Mod_R) \rightleftarrows \Alg_Q^{\stab,\posi}(\Mod_R) : U.$$
If $\psi$ is a weak equivalence, i.e., if $P_{s,w} \r Q_{s,w}$ is a weak equivalence in $\Mod_R^{\posi, \stab}$ for all $(s,w) \in \sSeq_W$, this is
a Quillen equivalence.
\xtheo

\exam
For example, there is a Quillen equivalence of algebras over the Barratt--Eccles operad
(i.e., $\Ei$-ring spectra) and commutative monoids in~$\Mod_R$ (i.e., commutative ring spectra).

Another obvious application is that $\Ai$-ring spectra can be rectified to strictly associative ring spectra.
See, e.g., \cop{\refsect{category.theory}} for a definition of $\Ai$.
\xexam

\pf
Again, this follows from \cop{\refth{O.Alg}} and \cop{\refth{cow.Coll}}, whose assumptions are satisfied by \refth{stable.R.general} and \refpr{stable.R.positive}.
\xpf

\rema
\label{rema--explain.symmetric.flat}
The main point in the abstract rectification results (for single-colored operads, in a symmetric monoidal model category $\C$)
is to show that for a free $P$-algebra $A := P \circ X$ on a cofibrant object $X \in \C$, the map $A \r Q \circ_P A = Q \circ X$ is a weak equivalence (in $\C$).
By definition, $\psi \circ X$ is a coproduct of maps of the form $\psi_n \t_{\Sigma_n} X^{\t n}$, which are weak equivalences by the symmetric flatness of $\Mod_R^{\stab, \posi}$.

For single-colored operads, it is therefore enough to know $y \t_{\Sigma_n} V^{\t n}$ is a weak equivalence
for any cofibrant object $V \in \Mod_R^{\stab, \posi}$ and any $\Sigma_n$-equivariant stable weak equivalence,
which is somewhat less demanding than the actual definition of symmetric flatness,
which requires that $y \pp_{\prod \Sigma_{n_i}} (v^{\pp n_1} \pp \dots \pp v^{\pp n_e})$ is a weak equivalence
for finite families of cofibrations $v = (v_1, \dots, v_e)$ and $\prod_{i=1}^e \Sigma_{n_i}$-equivariant weak equivalences $y$.
In the same vein as explained in \refre{explain.symmetric.h.monoidal}, the point in considering multi-indices in the definition of symmetric flatness
is to ensure that the condition is stable under composition (of the maps~$v$),
which is in its turn the crucial technical point in showing that symmetric flatness is preserved under left Bousfield localizations and transfers.
In other words, if we were to replace the statement of \refit{symseq.symmetric.flat} in the proof of \refpr{stable.R.positive} by the simpler,
but weaker statement that $y \t_{\Sigma_n} V^{\t n}$ is a weak equivalence for any cofibrant object $V \in \Sigma^\posi \C$ and any weak equivalence $y$,
we would be unable to promote this property to $\Mod_R^\posi$ and then to $\Mod_R^{\posi, \stab}$.
\xrema

A related result addresses the following question: when does the quasicategorical notion of an algebra over an operad in spectra
coincide (up to an equivalence) with the usual strict notion?

\theo
\label{theo--quasirect.spectra}
\def\HAlg{\mathop{\bf HAlg}\nolimits}
Let $(\Sigma \C, R)$ be a spectral context.
Let $O$ be an operad in simplicial sets.
Denote by $\CO_{\Mod_R}$ and $\CO_{\Alg_O(\Mod_R)}$ the full subcategories spanned by the corresponding classes of cofibrant objects, where in this theorem $\Mod_R$ is a shorthand for $\Mod_R^{\stab, \posi}$.
The canonical comparison functor $$\N(\CO_{\Alg_O(\Mod_R)})[\we_{\Alg_O(\Mod_R)}^{-1}]\to\HAlg_{\N^\otimes O}(\N(\CO_{\Mod_R})[\we_{\Mod_R}^{-1}])$$
is an equivalence of quasicategories.
Here $\HAlg$ is used in the sense of \cite[Definition~2.1.3.1]{Lurie:HA} (denoted $\Alg$ there)
and $\N^\otimes O$ denote the operadic nerve of~$O$, as explained in Definition~2.1.1.23 there.
\xtheo

\pf
This follows from the symmetric flatness of $\Mod_R$ (\refpr{stable.R.positive}) and \cop{\refth{quasicategorical.rectification}}.
\xpf

We finally give two transport results that describe the category of operadic algebras in different categories of spectra.
The first result is about a general weak monoidal Quillen adjunction.
In the special case of algebras in $R$-spectra and $S$-spectra, where $R \sim S$ are weakly equivalent commutative monoids in $\Sigma \C$, we get a stronger result.

Let
$$F : \C \leftrightarrows \D : G
\eqlabel{adjunction.F.G}$$
be an adjunction between closed symmetric monoidal categories where $G$ is symmetric lax monoidal.
We pick commutative monoids $S \in \Sigma \D$ and $R \in \Sigma \C$ and a map of commutative monoids $\varphi\colon R \r G(S)$.
Note that $G$ preserves commutative monoids since it is symmetric lax monoidal.
There are adjunctions
$$\eqalignno{F^\Mod : \Mod_R^{\C} & \leftrightarrows \Mod_S^{\D} : G,&\eqlabel{adjunction.F.G.Mod}\cr
(F^\Mod)^\sOper : \sOper_W (\Mod_R^{\C}) & \leftrightarrows \sOper_W(\Mod_S^{\D}) : G,\cr}$$
where $G$ is in both cases the obvious functor and $F^\Mod$ and $(F^\Mod)^{\sOper}$ are constructed in, say, \cite[Theorem~4.5.6]{Borceux:2}.

From now on, assume that $\C$ and $\D$ are nice (\refde{nice.model.category}).
To avoid technical statements around the cofibrancy hypothesis on $R$, we also assume $\C$ is combinatorial or $\C = \Top$ (see \refre{Top}) and likewise for $\D$.
We equip $\Sigma \C$ and $\Sigma \D$ with some admissible model structures and we consider the condition that this datum induces a weak monoidal Quillen adjunction \cite[Definition~3.6]{SchwedeShipley:Equivalences}
$$F^\Mod : \Mod_R^{\stab, \posi, \C} \leftrightarrows \Mod_S^{\stab, \posi, \D} : G,\eqlabel{adjunction.number.42}$$
that is
$$\eqalign{F^\Mod(\mcr R) & \r S, \cr
F^\Mod(C \t_R C') & \r F^\Mod(C) \t_S F^\Mod(C').\cr}$$
are weak equivalences for all cofibrant objects $C, C' \in \Mod_R^{\stab, \posi, \C}$.
Using the Quillen equivalences $\Mod_R^{\stab, \posi, \C} \sim \Mod_R^{\stab, \C}$ (\refpr{stable.we}), this condition is equivalent for the nonpositive or the positive stable model structures.
Since $F^\Mod(R \t -) = S \t F(-)$, the first condition holds if $1 \in \C$ is cofibrant.
Using pretty smallness (via \csy{\refle{sequential}}), the second condition can be reduced to free $R$-modules $C$ and $C'$, so that it holds provided that the original adjunction \refeq{adjunction.F.G} is weakly monoidal and that $\Sigma \C$ and $\Sigma \D$ both carry the projective model structure.

\theo \label{theo--transport.algebras.spectra}
Fix two spectral contexts $(\Sigma \C, R)$ and $(\Sigma \D, S)$ satisfying the above assumptions (in particular, \refeq{adjunction.number.42} is supposed to be a weak monoidal Quillen adjunction).
Then, for any operad $O$ in $\Mod_R$ and any operad $P$ in $\Mod_S$, there are Quillen adjunctions
$$\eqalign{F^\Alg : \Alg_O^{\stab, \posi} (\Mod_R^\C) & \rightleftarrows \Alg_{F^\sOper(O)}^{\stab, \posi} (\Mod_S^\D) : G, \cr
F_\Alg : \Alg_{G(P)}^{\stab, \posi} (\Mod_R^\C) & \rightleftarrows \Alg_{P}^{\stab, \posi} (\Mod_S^\D) : G.\cr}$$
They are Quillen equivalences if $(F^\Mod, G)$ is a weak monoidal Quillen equivalence and $O$ is cofibrant and $P$ is fibrant.
\xtheo

\pf
This is an immediate application of \cop{\refth{transport}} whose assumptions are satisfied by \refth{stable.R.general}.
\xpf

\rema
\label{rema--explain.transport}
The abstract transport result \cop{\refth{transport}} lies considerably deeper than the (strong) admissibility and rectification results used above.
The key point in the transport for cofibrant operads~$O$ is to use Spitzweck's and Berger--Moerdijk's analysis of the pushout of operads
and tracking the homotopy-theoretic properties of the maps arising in this description relative to a weak monoidal Quillen adjunction.
The statement for fibrant operads~$P$ uses in addition the rectification result \ref{theo--rectification.spectra}.
\xrema

In the special case $\C = \D$ and a weak equivalence $\varphi\colon R \r S$ in $\Sigma \C$, the transport of algebras applies to more general operads:

\coro
\label{coro--transport.algebras.spectra}
Suppose that for a spectral context $(\Sigma \C, R)$ and a weakly equivalent commutative monoid $S \in \Sigma \C$, there are Quillen equivalences
$$\varphi_* : \Mod_R^{\stab, \posi} \leftrightarrows \Mod_S^{\stab, \posi} : \varphi^*.\eqlabel{adjunction.Mod.R.S}$$
(See \refpr{stable.R.transport.2} for sufficient criteria.)
Then there are Quillen equivalences
$$\eqalignno{\varphi_* : \Alg_O^{\stab, \posi} (\Mod_R) & \leftrightarrows \Alg_{S \t_R O}^{\stab, \posi} (\Mod_S) : \varphi^*,&\eqlabel{adjunction.Alg.R.S}\cr
\varphi_* : \Alg_{\varphi^* P}^{\stab, \posi} (\Mod_R) & \leftrightarrows \Alg_P^{\stab, \posi} (\Mod_S) : \varphi^*&\eqlabel{adjunction.Alg.R.S.2}\cr}$$
for any operad $O$ in $\Mod_R$ whose levels $O_{s,w}$ are cofibrant in $\Mod_R^{\stab}$ and any operad $P$ in $\Mod_S$ whose levels $P_{s,w}$ are fibrant in $\Mod_S^{\stab, \posi}$.
\xcoro

\exam
\label{exam--transport.commutative.ring.spectra}
If $1_\C$ is cofibrant, $R=1_{\Mod_R}$ is cofibrant in $\Mod_R^{\stab}$.
The levels of the commutative operad $O=\Comm$ are given by $O_n=1_{\Mod_R} = R$.
We get $S \t_R O = \Comm$ and therefore a Quillen equivalence of commutative ring spectra.
\xexam

\pf
The left adjoint in \refeq{adjunction.Mod.R.S} is strong symmetric monoidal, the right adjoint is lax monoidal.
It therefore gives an adjunction whose left adjoint is again strong monoidal.
$$\varphi_* : (\sColl (\Mod_R), \circ) \leftrightarrows (\sColl (\Mod_S), \circ) : \varphi^*.$$
Therefore, by \cite[Proposition I.3.91]{AguiarMahajan:Monoidal}, both right and left adjoints preserve commutative monoids, i.e., operads, and the induced functors on operadic algebras form an adjoint pair.
In other words, on the underlying spectra, $(\varphi_*)^\sOper$ is just $\varphi_*$.
In the same vein, $(\varphi_*)^{\Alg}$ and $(\varphi_*)_\Alg$ are also just given by $\varphi_*$ on the underlying level.

Applying \refth{transport.algebras.spectra} to $\mcr^\sOper(O)$ (the cofibrant replacement of $O$, using the model structure on operads in $\Mod_R^{\stab, \posi}$ established in \refth{operad.admissible.spectra}) and $\mfr P$ (the fibrant replacement) and rectification (\refth{rectification.spectra}) we have to show weak equivalences of operads
$$\eqalign{\varphi_* (\mcr^\sOper(O)) & \sim \varphi_* (O) \cr
\varphi^* (P) & \sim \varphi^* (\mfr P).\cr}$$
The latter holds since $\mfr P \r P$ is a weak equivalence of operads whose levels $(\mfr P)_{s,w} \r P_{s,w}$ are a weak equivalence between fibrant objects in $\Mod_S^{\stab, \posi}$ (the latter by assumption).
Being a right Quillen functor, $\varphi^*$ preserves this weak equivalence by Brown's lemma.

The former weak equivalence is shown as follows: the weak equivalence $\mcr^\sOper(O) \r O$ gives a weak equivalence of the levels $\mcr^\sOper(O)_{s,w} \sim O_{s,w}$.
The operad $\mcr^\sOper(O)$ is cofibrant, and therefore its unit map $\eta \colon R[1] \r \mcr^\sOper(O)$ is a cofibration of symmetric collections.
For unit degrees, $(s,w) = (w,w)$, the unit map $R=1_{\Mod_R} \r \mcr^\sOper(O)_{w,w}$ is therefore a cofibration, whereas in degrees $(s,w) \ne (w,w)$, $\mcr^\sOper(O)_{s,w}$ is cofibrant in $\Mod_R^{\stab, \posi}$ and a fortiori in $\Mod_R^{\stab}$.
By \cop{\refle{monoidally.cofibrant.flat}}, using the flatness of $\Mod_R^{\stab}$, $\mcr^\sOper(O)_{s,w} \t_R - $ preserves stable weak equivalences in both cases.
Similarly for $O_{s,w}$, using the cofibrancy assumption on $O_{s,w}$.
Hence we get a chain of weak equivalences in $\Mod_R^{\stab}$ or equivalently in $\Mod_R^{\stab, \posi}$:
$$\mcr^\sOper(O)_{s,w} \t_R S \sim \mcr^\sOper(O)_{s,w} \sim O_{s,w} \sim O_{s,w} \t_R S.$$
\xpf

\section{Derived algebraic geometry over spectra}
\label{sect--DAG}

In this section, we prove that symmetric spectra in a symmetric monoidal model category form a \emph{homotopical algebraic context}
in the sense of To\"en and Vezzosi~\cite{ToenVezzosi:HomotopicalII}, so that one can do derived algebraic geometry over ring spectra.

\defi
A \emph{homotopical algebraic context} is a model category~$\D$ such that the following conditions are satisfied.
\begin{enumerate}[(i)]
\item \label{item--hac.stable}
$\D$ is a proper, pointed, combinatorial symmetric monoidal model category.
The canonical morphism from the homotopy coproduct to the homotopy product of any finite family of objects is a weak equivalence.
The homotopy category of~$\D$ is additive.
\item \label{item--hac.comm}
For any commutative monoid~$P$ in~$\D$
the transferred model structure on $\Mod_P(\D)$ exists and is a proper, flat, combinatorial symmetric monoidal model category.
\item \label{item--hac.P.Alg}
The transferred structure on commutative $P$-algebras and commutative nonunital $P$-algebras exists and is a proper combinatorial model category.
\item \label{item--hac.P.Alg.2}
Given a weak equivalence $f\colon E \r F$ in $\Mod_P(\D)$ and a cofibrant commutative $P$-algebra $Q$, $Q \t_P f$ is a weak equivalence in $\Mod_Q(\D)$.
\end{enumerate}
\xdefi

\rema
To\"en and Vezzosi also list Assumption~1.1.0.6 among the properties of a homotopical algebraic context,
which demands a t-structure whose positive part (i.e., connective objects) is closed under the derived tensor product.
Given a homotopical algebraic context in the sense of the previous definition,
such a t-structure can be constructed from \emph{any} set of objects by closing it under derived tensor products and homotopy colimits.
Lurie \cite[Proposition~1.4.4.11]{Lurie:HA} proves the axioms of a t-structure.
\xrema

\theo
\label{theo--Toen.Vezzosi}
Suppose $\C$ is a pointed, nice model category (\refde{nice.model.category}).
We fix an admissible model structure on $\Sigma \C$ and consider a commutative monoid~$R \in \Sigma \C$ such that $R_1$ is weakly equivalent to $S^1 \t B$,
where $S^1$ is a cofibrant representative of $* \sqcup_{1_\C}^h *$, the suspension of the monoidal unit, and $B \in \C$ is any cofibrant object.
Then the stable positive model structure~$\D := \Mod_R^{\stab,\posi}$ on the category of symmetric $R$-spectra defined in \refth{stable.R.general}
is a homotopical algebraic context, except possibly for the properness of the model categories.
\xtheo

\pf
\refit{hac.stable}~This is a restatement of \refth{stable.R.general}.
The last statement follows from the stability of $\D$, which holds by the assumption on $R_1$ and \refpr{stable.R.stable}.

\refit{hac.comm}~Let $P \in \Comm(\D)$, i.e., $P$ is a commutative ring spectrum.
The model structure on $\D$ transfers to a combinatorial, left proper, symmetric monoidal model structure on $\Mod_P$ by \csy{\refth{commutative.monoid}},
using that $\D$ satisfies the monoid axiom by \refpr{symseq.general}.
Likewise, the flatness of $\D$ transfers to $\Mod_P$ by \csy{\refpr{transfer.monoidal}\refit{transfer.flat}}.

\refit{hac.P.Alg}~The categories of (nonunital) commutative $P$-algebras are algebras over the operad $\Comm$ and $\Comm^+$ (which is given by $\Comm^+_n = \emptyset$ for $n=0$ and the monoidal unit $1$ for $n>0$), with values in $\Mod_P$.
Again by \csy{\refth{commutative.monoid}}, $\Mod_P$ is symmetric h-monoidal, so that any operad in $\Mod_P$, in particular $\Comm$ and $\Comm^+$ are admissible, so the transferred model structure on (nonunital) commutative $P$-algebras exists \cop{\refth{O.Alg}}.

\refit{hac.P.Alg.2}~As usual, we prove this by cellular induction.
The first case is when $Q = \Sym (P \t X)$, where $X$ is the (co)domain of a generating cofibration of $\D$ and $\Sym$ denotes the symmetric algebra on the $P$-module $P \t X$.
As above, we have a canonical isomorphism in $\C$:
$$Q \t_P f = \coprod_{t \ge0} ((P \t X)^{\t_P n})_{\Sigma_t} \t_P f=
\coprod_{t \ge0} f \t_{\Sigma_t} X^{\t t},$$
where $\Sigma_t$ acts trivially on $f$.
This is a weak equivalence in $\D$ since $\D$ is symmetric flat.
As $\D$ is h-monoidal, weak equivalences are closed under finite coproducts \cite[Proposition~1.15]{BataninBerger:Homotopy} and therefore, using the pretty smallness of $\D$, closed under countable coproducts.

Next, consider a cocartesian square in $\Alg_P$, where $i\colon X \r X'$ is a generating cofibration in $\C$,
$$\xymatrix{
\Sym(P \t X) \ar[r] \ar[d] &
\Sym(P \t X') \ar[d] \cr
Q \ar[r] &
Q',
}
\eqlabel{pf.99}$$
we want to show that our claim is true for $Q'$, provided that it holds for $Q$.
We again use the filtration that already appeared in the proof of \cop{\refth{O.Alg}}.
In the case considered here, $O = \Comm$, so that $\Env(O,Q)_t = Q$ (with the trivial $\Sigma_t$-action).
This description of the enveloping operad can be read off its explicit description in \cite[Proposition~7.6]{Harper:Symmetric}
(in loc.~cit., $\Env(O,Q)_t$ is denoted by $O_Q[t]$, and the formula for $O_Q[t]$ simplifies to $O_Q[t] = \colim (O \circ A \leftleftarrows O \circ (O \circ A))$ for $O = \Comm$).
As in \cop{\refth{O.Alg}}, we get a cocartesian square in $\Mod_P$,
$$\xymatrix{
Q \t (\ppdom^t i)_{\Sigma_t} = Q \t_P (\ppdom^t_P (P \t i))_{\Sigma_t} \ar[r] \ar[d] &
Q \t (X'^{\t t})_{\Sigma_t} = Q \t_P ((P \t X')^{\t_P t})_{\Sigma_t} \ar[d] \cr
Q_t \ar[r] &
Q_{t+1}.
}
$$
\vadjust{\filbreak}
We apply $f \t_P -$ to this square and get a cube whose front and back faces are cocartesian (in $\Mod_Q$, or in~$\C$):
$$\xymatrix{
& (F \t_P Q) \t_{\Sigma_t} \ppdom^t i \ar[rr] \ar'[d][dd] &
& (F \t_P Q) \t_{\Sigma_t} X'^t \ar[dd] \cr
(E \t_P Q) \t_{\Sigma_t} \ppdom^t i \ar[rr] \ar[ur]^\sim \ar[dd] &
& (E \t_P Q) \t_{\Sigma_t} X'^t \ar[ur]^\sim \ar[dd] \cr
& F \t_P Q_t \ar'[r][rr] &
& F \t_P Q_{t+1} & \cr
E \t_P Q_t \ar[rr] \ar[ur]^\sim &
& E \t_P Q_{t+1} \ar[ur] & \cr
}$$
The top horizontal arrows of this cube are h-cofibrations (in $\C$, say), since $i$ is a symmetric h-cofibration.
Consequently the front and back faces are homotopy pushout squares.
The three arrows labeled with $\sim$ are weak equivalences by induction and the case of free commutative $P$-algebras considered above.
Therefore, the map $f \t_P Q_{t+1}$ is also a weak equivalence.

Now, any cofibrant $P$-algebra is a retract of transfinite compositions (in $\Alg_P$) of maps as in \refeq{pf.99}.
The forgetful functor $\Alg_P \r  \Mod_P$ commutes with sifted colimits, therefore with transfinite compositions (and retracts).
Weak equivalences in $\C$ are stable under filtered colimits by \csy{\refle{sequential}}.
This finishes the proof of \refit{hac.P.Alg.2}.
\xpf

\section{Goerss--Hopkins obstruction theory for structured ring spectra}

In \cite{GoerssHopkins:Commutative}~and~\cite{GoerssHopkins:Structured}, Goerss and Hopkins formulated a number of axioms
that a category of spectra should satisfy in order to admit a good obstruction theory for lifting commutative monoid objects in the homotopy category of spectra to $\Ei$-spectra.
They pointed out that the stable positive model structure on topological spectra satisfies these properties
and raised the question whether the same property is true for spectra in a general model category.
This was shown for spectra with values in simplicial presheaves by Hornbostel \cite[\S3.3]{Hornbostel:Preorientations}.
In this section, we answer this question in the positive for spectra in the very broad class of {\it nice\/} model categories (\refde{nice.model.category}).

We summarize the more recent version of the axioms in~\cite{GoerssHopkins:Structured} in the following definition:

\defi \label{defi--Goerss.Hopkins}
A \emph{Goerss--Hopkins context} is a symmetric monoidal quasitractable combinatorial (or cellular) stable $V$-enriched model category~$\C$
($V$ is a quasitractable combinatorial (or cellular) symmetric monoidal model category)
such that every operad~$\cO$ in~$\C$ is admissible with the resulting model structure on $\cO$-algebras being quasitractable combinatorial (or cellular) and $V$-enriched
and every weak equivalence of operads induces a Quillen equivalence between their categories of algebras.
\xdefi

\theo
\label{theo--Goerss.Hopkins}
Suppose $(\Sigma \C, R)$ is a spectral context (\refde{context}) where $\C$ is in addition pointed and $\V$-enriched,
and the first level $R_1$ of $R$ is weakly equivalent to $S^1 \t B$, where $S^1$ is a cofibrant representative of the suspension of the monoidal unit
and $B$ is any cofibrant object.
The category of $R$-spectra, equipped with the stable positive model structure established in \refth{stable.R.general}, is a Goerss--Hopkins context.
\xtheo

\pf
The model structure $\Mod_R^{\stab, \posi}$ is stable, symmetric monoidal and tractable by \refth{stable.R.general}, \refpr{stable.R.stable}, and \refpr{symseq.general}
Every operad $\cO$ in $\Mod_R$ is admissible by \refth{operad.admissible.spectra}, and weak equivalences of operads induce Quillen equivalent categories of algebras by \refth{rectification.spectra}.
\xpf

\refde{Goerss.Hopkins} is slightly different from the list of properties mentioned in \cite[\S1.1, \S1.4]{GoerssHopkins:Commutative} and
\cite[Theorems 1.2.1 and 1.2.3]{GoerssHopkins:Structured}:
we omit the requirement that the homotopy category of~$\C$ is equivalent to the homotopy category of Bousfield--Friedlander spectra, i.e., nonsymmetric spectra.
The Quillen equivalence of symmetric and nonsymmetric spectra with values in an abstract model category is addressed by \cite[Corollary~10.4]{Hovey:Spectra}.
Axiom~1.2.1.3 is automatically satisfied in our setup.

\cite[Axiom~1.2.3.5]{GoerssHopkins:Structured} can be rephrased by requiring that the forgetful functor $\Alg_\cO (\Mod_R) \r \Mod_R$ preserves cofibrations.
In op.~cit.\ this is only used in Theorem~1.3.4.2, which in its turn is only used in Theorem~1.4.9 to establish cellularity, which can be replaced by combinatoriality.
Moreover, this property may fail for internal operads if, say, $\cO(1) \in \Mod_R$ is not cofibrant, so it is omitted in \refde{Goerss.Hopkins}.
A~positive result in this direction, for a general model category~$\C$, is given by \refth{operad.strongly.admissible.spectra}.

\cite [Axiom~1.2.3.6]{GoerssHopkins:Structured} states that for any~$n\ge0$, and any cofibrant object $X \in \Mod_R$,
the functor $\Sigma_n^\inje \V \r \Mod_R$, $K \mapsto K \t_{\Sigma_n} X^{\t n}$ preserves weak equivalences and cofibrations.
The cofibration part of this condition is again not present in \refde{Goerss.Hopkins}.
It is used only in \cite[Theorem~1.2.4]{GoerssHopkins:Structured} (rectification for operads in~$R$-spectra).
Our proof of this statement is based on the symmetric flatness of the stable positive model structure on $\Mod_R$,
which is a generalization of the preservation of weak equivalences by the above functor.

\section{Construction of commutative ring spectra}
\label{sect--construction}

In this section, we apply the results of \refsect{symmetric.spectra} to the construction of strictly symmetric ring spectra representing a given cohomology theory.
We will apply this in \refsect{Deligne} to the case of Deligne cohomology.

We recall two technical tools: first, we study nonsymmetric lax monoidal right adjoints, such as the Dold--Kan functor $\Gamma\colon \cCh \r \simpl\Ab$, and the endomorphism operad associated to such a functor.
This is due to Richter \cite[Definition~3.1]{Richter:Symmetry} (also see \cite[\S4.3.2]{AguiarMahajan:Monoidal}).
Secondly, in order to capture the maximal information from the ring spectra constructed in \refth{construction.spectra}, we will not only consider mapping spaces, but convolution algebras, which encode the multiplication on mapping spaces (see for example \cite[\S3.4.5]{AguiarMahajan:Monoidal}).

\defi
Let $G\colon \D \r \C$ be a lax monoidal (but not necessarily symmetric lax monoidal) functor between two symmetric monoidal categories, where $\C$ is enriched over a symmetric monoidal category~$\V$.
The \emph{endomorphism operad} of~$G$ is the operad in~$V$ defined by
$$\cO_G(n) = \Hom_{\Fun(\D^n,\C)}(G(-) \t \cdots \t G(-), G(- \t \cdots \t -)).$$
We say that $G$ is $\cO$-lax monoidal for some operad $\cO$ in $V$ if there is a natural map $\cO \r \cO_G$.
\xdefi

For example, a symmetric lax monoidal functor $G$ is just the same as a $\Comm$-lax monoidal functor \cite[Table~4.2]{AguiarMahajan:Monoidal}.

\lemm \label{lemm--adjunction.Alg}
Let
$$F : \C \rightleftarrows \D : G\eqlabel{adjunction.F.G.spectra}$$
be an adjunction of symmetric monoidal categories, where $G$ is $\cO$-lax monoidal for some operad $\cO$.
Also suppose that $\C$ and $\D$ are accessible.
\begin{enumerate}
\item
There is an adjunction
$$F^\Alg : \Alg_{\cO} \C \rightleftarrows \Alg_\Comm \D : G,\eqlabel{adjunction.alg}$$
where $G$ sends a commutative algebra $D \in \D$ to $G(D)$ with the $\cO$-algebra structure defined by
$$\cO(n) \t G(D)^{\t n} \r \cO_G(n) \t G(D)^{\t n} \r G(D^{\t n}) \r G(D).$$
\item \label{item--F.strong}
\cite[Proposition~3.91]{AguiarMahajan:Monoidal}
If $G$ is symmetric monoidal (so that $\cO = \Comm$) and $F$ is strong symmetric monoidal, then $F^\Alg$ sends a commutative algebra $C \in \C$ to $F(C)$ with the commutative algebra structure
$$F(C) \t F(C) \stackrel \cong \r F(C \t C) \r F(C),$$
where the first map is the isomorphism that is part of the strong symmetric monoidal functor.
\end{enumerate}
\xlemm
\pf
The functor~$G$ preserves limits and filtered colimits of algebras,
since these are created by the functor forgetting the algebra structure.
Since $G$ is a functor between locally presentable categories, it therefore has a left adjoint $F^\Alg$.
\xpf

\defi
Suppose that $\C$ is a closed symmetric monoidal category.
The internal hom functor $\C^\op \x \C \r \C$ is symmetric monoidal.
The induced functor
$$\Hom^\Alg\colon \Alg_\Comm (\C^\op) \x \Alg_\Comm(\C) = \Alg_\Comm (\C^\op \x \C) \r \Alg_\Comm (\C)$$
is called the \emph{convolution algebra}.
More generally, given an operad $\cO$ in $\C$, the \emph{convolution $\cO$-algebra} is the functor
$$\Conv \colon \Alg_\Comm(\C^\op) \x \Alg_\cO(\C) \r \Alg_\cO(\C),\eqlabel{Conv}$$
which sends $(X, Y)$ to the internal $\Hom(X, Y) \in \C$ equipped with the $\cO$-algebra structure induced by the comultiplication on~$X$ and the $\cO$-algebra structure on~$Y$.
Explicitly, it is defined by
$$\cO(n) \t \Hom(X, Y)^{\t n} \r \cO(n) \t \Hom(X^{\t n}, Y^{\t n}) \r \Hom(X^{\t n}, \cO(n) \t Y^{\t n}) \r \Hom(X, Y).$$
\xdefi

\lemm \label{lemm--monoidal} (cf.\ \cite[3.83]{AguiarMahajan:Monoidal})
In the situation of \refle{adjunction.Alg}, let $\C' \subset \C$ be a full subcategory such that $F' := F|_{\C'}$ is symmetric oplax monoidal (so that $F'$ preserves commutative coalgebras).
The natural transformation
$$\Conv_\C (-, G (-)) \r G (\Conv_\D(F'(-), D))\eqlabel{Conv.F.G}$$
is a morphism of functors $\Alg_\Comm(\C'^\op) \x \Alg_\Comm(\D)  \r \Alg_{\cO_G}(\C)$.
It is an isomorphism if the oplax structural map
$$F(T \t X) \r F(T) \t F(X),\eqlabel{left.adjoint.monoidal}$$
is an isomorphism for any $T \in \C$ and any $X \in \C'$.
\xlemm
\pf
The underlying internal hom objects are given by the compositions
$$\Phi = \Hom(-, G-) \colon \C'^\op \x \D \stackrel{\id \x G}\lr \C'^\op \x \C \stackrel {\Hom_\C} \lr \C,$$
$$\Psi = G \Hom(F'-, -) \colon \C'^\op \x \D \stackrel{F' \x \id}\lr \D^\op \x \D \stackrel {\Hom_\D} \lr \D \stackrel G \lr \C.$$
The functors $G$ and $\id \x G$ are $\cO_G$-monoidal, and all other functors are symmetric lax monoidal, i.e., $\Comm$-monoidal.
Thus, their composition is $\cO_G \boxtimes \Comm = \cO_G$-monoidal.
Here $- \boxtimes -$ denotes the Hadamard product of operads \cite[Theorem~4.28]{AguiarMahajan:Monoidal}.
The natural transformation $\Phi \r \Psi$ induces the transformation in \refeq{Conv.F.G}, which is therefore a map of $\cO_G$-algebras.
For the second claim, $\Phi \r \Psi$ is an isomorphism in this case, hence so is the transformation in \refeq{Conv.F.G}.
\xpf

We now consider the interaction of $\Conv$ and model structures.
Suppose $\C$ is a symmetric monoidal model category.
Then the convolution algebra \refeq{Conv} is a functor between categories with weak equivalences.
To get homotopically meaningful information, we therefore have to derive it.
A natural strategy to compute this (right) derived functor would be to endow the category of commutative coalgebras in $\C$ (=commutative algebras in $\C^\op$) with a model structure.
The standard choice of such a model structure is the transferred structure along the forgetful functor
$$\C^\op \leftarrow \Alg_\Comm(\C^\op).$$
However, this is a notoriously difficult task (see, e.g., \cite{BayehHessKarpovaKedziorekRiehlShipley:Left}), which we will not undertake in this paper.
Instead we use the following fact:

\lemm \label{lemm--convolution.algebras}
Let $\C$ be a symmetric monoidal model category.
Let $X \in \C$ be a cofibrant object that is also endowed with a commutative coalgebra structure.
The functor
$$\Conv (X, -) \colon \Alg_\cO(\C) \r \Alg_\cO(\C)$$
is a right Quillen functor.
Its derived functor will be denoted by $\hConv(X, -)$.
\xlemm
\pf
We have to check $\Conv (X, -)$ preserves (acyclic) fibrations.
These are created by the forgetful functor to $\C$.
Forgetting the $\cO$-algebra structure, $\Conv(X, -)$ is just the internal $\Hom(X, -)$, which is a right Quillen functor since $X$ is cofibrant and $\C$ is a monoidal model category.
\xpf

We now upgrade \refle{monoidal} to model categories.
We use the notation of \refle{adjunction.Alg} and \refle{monoidal}.
\prop \label{prop--monoidal}
Suppose that \refeq{adjunction.F.G.spectra} is a Quillen adjunction between combinatorial model categories and the transferred model structures on the categories of algebras in \refeq{adjunction.alg} exist.
Also suppose $X$ is an object of $\C'$, which is cofibrant in $\C$ and such that the lax monoidal structural map \refeq{left.adjoint.monoidal} is a weak equivalence for all cofibrant objects $T \in \C$.
\begin{enumerate}[(1)]
\item
\label{item--alg.adjunction}
The adjunction \refeq{adjunction.alg}, which exists by \refle{adjunction.Alg}, is a Quillen adjunction.
The map
$$\hConv_\C (X, \rdf G D) \stackrel \sim \lr \rdf G (\hConv_\D(F(X), D))\eqlabel{hConv}$$
is a weak equivalence in $\Alg_\cO \C$.
\item \label{item--alg.3}
In the situation of \refle{adjunction.Alg}\refit{F.strong}, suppose that \refeq{adjunction.F.G.spectra} and \refeq{adjunction.alg} are Quillen equivalences.
Then, for any object $C \in \Alg_\Comm \C$ there is a weak equivalence in $\Alg_\Comm \C$
$$\ldf F^\Alg \hConv_\C(X, C) \stackrel \sim \lr \hConv_\D(F^\Alg  X, \ldf F^\Alg  C).$$
\end{enumerate}
\xprop

\pf
\refit{alg.adjunction}~\refeq{adjunction.alg} is a Quillen adjunction since (acyclic) fibrations are created by the functors forgetting the respective operadic algebra structures.

By \refle{monoidal}, \refeq{hConv} is a map of $\cO$-algebras.
It is therefore enough to show that \refeq{hConv} is a weak equivalence in $\C$, i.e., after forgetting the $\cO$-algebra structure.
This is an easy consequence of the assumption that \refeq{left.adjoint.monoidal} is a weak equivalence.

\refit{alg.3}~In \refeq{hConv}, put $D = \ldf F^\Alg (C)$.
As \refeq{adjunction.alg} is a Quillen equivalence, there is a weak equivalence $C \r \rdf G (\ldf F^\Alg (C))$.
Hence, we get a weak equivalence $\hConv_\C(X, C) \stackrel \sim \lr \rdf G \hConv_\D(FX, \ldf F^\Alg C)$, which implies our claim again using the Quillen equivalence \refeq{adjunction.alg}.
\xpf

We now prepare for \refth{construction.spectra} by fixing some notation related to the Dold--Kan equivalence.
Let $\cA$ be a symmetric monoidal Grothendieck abelian category.
We fix a model structure on the category $\simpl \cA$ of simplicial objects.
We assume that this model structure transfers, via the Dold--Kan equivalence,
$$N : \simpl \cA \rightleftarrows \cCh \cA : \Gamma\eqlabel{Dold.Kan}$$
to a model structure on connective chain complexes.
Finally, we pick a model structure on chain complexes, such that the adjunction of the good truncation functor and the inclusion of chain complexes in nonnegative degrees,
$$\iota : \cCh \cA \rightleftarrows \Ch \cA : \tau\eqlabel{iota.tau}$$
(with $\tau A_* := [\cdots \r A_1 \r \ker d_0]$) is a Quillen adjunction.
We assume that these three model categories are nice (\refde{nice.model.category}),
i.e., they are pretty small, tractable, flat, h-monoidal, symmetric monoidal model categories.
We also assume that the monoidal unit is cofibrant in $\cCh \cA$.
(This is needed to apply \refex{transport.commutative.ring.spectra}.)

We fix the projective model structure on symmetric sequences with values in $\simpl \cA$ etc.
This is the model structure transferred from $\simpl \cA$.
Let $\tilde R_1 \in \cCh \cA$ and $R_1 \in \cA$ be any objects.
We regard $R_1$ as simplicially constant object in $\simpl \cA$.
We write $R$ for the commutative monoid in $\Sigma \cA \subset \Sigma \simpl \cA$ whose $n$th level is given by $R^{\t n}$ and likewise for $\tilde R \in \Sigma \cCh \cA$.
We suppose there is a weak equivalence
$$\varphi\colon \tilde R \r N (R)$$
in $\Sigma \simpl \cA$.
The categories of $R$- and $\tilde R$-modules are equipped with their stable positive model structure (\refth{stable.R.general}).
To simplify the notation, we will write $\Modsi_R := \Mod_R^{\stab, \posi} (\Sigma \simpl \cA)$, $\ModcCh_{\tilde R} := \Mod_{\tilde R}^{\stab, \posi}(\Sigma \cCh \cA)$, and similarly with $\ModCh_{\tilde R}$.

The normalization functor $N$ in the Dold--Kan equivalence \refeq{Dold.Kan}, applied to $\Sigma \cA$ instead of $\cA$, is symmetric lax monoidal and (nonsymmetric) oplax monoidal by means of the Alexander--Whitney and Eilenberg--Zilber maps (see for example \cite[\S5.4]{AguiarMahajan:Monoidal}).
Therefore the right adjoint $\Gamma$ is symmetric oplax monoidal and (nonsymmetric) lax monoidal.
However, the lax monoidal structural map
$$\Gamma(A)^{\t n} \t \Gamma(B) \r \Gamma(A^{\t n} \t B)$$
is a $\Sigma_n$-equivariant isomorphism for any $B \in \cCh \Sigma \cA$ provided that $A$ is a chain complex concentrated in degree $0$.
This can be checked using the explicit description of this map.
Dually, there is a lax monoidal map for $N$,
$$N(A \t B) \r \N(A) \t \N(B),$$
which is an isomorphism if $A$ is a constant simplicial object.
Applying this to $A = R$, we obtain an adjunction
$$N : \sMod_{R} \rightleftarrows \ModcCh_{N(R)} : \Gamma.\eqlabel{Dold.Kan.Mod.R}$$

With these preparations, we can now state the construction of commutative ring spectra.
Similar methods have been employed by Shipley to construct (noncommutative) ring spectra \cite[Theorem~1.1]{Shipley:HZ-algebras}.

\theo \label{theo--construction.spectra}
With the notation and assumptions fixed above, there is a functor
$$\H\colon \Alg_\Comm(\ModCh_{\tilde R}) \r \Alg_\Comm(\Modsi_R)$$
defined by
$$\H(A) := \Comm \circ^\ldf_\cO (\Gamma (R \t^\ldf_{\tilde R} \rdf \tau A)).$$
The spectrum $\H(A)$ represents the same cohomology as $A$ in the sense that the following derived mapping spaces are weakly equivalent, where $X$ is any object in~$\cA$:
$$\hMap_{\Modsi_R}(R \t X, \H(A)) \sim \Gamma \rdf \tau \hMap_{\ModCh_{\tilde R}}(\iota (\tilde R \t N(X)), A).$$
Moreover, the multiplicative structure is preserved in the strongest possible sense: if $X \in \cA (\subset \simpl \cA)$ is cofibrant and in addition a commutative coalgebra, there is a weak equivalence of convolution algebras
$$\hConv_{\Modsi_R}(R \t X, \H(A)) \sim \Comm \circ^\ldf_\cO \Gamma R \t^\ldf_{\tilde R} \rdf \tau \hConv_{\ModCh_{\tilde R}}(\iota (\tilde R \t N(X)), A).$$
\xtheo

\pf
We prove this using \refpr{monoidal}, a theorem of Richter \cite{Richter:Symmetry}, and the rectification theorem~\ref{theo--rectification.spectra}.

The functor $\iota$ is strong monoidal and $\tau$ is symmetric lax monoidal (because of the Leibniz rule).
Therefore, \refeq{iota.tau} induces a similar adjunction
$$\iota : \ModcCh_{\tilde R} \rightleftarrows \ModCh_{\tilde R} : \tau.$$
The unstable positive model structures on $\tilde R$-modules (\refde{Mod.R.unstable}) are transferred from \refeq{iota.tau}, which is a Quillen adjunction by assumption.
Therefore, by the universal property of the Bousfield localization, the stable positive model structures are also related by a Quillen adjunction.
Thus \refpr{monoidal}\refit{alg.adjunction} yields a weak equivalence
$$\hConv_{\ModCh_{\tilde R}}(\tilde R \t N(X), \rdf \tau A)\sim\rdf \tau \hConv_{\ModCh_{\tilde R}}(\iota (\tilde R \t N(X)), A).\eqlabel{we.1}$$

The map $\varphi\colon \tilde R \r N(R)$ induces a Quillen adjunction
$${-} \t_{\tilde R} N(R) : \ModcCh_{\tilde R} \rightleftarrows \ModcCh_{N(R)} : \text{restriction}.$$
The left adjoint is strong monoidal, the right adjoint is symmetric lax monoidal.
Since $\varphi$ is a weak equivalence by assumption, both this adjunction,
as well as the induced adjunction of commutative algebra objects are Quillen equivalences (\refex{transport.commutative.ring.spectra}, using the cofibrancy of the unit in $\cCh(\cA)$).
\refpr{monoidal}\refit{alg.3} gives a weak equivalence
$$\hConv_{\Mod_R(\cCh \Sigma \cA)}(N(R) \t X, R \t^\ldf_{\tilde R} \rdf \tau A)\sim R \t^\ldf_{\tilde R} \hConv_{\ModCh_{\tilde R}}(\tilde R \t N(X), \rdf \tau A).\eqlabel{we.2}$$

The next step is the Dold--Kan equivalence.
\refeq{Dold.Kan} is a Quillen adjunction by assumption.
Therefore so is \refeq{Dold.Kan.Mod.R} (where both sides carry the stable positive model structures of \refth{stable.R.general}).
Let $\cO = \cO_\Gamma$ be the endomorphism operad of $\Gamma$.
Using \refpr{monoidal}, we get a weak equivalence
$$\hConv_{\Modsi_R}(R \t X, \Gamma (R \t^\ldf_{\tilde R} \rdf \tau A))\sim\Gamma \hConv_{\ModcCh_{N(R)}}(N(R) \t X, N(R) \t^\ldf_{\tilde R} \rdf \tau A).\eqlabel{we.3}$$

Given a commutative monoid object $Z \in \Modsi_R$, it is easy to check that there is an isomorphism of $\cO$-algebras,
$$\Conv_{\Modsi_R}(R \t X, U Z) \stackrel \cong \lr U \Conv_{\Modsi_R}(R \t X, Z).$$
Here $U$ denotes the forgetful functors from commutative to $\cO$-algebras, by means of the unique map of operads $\cO \r \Comm$.
This passes to a weak equivalence
$$\hConv_{\Modsi_R}(R \t X, \rdf U Z) \stackrel \sim \lr \rdf U \hConv_{\Modsi_R}(R \t X, Z).\eqlabel{hConv.forget}$$

Using that $R$ is simplicially constant and therefore $N(R)$ is concentrated in degree~$0$, we can rewrite the adjunction \refeq{Dold.Kan.Mod.R} as the Dold--Kan equivalence applied to the abelian category $\Mod_R (\Sigma \cA)$:
$$N : \simpl \Mod_R (\Sigma \cA) \rightleftarrows \cCh \Mod_{N(R)} (\Sigma \cA) : \Gamma.$$
According to Richter's theorem \cite[Theorem~4.1]{Richter:Symmetry}, $\cO \r \Comm$ is a levelwise weak equivalence for the Dold--Kan equivalence on the abelian category $\Ab$.
The proof of loc.~cit.\ readily generalizes to a general abelian category such as $\Mod_R (\Sigma \cA)$.
Thus, \refth{rectification.spectra} establishes a Quillen equivalence
$$\Comm \circ_\cO - : \Alg_\cO(\Modsi_R) \leftrightarrows \Alg_\Comm (\Modsi_R) : U.$$
This Quillen equivalence and \refeq{hConv.forget}, applied to $Z = \Comm \circ^\ldf_\cO Y$ gives the following chain of weak equivalences of convolution algebras, i.e., commutative algebras in $\Mod_R$:
$$\eqalignno{
\Comm \circ^\ldf_\cO \hConv (X, Y) & \stackrel \sim \r \Comm \circ^\ldf_\cO \hConv (X, \rdf U_{\cO \r \Comm} \Comm \circ^\ldf_\cO Y)\cr
& \stackrel \sim \r \Comm \circ^\ldf_\cO \rdf U_{\cO \r \Comm} \hConv (X, \Comm \circ^\ldf_\cO Y)\cr
& \stackrel \sim \r \hConv (X, \Comm \circ^\ldf_\cO Y).&\eqlabel{we.4}\cr}$$

Combining \refeq{we.1}, \refeq{we.2}, \refeq{we.3}, and \refeq{we.4}, we obtain the desired weak equivalence
$$\eqalign{\hConv_{\Modsi_R}(R \t X, \H(A)) & = \hConv_{\Modsi_R}(R \t X, \Comm \circ^\ldf_\cO (\Gamma (R \t^\ldf_{\tilde R} \rdf \tau A))) \cr
&\sim\Comm \circ^\ldf_\cO \Gamma R \t^\ldf_{\tilde R} \rdf \tau \hConv_{\ModCh_{\tilde R}}(\iota (\tilde R \t N(X)), A).\cr}$$
\xpf

\section{A commutative ring spectrum for Deligne cohomology}
\label{sect--Deligne}

In this section, we construct a strictly commutative ring spectrum representing Deligne cohomology with integral coefficients.
For a smooth projective variety $X / \CC$, Deligne cohomology is defined as the hypercohomology group
$$\HD^n(X, \Z(p)) := \HH^n(X^\an, \underbrace{[\Z(p)\to\Omega^0_X\to\Omega^1_X\to\cdots\to\Omega^{p-1}_X]}_{=:\Z(p)_D}),\eqlabel{Deligne}$$
where $X^\an$ is the smooth complex manifold associated to $X$, $\Z(p) := (2\pi i)^p \Z$ sits in degree $0$ and $\Omega^*_X$ is the complex of holomorphic forms on $X^\an$.
Applications of Deligne cohomology range from arithmetic geometry, notably Beilinson's conjecture on special values of $L$-functions \cite{Beilinson:Higher} to higher Chern--Simons theory \cite[\S5.5.8]{Schreiber:Differential}.
The product structure on Deligne cohomology is surprisingly subtle.
It was defined by Beilinson by certain maps
$$-\cup_\alpha - \colon \Z(p)_D \t \Z(q)_D \r \Z(p+q)_D$$
that depend on a parameter $\alpha \in \CC$ \cite[Definition~3.2]{EsnaultViehweg:Deligne}.
This parameter is used to show that the product on the complexes of sheaves is commutative and associative up to homotopy.
In particular, $\bigoplus_{n,p} \HD^n(X, \Z(p))$ is a commutative ring.
This was used by Holmstrom and the second author to construct an unstructured commutative ring spectrum representing Deligne cohomology~\cite{HolmstromScholbach:Arakelov}.
This is the weakest possible requirement on the product operation on a spectrum: the multiplication is only commutative and associative up to homotopy.
In a somewhat similar vein, Hopkins and Quick studied ring spectra that result from replacing the Betti cohomology part in Deligne cohomology by a different ring spectrum, such as complex cobordism \cite{HopkinsQuick:Hodge}.
In this section, we provide a structured strictly commutative model for Deligne cohomology, which is the strongest possible multiplicative structure on such a spectrum.

We emphasize that we are working with integral coefficients.
For rational coefficients (i.e., with $\Q(p)$ instead of $\Z(p)$), all model categories in sight are freely powered, so it is possible to use Lurie's rectification result \cite[Theorem~4.5.4.7]{Lurie:HA} to obtain a strictly commutative ring spectrum.
However, integral coefficients are interesting from many points of view.
To refine the treatment of special $L$-values, which is up to rational factors in \cite{Scholbach:SpecialL},
it will be necessary to have the integral structure available.
One motivation for Hopkins and Quick's work is to find new torsion algebraic cycles, which also requires integral coefficients.
In yet another direction, one may speculate about the relation of modules over the Deligne cohomology spectrum and mixed Hodge modules by Saito \cite{Saito:Mixed}.
Again, for such considerations, it would be unnatural to throw away torsion.

Before discussing Deligne cohomology proper, we show how to turn a certain product structure on a fiber product of commutative differential graded algebras (cdga's) into a strictly commutative and associative one.
As in \refsect{construction}, our complexes are regarded as chain complexes, i.e., $\deg d = -1$.
Consider a diagram of cdga's, where we suppose that $B$ takes values in $\Q$-vector spaces:
$$A \stackrel a \r B \stackrel c \leftarrow C.$$
Because of rational coefficients, a path object for $B$ is given by $B \t Q$ \cite[Lemma~1.19]{Behrend:Differential}, where $Q$ is the chain complex of polynomial differential forms on $\Delta^1$ familiar from rational homotopy theory.
It is the complex in the left column, where the terms are in degrees $0$ and $-1$, respectively.
The complex~$R$ at the right is quasiisomorphic to~$Q$:
$$\vcd{5em}{
\Q[t]&\mapright{\ev(0),\,\ev(1)}&\Q\oplus \Q\cr
\mapdown{d}&&\mapdown{(a,b)\mapsto a-b}\cr
\Q[t]dt&\mapright{\int_0^1-}&\Q.\cr}$$
We endow $R$ with the multiplication $R \t R \r R$ given by the following matrix in terms of the standard basis $e_1, e_2 \in R_0$, $f \in R_{-1}$:
$e_1 \cdot e_1 = e_1$, $e_2 \cdot e_2 = e_2$, $f \cdot e_2 = f$, $e_1 \cdot f = f$, and all other products of basis vectors are $0$.
This product is associative and left unital, but not commutative.
Because of the latter defect, we consider the following diagram of associative left unital differential graded algebras
$$Q = Q \t \Q[0] \stackrel{\id\t1}\lr S := Q \t R \stackrel{(1,1) \t \id} \longleftarrow \Q[0] \t R = R.$$
The horizontal maps are induced by the unit elements of $Q$ and $R$, respectively.
These maps are quasiisomorphisms.
In addition, the augmentation maps $Q_0 = \Q[t] \r \Q^2$ and similarly for $S$ and $R$ commute with these quasiisomorphisms.
Therefore, there is a zigzag of quasiisomorphisms of associative (noncommutative, except for $D$) left unital differential graded algebras
$$D := A \x_B (B \t Q) \x_B C \stackrel \sim \r A \x_B (B \t S) \x_B C \stackrel \sim \gets A \x_B (B \t R) \x_B C.\eqlabel{associative.algebras}$$
The right object is just $E := \cone(A\oplus C \stackrel{a-c}\r B)[1]$ or, equivalently, the homotopy pullback $A \x_C^h B$ in the model category of chain complexes, while $D$ is the homotopy pullback in the (much more natural) model category of cdga's.
The quasiisomorphisms in \refeq{associative.algebras} are compatible with the respective product structures.
In particular, the induced product on $\H^*(D)$ agrees with the one on $\H^*(E)$.
Moreover, higher order multiplications, such as Massey products, also agree.

In the sequel, we just write $\PSh_\bullet := \PSh(\Sm / \CC, \Set_\bullet)$ for the presheaves of pointed sets on the site of smooth schemes over $\CC$.
We write $\simpl\PSh\Ab$ for simplicial presheaves of abelian groups and $\Ch \PSh$ for chain complexes of presheaves of abelian groups and likewise $\cCh \PSh$ for presheaves of chain complexes in degrees $\ge0$.
We equip the categories $\simpl \PSh_\bullet$, $\simpl \PSh \Ab$, $\cCh \PSh$, $\Ch \PSh$
with the local projective model structure.
They are the left Bousfield localizations of the projective model structures with respect to the covers
$$[(U \x_X V)_+ \rightrightarrows U_+ \sqcup V_+] \r X_+$$
(for $\simpl \PSh_\bullet$, and likewise for the three other categories).
Here $U \sqcup V \r X$ is a covering in the Zariski topology.

The corresponding categories of presheaves on the site $\SmAn$ of smooth complex manifolds are endowed with the local projective model structures with respect to the usual topology on $\SmAn$.

\rema
The results of this section hold unchanged if we replace the Zariski by the Nisnevich or etale topology on $\Sm / \CC$.
We could also furthermore localize with respect to $\A^1$ on the algebraic side and with respect to the disk $D^1$ on the analytic side.
\xrema

\lemm \label{lemm--four.cats}
There is a chain of Quillen adjunctions of the model categories mentioned above
$$\simpl \PSh_\bullet \mathrel{\mathop{\rightleftarrows}^{\Z[-]}}
\simpl \PSh \Ab \mathrel{\mathop{\rightleftarrows}^N_\Gamma}
\cCh \PSh \mathrel{\mathop{\rightleftarrows}^\iota_\tau}
\Ch \PSh.$$
The analogous categories for the site $\SmAn$ are related to these categories by Quillen adjunctions, for example
$$\Ch \PSh \mathrel{\mathop{\rightleftarrows}^{\an^*}_{\an_*}} \Ch \PSh(\SmAn, \Ab).$$
All these model categories are nice (\refde{nice.model.category}).
Moreover, their monoidal units are cofibrant.
\xlemm

\pf
The Quillen adjunctions of these categories, equipped with the projective model structure, transfer from the standard Quillen adjunctions for simplicial sets etc.
It passes to adjunctions of the local structures by the universal property of the Bousfield localization.
The Quillen adjunctions to presheaves on $\SmAn$ hold since $\an\colon \Sm \r \SmAn$ sends Zariski covers to analytic covers.
The properties required in \refde{nice.model.category} are discussed in \csy{\refsect{simplicial.presheaves}}.
Like any representable presheaf, the monoidal units, which are the representable presheaves associated to $\Spec \CC$ (or $\an(\Spec \CC)$) are cofibrant.
\xpf

We now turn toward the construction of our Deligne cohomology spectrum.
The cdga corresponding to Betti cohomology is defined by $A = \bigoplus_{p \in \Z} \rdf \an_* \Z(p)[-2p]$.
Similarly, let $B := \bigoplus_p \rdf \an_* \Omega^*[-2p]$, where $\Omega^*$ denotes the cdga of holomorphic differential forms.
Finally, let
$$C\colon X \mapsto \bigoplus_p \left (\colim_{\ol X} F^p \Omega^*_{\ol X}(\log (\ol X \backslash X)) \right )[-2p],$$
be the Hodge filtration, i.e., the stupid truncation $\sigma^{\ge p}$ of the complex of meromorphic forms on $\ol X^\an$, which are holomorphic on $X^\an$ and have at most logarithmic poles at $\ol X \backslash X$.
The colimit runs over all smooth compactifications $j\colon X \r \ol X$ such that $\ol X \backslash X$ is a strict normal crossings divisor.

We have obvious maps $A \stackrel a \r B \stackrel c \leftarrow C$ of cdga's of presheaves on $\Sm/\CC$ and consider the cdga $D$ and the weakly equivalent
dga $E = \cone (A \oplus C \stackrel{a-c}\lr B)[1]$ defined above.
On the other hand, we have the associative (but noncommutative) product on~$E$,
which is the particular case $\alpha=0$ of the classical product on the Deligne complexes \cite[Definition~3.2]{EsnaultViehweg:Deligne}.
The following result, which was already pointed out by Beilinson \cite[Remark~1.2.6]{Beilinson:Higher}, relates the two products:

\prop \label{prop--Deligne.complexes}
The cdga $D$ of presheaves on $\Sm / \CC$ defined above represents Deligne cohomology with integral coefficients in the sense that there is a functorial isomorphism, for any $X \in \Sm / \CC$:
$$\bigoplus_{p \in \Z} \HD^{2p-n}(X, \Z(p)) = \H_n\hMap_{\Ch\PSh}(X,D).$$
Under this isomorphism, the product on the right term induced by the multiplication on $D$ and the comultiplication given by the diagonal map $X \r X \x X$ agrees with the classical product on Deligne cohomology.
Moreover, all classical higher order products induced by the multiplication on $E$, such as Massey products \cite{Deninger:Higher}, agree with the corresponding higher order products on the cdga $D$, in the sense that the derived convolution algebras are weakly equivalent differential graded algebras of presheaves:
$$\hConv_{\Ch \PSh}(X, D) \sim \hConv_{\Ch \PSh}(X, E).$$
\xprop
\pf
The identification of Deligne cohomology with the right side is well known, see for example \cite[Lemma~3.2]{HolmstromScholbach:Arakelov} for a very similar statement.
Note that $X \in \Ch \PSh$ (i.e., the free abelian representable presheaf $\Z[X]$) is cofibrant in the projective model structure.
Hence the derived convolution algebras are defined (\refle{convolution.algebras}).
The extra information concerning the products follows immediately from the above discussion.
\xpf

In order to connect the cdga $D$ of \refpr{Deligne.complexes} to, say, algebraic $K$-theory, it is necessary to work with presheaves of simplicial sets.
As is well known to the experts (we learned it from Denis-Charles Cisinski), it is not possible to construct a strictly commutative simplicial abelian group representing Deligne cohomology or even Betti cohomology with integral coefficients.
In fact, Steenrod operations preclude the existence of a strictly commutative simplicial abelian (pre)sheaf representing Betti cohomology with integral coefficients.
This problem gives rise to an application of the operadic rectification for which we need to work in some category of symmetric spectra.
Because of its interest from the viewpoint of motivic homotopy theory, we work in the category of symmetric $\P^1$-spectra.

The category of motivic symmetric $\P^1$-spectra is the categories of modules over the monoid $R \in \Sigma \simpl \PSh_\bullet$
whose $n$th level is $R_n = (\P^1,1)^{\t n}$, i.e., the $n$th smash power of $\P^1$, pointed by~$1$.
Here and below we identify any scheme over~$\CC$ with its representable presheaf.
Note that $R_n$ is a constant simplicial presheaf.
We will abbreviate $\Modsi_{\P^1} := \Mod_R(\Sigma \simpl \PSh_\bullet)$.
We have a similar category $\ModsAb_{\P^1} := \Mod_{\Z[R]}(\Sigma \simpl \PSh \Ab)$ of modules over the monoid $\Z[R]$ whose $n$th level is $(\coker (\Z \stackrel 1 \r \Z[\P^1]))^{\t n}$.

Given the cdga $D = \bigoplus_p D_p$ of \refpr{Deligne.complexes}, we consider the symmetric sequence, again denoted by~$D$, whose $l$th level is given by $D(l) := \bigoplus_{p} D_{p+l}$, with a trivial $\Sigma_l$-action.
Then $D$ is a commutative monoid object in $\Sigma \Ch \PSh$.
Turning $D$ into a commutative monoid object in $\Modsi_{\P^1}$, i.e., a commutative symmetric $\P^1$-spectrum is equivalent to specifying a monoid map $R \r D$ in $\Sigma \simpl \PSh_\bullet$,
which is equivalent to specifying a pointed map $(\P^1,1) \r D(1) = \bigoplus_p D_{p+1}$ in $\simpl \PSh$ or, equivalently, a section on $\P^1$ of $D$ whose restriction to the point $1 \in \P^1$ vanishes.
In other words, we need to specify a line bundle with a flat connection on $\P^1$.
As is well known, a nontrivial line bundle (more precisely, a generator of $\HD^2(\P^1, \Z(1))$) is not representable by a global section, but has to be constructed by patching local data.
In the parlance of model categories, the nonfibrancy of~$D$ obstructs the existence of the required map.
We therefore replace~$\P^1$ by a weakly equivalent model.
This amounts to the standard idea of representing cohomology classes by \v Cech covers.
Consider the object $\widetilde {\P^1} \in \simpl \PSh$ defined as
$$\widetilde{\P^1} := [G \rightrightarrows \P^1 \backslash 0 \sqcup \P^1 \backslash \infty],$$
where the simplicial presheaf $G$ is defined by the homotopy pullback diagram
$$\xymatrix{
G \ar[r] \ar[d]^\sim & \rdf \an_* \mathcal U \ar[d]^\sim \cr
\Gm \ar[r] & \rdf \an_* \an^* \Gm
}$$
where $\mathcal U = [U^\pm \rightrightarrows U^+ \sqcup U^-]$ is the simplicial scheme whose only nondegenerate simplices are in degrees $1$~and~$0$, which is the \v Cech cover of $\Gm^\an$
arising from the cover
$\Gm^\an = U^+ \cup U^-$, where $U^+ = \{ z \in \CC, + z \notin \RR^{\ge0}\}$ and similarly with $U^-$.
The map $\mathcal U \r \an^* \Gm = \Gm^\an$ is a weak equivalence in the local model structure (with respect to the usual topology on $\SmAn$).
Hence the map $G \r \Gm$ is a weak equivalence.
Likewise, $[\Gm \rightrightarrows \P^1 \backslash 0 \sqcup \P^1 \backslash \infty] \r \P^1$ is a weak equivalence in the (Zariski) local model structure.
Therefore, the composition of these maps yields a weak equivalence $\widetilde{\P^1} \stackrel \sim \lr \P^1$.
It induces a Zariski-local weak equivalence
$$(\widetilde {\P^1})^{\t n} \stackrel \sim \r (\P^1)^{\t n}$$
in $\simpl \PSh$: for this we may use the local injective model structure on $\simpl\PSh_\bullet$
in which all objects are cofibrant, so weak equivalences are stable under tensor products.
Note that the local injective model structure is a monoidal localization of the injective model structure, which is monoidal, since $\sSet$ is monoidal.

Let $\widetilde{\Z \twi 1} := \coker (\Z \stackrel{1} \r N\Z[\widetilde{\P^1}]) \in \cCh \PSh$.
As in \refex{monoids}, we consider the free commutative monoid $\widetilde R$ on $\widetilde{\Z \twi 1}$, i.e., $\widetilde{R}_n = (\widetilde{\Z \twi 1})^{\t n}$.
By the above, the natural map
$$\tilde R \r N \Z[R]$$
is a Zariski-local weak equivalence in $\cCh \PSh$.

We now specify the $\tilde R$-module structure on the symmetric sequence $D$ defined above.
As for Betti cohomology, the map $N\Z[\widetilde{\P^1}] \r \an_* \Z(1)[-2] \subset A_1$ is given by the map $\an_* U^\pm = \an_* (\HH \sqcup \HH') \r \an_* \Z(1)$,
which is given by the section $2\pi i$ on $\HH := \{\Im z>0\}$ and $0$ on $\HH' := \{\Im z<0\}$.
The map $N\Z[\widetilde{\P^1}] \r C_1$ is determined by the map $\Gm \r F^1 \Omega^1_{\P^1}(\log (\{0, \infty\}))$
given by the section $d\log z = dz / z \in \Omega^1_{\P^1}(\log (\{0, \infty\}))$.

Finally, the map $\widetilde{\P^1} \r \Omega^* \t Q[-2] \subset B_1$ is given by the following map of complexes (the leftmost term lies in degree~2):
$$\xymatrix{
\an_* U^\pm \ar[r] \ar[d]_{{2 \pi i |_{\HH} \choose 0 |_{\HH'} } \t (1-t)} &
\Gm \sqcup \an_* (U^+ \sqcup U^-) \ar[r] \ar[d]^{dz/z \t t + {\log^+ z |_{U^+} \choose  \log^- z |_{U^-} }\t dt} &
\an_* \an^* \Gm \sqcup \P^1 \backslash 0 \sqcup \P^1 \backslash \infty \ar[d]^0 \cr
\Omega^0 \t \Q[t] \ar[r] &
\Omega^1 \t \Q[t] \oplus \Omega^0 \t \Q[t]dt \ar[r] &
(\Omega^* \t Q)^2 \ar[r] & \cdots
}$$
Here, $\log^+ z$ and $\log^- z$ are two branches of the complex logarithm (defined on $U^+$ and $U^-$, respectively) that agree on $\HH'$ and satisfy $\log^+ z - \log^- z = 2 \pi i$ on $\HH$.
One easily checks that this defines a map of complexes that yields a map
$$\widetilde{\Z\twi 1} \r D_1 = A_1 \x_{B_1} (B_1 \t Q) \x_{B_1} C_1.\eqlabel{bonding.map}$$
This defines an $\tilde R$-module structure on the symmetric sequence $D = (D(l))_{l \ge0}$ defined above.
Therefore, we obtain a strictly commutative motivic $\widetilde{\P^1}$ ring spectrum, which we denote by $\widetilde \HD$.
As above, we write $\ModcCh_{\widetilde {\P^1}} := \Mod_{\tilde R} (\Sigma \cCh \PSh)$ and likewise for $\ModCh_{\widetilde {\P^1}}$.
The cohomology represented by $\widetilde{\HD}$ is Deligne cohomology, including all higher product operations:

\prop \label{prop--penultimate}
The strictly commutative $\widetilde{\P^1}$ ring spectrum
$$\widetilde \HD \in \Alg_\Comm (\ModCh_{\widetilde {\P^1}}),$$
defined above is such that, for any smooth scheme $X / \CC$, there is a natural isomorphism of derived convolution algebras
$$\hConv_{\ModCh_{\widetilde {\P^1}}}(\tilde R \t X, \widetilde \HD) \sim \rdf U \hConv_{\Ch \PSh}(X, D),$$
where
$$U \colon \ModCh_{\widetilde {\P^1}} \r \Sigma \Ch \PSh \stackrel{\ev_0}\r \Ch \PSh$$ is the forgetful functor
and
$$U \colon \Alg_\Comm (\ModCh_{\widetilde {\P^1}}) \r \Alg_\Comm(\Ch \PSh)$$ is the induced functor (\refle{adjunction.Alg}).
In particular, by \refpr{Deligne.complexes}, all products and higher order operations such as Massey products are computed by $\widetilde \HD$.
\xprop
\pf
Again, $X \in \Ch (\PSh)$ is cofibrant, hence so is $\tilde R \t X$ as an $\tilde R$-module.
Therefore the derived convolution algebras are well defined.
By \refpr{monoidal}\refit{alg.adjunction}, we have to check that $\rdf U \widetilde \HD \r U \widetilde \HD = D$ is a weak equivalence.
This is implied by the fibrancy of $\widetilde \HD$, which by \refth{stable.R.general} follows from the fact that the maps
$$D(l) \r \hHom(\widetilde{\Z\twi 1}, D(l+1))$$
are weak equivalences.
This can be checked by applying the derived mapping space $\H_n \hMap(X, -)$ for any $X \in \Sm / \CC$ and any $n \in \Z$.
By \refpr{Deligne.complexes}, we get
$$\bigoplus_{p}\HD^{2p-n}(X, \Z(p)) \lr \ker \bigoplus_{p} \left (\HD^{2p+2-n}(\P^1 \x X, \Z(p+1)) \r \HD^{2p+2-n}(X, \Z(p+1)) \right ).$$
The map between them is the cup product with the element in $\zeta \in \HD^2(\P^1, \Z(1))$ represented by the map~\refeq{bonding.map}.
The element $\zeta$ generates this cohomology group, since the forgetful map to Betti cohomology $\HD^2(\P^1, \Z(1)) \r \H^2(\P^1, \Z(1)) \cong \Z$ is an isomorphism that sends $\zeta$ to $1$.
By the projective bundle formula for Deligne cohomology \cite[Proposition~8.5]{EsnaultViehweg:Deligne} the map above is an isomorphism.
\xpf

Finally, we construct the strictly commutative Deligne cohomology spectrum:

\theo \label{theo--ultimate}
There is a strictly commutative $\P^1$-spectrum with values in simplicial presheaves on $\Sm / \CC$,
$$\HD \in \Alg_\Comm (\Modsi_{\P^1})$$
defined by
$$\HD := \Comm \circ^\ldf_\cO \Gamma R \t^\ldf_{\tilde R} \rdf \tau \widetilde \HD,$$
which represents Deligne cohomology with integral coefficients, i.e., for any smooth algebraic variety $X / \CC$, there is an isomorphism
$$\pi_n \hMap_{\Modsi_{\P^1}} (R \t X, \HD) = \bigoplus_{p\in\Z} \HD^{2p-n}(X, \Z(p)).$$
The multiplication on the left induced by the ring spectrum structure on $\HD$ agrees with the classical product on Deligne cohomology.
Moreover, the convolution algebras are related by the following natural weak equivalence:
$$\Comm \circ^\ldf_\cO \Gamma R \t^\ldf_{\tilde R} \rdf \tau \hConv_{\ModCh_{\widetilde {\P^1}} }(\tilde R \t X, \HD) \sim
\hConv_{\Modsi_{\P^1}}(R \t X, \HD).$$
In particular, all higher order products on Deligne cohomology, such as Massey products, are represented by the commutative ring spectrum $\HD$.
\xtheo
\pf
This follows from Propositions \ref{prop--Deligne.complexes}, \ref{prop--penultimate} and \refth{construction.spectra}, applied to the Grothen\-dieck abelian category $\cA = \PSh(\Sm/\CC, \Ab)$ and the model structures mentioned in \refle{four.cats}.
\xpf

\rema
In the context of complex-analytic smooth manifolds, a variant of the Deligne complexes above is given by replacing the Hodge filtration as defined above by $F^p \Omega^*_X$.
The resulting groups (called \emph{analytic Deligne cohomology}) are the ones defined in \refeq{Deligne} for all (including noncompact) manifolds.
The above technique of rectifying this spectrum works essentially the same way.
An even more basic case covered by the techniques above is a strictly commutative ring spectrum representing Betti cohomology with integral coefficients.
\xrema

\section{Appendix}
\label{sect--appendix}

This appendix shows how to adapt the main results of this paper to even more general basic model categories than nice model categories.
Suppose that $\C$ is a symmetric monoidal model category that is h-monoidal, flat (as in \refde{nice.model.category}),
but instead of tractable (i.e., locally presentable, cofibrantly generated, and quasitractable) and pretty small,
we assume $\C$ is cellular \cite[Definition~12.1.1]{Hirschhorn:Model}, cofibrantly generated, quasitractable \csy{\refde{model.category.basic}},
and strongly admissibly generated \csy{\refde{strongly.admissibly.generated}}.
We will use the projective model structure on $\Sigma \C$ (\refle{projective.admissible}).
We also require that $R_n\in\C$ is cofibrant for all $n\ge0$.

As is explained in \csy{\refsect{topological.spaces}}, topological spaces satisfy all the properties above (but fail to be pretty small or combinatorial).

With these definitions, all the results in this paper remain valid, except where we use combinatoriality to construct model structures (this only occurs in \refle{injective.admissible} and \refpr{minimal.admissible}).
In particular, all operads in topological spectra are admissible, and a rectification result for weakly equivalent operads in topological spectra holds true.

We now explain where the combinatoriality and pretty smallness have been used in this paper and show how to modify the proofs using the above set of assumptions.
The point of combinatoriality is to ensure the existence of the Bousfield localization $\Mod_R^{\stab, \ppos}$ in the proof of \refth{stable.R.general}.
We have the following replacement for this:

\lemm
Under the above assumptions, the model structure $\Mod_R^{\stab, \ppos}$ exists.
\xlemm

\pf
We first show that $\Mod_R^\ppos$ is cellular.

As $\Sigma_m$ is a finite group, an object $X \in \Sigma_m \C$ is small (or compact in the sense of \cite[Definition~10.8.1]{Hirschhorn:Model})
relative to a class~$I$ of maps if $\varphi(X) \in \C$ is small (or compact) relative to $\varphi(I)$.
The cofibrations in $\Sigma^\ppos \C$ are effective monomorphisms since they are injective cofibrations by \ref{defi--admissible}\refit{admissible.injective} and the evaluation functor $\prod \ev_m\colon \Sigma \C \r \prod_m \C$ creates all limits and colimits.
Therefore $\Sigma^\ppos \C$ is cellular.

To show that $\Mod_R^\ppos$ is cellular note that the forgetful functor $U$ sends any cofibration $f$ in $\Mod_R^\ppos$ to an injective cofibration in $\Sigma^\ppos \C$ \csy{\refpr{transfer.basic}\refit{G.preserves.cofibrations}}, since $R \t -$ sends cofibrations in $\Sigma^{\ppos} \C$ to injective cofibrations in $\Sigma^{\ppos} \C$ by \ref{defi--admissible}\refit{admissible.injective}--\refit{admissible.product.2} and the assumption that $R$ is levelwise injectively cofibrant.
By cellularity of $\Sigma^\ppos \C$, $U(f)$ is thus an effective monomorphism in $\Sigma^{\ppos} \C$.
Since $U$ preserves limits and colimits, $f$ is itself an effective monomorphism.
The remaining conditions on a cellular model category \cite[Definition 12.1.1.(1),~(2)]{Hirschhorn:Model} follow since generating (acyclic) cofibrations in $\Mod_R^{\ppos}$ are of the form $R\t$ a generating (acyclic) cofibration of $\Sigma^\ppos \C$.

Moreover $\Mod_R^\ppos$ is h-monoidal, therefore left proper.
Hence the left Bousfield localization $\Mod_R^{\stab, \ppos}$ exists by \cite[Theorem~4.1.1]{Hirschhorn:Model}.
\xpf

\rema
\label{rema--Top}
The cofibrancy assumption on $R$ was used in the proof above
to show that $U$ sends cofibrations $f$ in $\Mod_R^{\ppos}$ to effective monomorphisms in $\Sigma \C$.
For $\C = \Top$, this is true for arbitrary~$R$.
Indeed, $U(f)$ is an inclusion of topological spaces (in all spectral levels) and therefore an effective monomorphism:
inclusions are stable under $A \x -$, pushouts, transfinite compositions, and retracts.
\xrema

Recall from \csy{\refde{pretty.small}} that a model category $\C$ is pretty small, if there is another ``auxiliary'' cofibrantly generated model structure on the same category with the same weak equivalences, (possibly) fewer cofibrations, and whose generating cofibrations have compact (co)domains.
Pretty smallness of a model category implies that (acyclic) h-cofibrations and weak equivalences are stable under transfinite compositions.
This point is helpful in the proofs of admissibility of operads and their rectification \cop{Theorems \ref{theo--O.Alg}, \ref{theo--cow.Coll}}
because the homotopic properties of transfinite compositions are easy to handle then.
A closer look, however, reveals that one does not need compactness of the (co)domains of the (auxiliary) generating cofibrations with respect to \emph{all} maps, but just with respect to the weak saturation of the class $Y \t_{\Sigma_n} v^{\pp n}$, where $Y$ is arbitrary and $v$ is a finite family of cofibrations. This property is called \emph{strongly admissibly generated} \csy{\refde{strongly.admissibly.generated}}.
Arguing as in \csy{\refth{symmetric.weakly.saturated}\refit{admissibly.generated.weakly.sat}}, it is enough to check this for finite families of generating cofibrations.
Using this property instead of pretty smallness, the above admissibility and rectification results can still be applied.

\lemm
Under the above assumptions on $\C$, $\Mod_R^{\stab, \posi}$ is strongly admissibly generated.
\xlemm

\pf
To show that $\Sigma^\ppos \C$ (with the projective model structure) is strongly admissibly generated, we have to show
that the (co)domain $X = \Sigma_k \cdot X'$ of a generating cofibration in $\Sigma_k \C$ is compact with respect to
$$\cell(\ev_k(Y \t_{\Sigma_n} v^{\pp n})). \eqlabel{proof.93}$$
Note that $X'$ is cofibrant by quasitractability.
One reduces to the case $Y = G_u(Z)$ and to generating cofibrations $v_i=G_{t_i}(h_i)$, $h_i \in \cof_{\Sigma_{t_i} \C}$.
By adjunction, $X$ is compact relative to some subcategory $\D \subset \Sigma_k \C$ if $X' \in \C$ is compact relative to $\varphi(\D) \subset \C$.
By \ref{defi--admissible}\refit{admissible.injective}, $\varphi(X)$ is cofibrant in $\C$.
As $\C$ is strongly admissibly generated, $X'$ is indeed compact with respect to $\varphi$ applied to the class \refeq{proof.93}
(if all $t_i=0$, the expression reduces to $Z \t_{\Sigma_n} h^{\pp n}$ otherwise the expression is given by the right side of \refeq{proof.77}).

The property of being strongly admissibly generated is stable when passing from $\Sigma \C$ to $\Mod_R^\ppos$ since we only have to consider generating cofibrations of $\Mod_R^\ppos$. Finally, the property is trivially stable under Bousfield localizations.
\xpf

\raggedright\rightskip0em plus \maxdimen

\bibliographystyle{dp}
\def\ZM#1{\href{https://zbmath.org/?q=an:#1}{Zbl #1}}
\def\MRx#1 #2\relax{\href{https://mathscinet.ams.org/mathscinet-getitem?mr=#1}{MR#1}}
\def\MR#1{\MRx#1 \relax}
\bibliography{bib}

\newcommand{\etalchar}[1]{$^{#1}$}
\providecommand{\bysame}{\leavevmode\hbox to3em{\hrulefill}\thinspace}
\providecommand{\ZM}{\relax\ifhmode\unskip\space\fi Zbl }
\providecommand{\MR}{\relax\ifhmode\unskip\space\fi MR }
\providecommand{\arXiv}[1]{\relax\ifhmode\unskip\space\fi\href{http://arxiv.org/abs/#1}{arXiv:#1}}
\providecommand{\MRhref}[2]{%
  \href{http://www.ams.org/mathscinet-getitem?mr=#1}{#2}
}
\providecommand{\href}[2]{#2}
\begin{thebibliography}{BHK{\etalchar{+}}15}

\bibitem[AM10]{AguiarMahajan:Monoidal}
Marcelo Aguiar and Swapneel Mahajan, \emph{Monoidal functors, species and
  {H}opf algebras}, CRM Monograph Series, vol.~29, American Mathematical
  Society, Providence, RI, 2010. \MR{2724388}, \ZM{1209.18002}.
  \url{https://pi.math.cornell.edu/~maguiar/a.pdf}.

\bibitem[Bar10]{Barwick:Left}
Clark Barwick, \emph{On left and right model categories and left and right
  {B}ousfield localizations}, Homology Homotopy Appl. \textbf{12} (2010),
  no.~2, 245--320. \MR{2771591}, \ZM{1243.18025}, \arXiv{0708.2067v2}.
  \url{http://projecteuclid.org/euclid.hha/1296223884}.

\bibitem[BB13]{BataninBerger:Homotopy}
Michael Batanin and Clemens Berger, \emph{{Homotopy theory for algebras over
  polynomial monads}}, 2013. \arXiv{1305.0086v6}.

\bibitem[Beh02]{Behrend:Differential}
Kai Behrend, \emph{{Differential Graded Schemes I: Perfect Resolving
  Algebras}}, 2002. \arXiv{math/0212225v1}.

\bibitem[Be{\u\i }84]{Beilinson:Higher}
A.~A. Be{\u\i }linson, \emph{Higher regulators and values of {$L$}-functions},
  Current problems in mathematics, {V}ol. 24, Itogi Nauki i Tekhniki, vol.~30,
  Akad. Nauk SSSR, Vsesoyuz. Inst. Nauchn. i Tekhn. Inform., Moscow, 1984,
  181--238. \MR{760999}, \ZM{0588.14013}.

\bibitem[BHK{\etalchar{+}}15]{BayehHessKarpovaKedziorekRiehlShipley:Left}
Marzieh Bayeh, Kathryn Hess, Varvara Karpova, Magdalena Kȩdziorek, Emily
  Riehl, and Brooke Shipley, \emph{Left-induced model structures and diagram
  categories}, Women in topology: collaborations in homotopy theory, Contemp.
  Math., vol. 641, Amer. Math. Soc., Providence, RI, 2015, 49--81.
  \MR{3380069}, \ZM{1346.18023}, \arXiv{1401.3651v2}.

\bibitem[BM03]{BergerMoerdijk:Axiomatic}
Clemens Berger and Ieke Moerdijk, \emph{Axiomatic homotopy theory for operads},
  Comment. Math. Helv. \textbf{78} (2003), no.~4, 805--831. \MR{2016697},
  \ZM{1041.18011}, \arXiv{math/0206094v3}.

\bibitem[BM09]{BergerMoerdijk:Derived}
Clemens Berger and Ieke Moerdijk, \emph{On the derived category of an algebra
  over an operad}, Georgian Math. J. \textbf{16} (2009), no.~1, 13--28.
  \MR{2527612}, \ZM{1171.18005}, \arXiv{0801.2664v2}.

\bibitem[Bor94]{Borceux:2}
Francis Borceux, \emph{Handbook of categorical algebra. 2}, Encyclopedia of
  Mathematics and its Applications, vol.~51, Cambridge University Press,
  Cambridge, 1994. \MR{1313497}, \ZM{1143.18002}.
  \url{http://gen.lib.rus.ec/book/index.php?md5=25437F068C15E3852CE4583EEACC7367}.

\bibitem[CD09]{CisinskiDeglise:Triangulated}
Denis-Charles Cisinski and Frédéric Déglise, \emph{{Triangulated categories
  of mixed motives}}, 2009. \arXiv{0912.2110v3}.

\bibitem[Den95]{Deninger:Higher}
Christopher Deninger, \emph{Higher order operations in {D}eligne cohomology},
  Invent. Math. \textbf{120} (1995), no.~2, 289--315. \MR{1329043},
  \ZM{0847.55014}.

\bibitem[EM06]{ElmendorfMandell:Rings}
A.~D. Elmendorf and M.~A. Mandell, \emph{Rings, modules, and algebras in
  infinite loop space theory}, Adv. Math. \textbf{205} (2006), no.~1, 163--228.
  \MR{2254311}, \ZM{1117.19001}, \arXiv{math/0403403v1}.

\bibitem[EV88]{EsnaultViehweg:Deligne}
H\'el\`ene Esnault and Eckart Viehweg, \emph{Deligne-{B}e\u\i linson
  cohomology}, Be\u\i linson's conjectures on special values of
  {$L$}-functions, Perspect. Math., vol.~4, Academic Press, Boston, MA, 1988,
  43--91. \MR{944991}, \ZM{0656.14012}.
  \url{https://page.mi.fu-berlin.de/esnault/preprints/ec/deligne_beilinson.pdf}.

\bibitem[Fre09]{Fresse:Modules}
Benoit Fresse, \emph{Modules over operads and functors}, Lecture Notes in
  Mathematics, vol. 1967, Springer-Verlag, Berlin, 2009. \MR{2494775},
  \ZM{1178.18007}, \arXiv{0704.3090v4}.

\bibitem[GG16]{GorchinskiyGuletskii:Symmetric}
S.~Gorchinskiy and V.~Guletski{\u\i }, \emph{Symmetric powers in abstract
  homotopy categories}, Adv. Math. \textbf{292} (2016), 707--754. \MR{3464032},
  \arXiv{0907.0730v4}.

\bibitem[GG17]{GorchinskiyGuletskii:Positive}
S.~Gorchinskiy and V.~Guletski{\u{\i}}, \emph{Positive model structures for
  abstract symmetric spectra}, Applied Categorical Structures \textbf{26}
  (2017), no.~1, 29--46. \arXiv{1108.3509v3}.

\bibitem[GH]{GoerssHopkins:Structured}
Paul~G. Goerss and Michael~J. Hopkins, \emph{{Moduli problems for structured
  ring spectra (June~8, 2005)}}.
  \url{https://math.northwestern.edu/~pgoerss/spectra/obstruct.pdf}.

\bibitem[GH04]{GoerssHopkins:Commutative}
P.~G. Goerss and M.~J. Hopkins, \emph{Moduli spaces of commutative ring
  spectra}, Structured ring spectra, London Math. Soc. Lecture Note Ser., vol.
  315, Cambridge Univ. Press, Cambridge, 2004, 151--200. \MR{2125040},
  \ZM{1086.55006}. \url{https://math.northwestern.edu/~pgoerss/papers/sum.pdf}.

\bibitem[Har09]{Harper:Symmetric}
John~E. Harper, \emph{Homotopy theory of modules over operads in symmetric
  spectra}, Algebr. Geom. Topol. \textbf{9} (2009), no.~3, 1637--1680.
  \MR{2539191}, \ZM{1235.55004}, \arXiv{0801.0193v3}.

\bibitem[Har10]{Harper:Monoidal}
John~E. Harper, \emph{Homotopy theory of modules over operads and
  non-{$\Sigma$} operads in monoidal model categories}, J. Pure Appl. Algebra
  \textbf{214} (2010), no.~8, 1407--1434. \MR{2593672}, \ZM{1231.55011},
  \arXiv{0801.0191v2}.

\bibitem[Hin97]{Hinich:Homological}
Vladimir Hinich, \emph{Homological algebra of homotopy algebras}, Comm. Algebra
  \textbf{25} (1997), no.~10, 3291--3323. \MR{1465117}, \ZM{0894.18008},
  \arXiv{q-alg/9702015v1}.

\bibitem[Hir03]{Hirschhorn:Model}
Philip~S. Hirschhorn, \emph{Model categories and their localizations},
  Mathematical Surveys and Monographs, vol.~99, American Mathematical Society,
  Providence, RI, 2003. \MR{1944041}, \ZM{1017.55001}.
  \url{http://gen.lib.rus.ec/book/index.php?md5=23EF8741E363DBAE39D1E7BA13F134E6}.

\bibitem[Hor13]{Hornbostel:Preorientations}
Jens Hornbostel, \emph{Preorientations of the derived motivic multiplicative
  group}, Algebr. Geom. Topol. \textbf{13} (2013), no.~5, 2667--2712.
  \MR{3116300}, \ZM{1281.55009}, \arXiv{1005.4546v2}.

\bibitem[Hov99]{Hovey:Model}
Mark Hovey, \emph{Model categories}, Mathematical Surveys and Monographs,
  vol.~63, American Mathematical Society, Providence, RI, 1999. \MR{1650134},
  \ZM{0909.55001}.
  \url{http://gen.lib.rus.ec/book/index.php?md5=229756633EF1A320BED3109B1AC1AB52}.

\bibitem[Hov01]{Hovey:Spectra}
Mark Hovey, \emph{Spectra and symmetric spectra in general model categories},
  J. Pure Appl. Algebra \textbf{165} (2001), no.~1, 63--127. \MR{1860878},
  \ZM{1008.55006}, \arXiv{math/0004051v3}.

\bibitem[HQ15]{HopkinsQuick:Hodge}
Michael~J. Hopkins and Gereon Quick, \emph{Hodge filtered complex bordism}, J.
  Topol. \textbf{8} (2015), no.~1, 147--183. \MR{3335251}, \ZM{1349.32009},
  \arXiv{1212.2173v3}.

\bibitem[HS15]{HolmstromScholbach:Arakelov}
Andreas Holmstrom and Jakob Scholbach, \emph{{Arakelov motivic cohomology I}},
  Journal of Algebraic Geometry \textbf{24} (2015), no.~4, 719--754.
  \arXiv{1012.2523v4}.

\bibitem[HS16]{HessShipley:Waldhausen}
Kathryn Hess and Brooke Shipley, \emph{{Waldhausen K-theory of spaces via
  comodules}}, Advances in Mathematics \textbf{290} (2016), 1079--1137.
  \arXiv{1402.4719v2}.

\bibitem[HSS00]{HoveyShipleySmith:Symmetric}
Mark Hovey, Brooke Shipley, and Jeff Smith, \emph{Symmetric spectra}, J. Amer.
  Math. Soc. \textbf{13} (2000), no.~1, 149--208. \MR{1695653},
  \ZM{0931.55006}, \arXiv{math/9801077v2}.

\bibitem[Jar00]{Jardine:Motivic}
J.~F. Jardine, \emph{Motivic symmetric spectra}, Doc. Math. \textbf{5} (2000),
  445--553. \MR{1787949}, \ZM{0969.19004}.
  \url{https://www.emis.de/journals/DMJDMV/vol-05/15.html}.

\bibitem[Lur]{Lurie:HA}
Jacob Lurie, \emph{{Higher algebra (September 18, 2017)}}.
  \url{https://www.math.ias.edu/~lurie/papers/HA.pdf}.

\bibitem[Lur09]{Lurie:HTT}
Jacob Lurie, \emph{Higher topos theory}, Annals of Mathematics Studies, vol.
  170, Princeton University Press, Princeton, NJ, 2009. \MR{2522659},
  \ZM{1175.18001}, \arXiv{math/0608040v4}.
  \url{https://www.math.ias.edu/~lurie/papers/HTT.pdf}.

\bibitem[MMSS01]{MandellMaySchwedeShipley:Model}
M.~A. Mandell, J.~P. May, S.~Schwede, and B.~Shipley, \emph{Model categories of
  diagram spectra}, Proc. London Math. Soc. (3) \textbf{82} (2001), no.~2,
  441--512. \MR{1806878}, \ZM{1017.55004}.

\bibitem[Per16]{Pereira:Cofibrancy}
Lu{\'{\i}}s~Alexandre Pereira, \emph{Cofibrancy of operadic constructions in
  positive symmetric spectra}, Homology, Homotopy and Applications \textbf{18}
  (2016), no.~2, 133--168. \arXiv{1410.4816v2}.

\bibitem[PS18a]{PavlovScholbach:Operads}
Dmitri Pavlov and Jakob Scholbach, \emph{Admissibility and rectification of
  colored symmetric operads}, Journal of Topology \textbf{11} (2018), no.~3,
  559--601. \arXiv{1410.5675v3}.

\bibitem[PS18b]{PavlovScholbach:Symmetry}
Dmitri Pavlov and Jakob Scholbach, \emph{Homotopy theory of symmetric powers},
  Homology, Homotopy and Applications \textbf{20} (2018), no.~1, 359--397.
  \arXiv{1510.04969v3}.

\bibitem[Ric03]{Richter:Symmetry}
Birgit Richter, \emph{Symmetry properties of the {D}old-{K}an correspondence},
  Math. Proc. Cambridge Philos. Soc. \textbf{134} (2003), no.~1, 95--102.
  \MR{1937795}, \ZM{1033.18009}.

\bibitem[Sai91]{Saito:Mixed}
Morihiko Saito, \emph{Mixed {H}odge modules and applications}, Proceedings of
  the {I}nternational {C}ongress of {M}athematicians, {V}ol.\ {I}, {II}
  ({K}yoto, 1990), Math. Soc. Japan, Tokyo, 1991, 725--734. \MR{1159259},
  \ZM{0826.32029}.
  \url{https://www.mathunion.org/fileadmin/ICM/Proceedings/ICM1990.1/ICM1990.1.ocr.pdf}.

\bibitem[Sch97]{Schwede:Spectra}
Stefan Schwede, \emph{Spectra in model categories and applications to the
  algebraic cotangent complex}, J. Pure Appl. Algebra \textbf{120} (1997),
  no.~1, 77--104. \MR{1466099}, \ZM{0888.55010}.

\bibitem[Sch13]{Schreiber:Differential}
Urs Schreiber, \emph{{Differential cohomology in a cohesive infinity-topos}},
  2013. \arXiv{1310.7930v1}.

\bibitem[Sch17]{Scholbach:SpecialL}
Jakob Scholbach, \emph{Special {$L$}-values of geometric motives}, Asian
  Journal of Mathematics \textbf{21} (2017), no.~2, 225--264.
  \arXiv{1003.1215v3}.

\bibitem[Shi04]{Shipley:Convenient}
Brooke Shipley, \emph{A convenient model category for commutative ring
  spectra}, Homotopy theory: relations with algebraic geometry, group
  cohomology, and algebraic {$K$}-theory, Contemp. Math., vol. 346, Amer. Math.
  Soc., Providence, RI, 2004, 473--483. \MR{2066511}, \ZM{1063.55006}.

\bibitem[Shi07]{Shipley:HZ-algebras}
Brooke Shipley, \emph{{$H\Bbb Z$}-algebra spectra are differential graded
  algebras}, Amer. J. Math. \textbf{129} (2007), no.~2, 351--379. \MR{2306038},
  \arXiv{math/0209215v4}.

\bibitem[{Spi}01]{Spitzweck:Operads}
Markus {Spitzweck}, \emph{{Operads, algebras and modules in model categories
  and motives.}}, Bonn: Univ. Bonn. Mathematisch-Naturwissenschaftliche
  Fakult\"at (Dissertation), 2001. \ZM{1103.18300}.
  \url{https://d-nb.info/970107374/34}.

\bibitem[SS00]{SchwedeShipley:Algebras}
Stefan Schwede and Brooke~E. Shipley, \emph{Algebras and modules in monoidal
  model categories}, Proc. London Math. Soc. (3) \textbf{80} (2000), no.~2,
  491--511. \MR{1734325}, \ZM{1026.18004}, \arXiv{math/9801082v1}.

\bibitem[SS03a]{SchwedeShipley:Equivalences}
Stefan Schwede and Brooke Shipley, \emph{Equivalences of monoidal model
  categories}, Algebr. Geom. Topol. \textbf{3} (2003), no.~1, 287--334.
  \MR{1997322}, \ZM{1028.55013}, \arXiv{math/0209342v2}.

\bibitem[SS03b]{SchwedeShipley:Stable}
Stefan Schwede and Brooke Shipley, \emph{Stable model categories are categories
  of modules}, Topology \textbf{42} (2003), no.~1, 103--153. \MR{1928647},
  \ZM{1013.55005}, \arXiv{math/0108143v1}.

\bibitem[SS12]{SagaveSchlichtkrull:Diagram}
Steffen Sagave and Christian Schlichtkrull, \emph{Diagram spaces and symmetric
  spectra}, Adv. Math. \textbf{231} (2012), no.~3-4, 2116--2193. \MR{2964635},
  \ZM{1315.55007}, \arXiv{1103.2764v2}.

\bibitem[TV08]{ToenVezzosi:HomotopicalII}
Bertrand To{\"e}n and Gabriele Vezzosi, \emph{Homotopical algebraic geometry.
  {II}. {G}eometric stacks and applications}, Mem. Amer. Math. Soc.
  \textbf{193} (2008), no.~902, x+224. \MR{2394633}, \ZM{1145.14003},
  \arXiv{math/0404373v7}.

\bibitem[Whi14]{White:Monoidal}
David White, \emph{{Monoidal Bousfield Localizations and Algebras over
  Operads}}, 2014. \arXiv{1404.5197v1}.

\bibitem[Whi17]{White:Model}
David White, \emph{Model structures on commutative monoids in general model
  categories}, Journal of Pure and Applied Algebra \textbf{221} (2017), no.~12,
  3124--3168. \arXiv{1403.6759v2}.

\end{thebibliography}

\end{document}